\documentclass[sts,preprint]{imsart}

\RequirePackage{amsthm,amsmath,amsfonts,amssymb}
\RequirePackage[authoryear]{natbib}

\startlocaldefs
\RequirePackage[OT1]{fontenc}
\usepackage{amsthm,amsmath,natbib}
\usepackage[textwidth=25mm,textsize=small,colorinlistoftodos]{todonotes}
\RequirePackage[colorlinks,citecolor=blue,urlcolor=blue]{hyperref}

\pubyear{2020}
\volume{0}
\issue{0}
\firstpage{1}
\lastpage{48}

\startlocaldefs
\numberwithin{equation}{section}
\theoremstyle{plain}

\endlocaldefs

\usepackage{thmtools}

\usepackage{mathtools,amsmath,physics,framed,enumitem}
\usepackage{arydshln}
\usepackage{graphicx,float}
\usepackage{xcolor}
\usepackage{soul}
\usepackage{comment}
\usepackage{natbib}
\usepackage{url}
\usepackage{bm}
\newcommand{\distas}[1]{\mathbin{\overset{#1}{\kern\z@\sim}}}%

\newcommand{\diag}{\mbox{diag}}

\usepackage{hyperref}
\hypersetup{
    colorlinks=true,
    linkcolor=blue,
    filecolor=magenta,      
    urlcolor=red
}

\newsavebox{\mybox}\newsavebox{\mysim}
\usepackage{
amsmath, 
amsthm ,
amssymb, 
}
\theoremstyle{definition}
\usepackage{fixltx2e,fix-cm}
\usepackage{graphicx}
\usepackage{multicol}
\usepackage{bm}
\usepackage{mathptmx} 

 \usepackage{natbib}
\usepackage{showidx}
\usepackage{color}

\usepackage{changepage}
\usepackage{subfig}   

\usepackage[english]{babel}
\usepackage[T1]{fontenc}
\usepackage[ansinew]{inputenc}
\usepackage{makeidx}

\makeindex 

\definecolor{gray}{rgb}{0.95,0.95,0.95}
\definecolor{brightskyblue}{rgb}{0.9, 0.95, 0.99}

\newcommand{\DHo}{dimensional homogeneity}

\newcommand{\Real}{{\rm I\!R}}

\usepackage[normalem]{ulem}   

\renewcommand{\Xi}{\ensuremath{X_i}}

\newcommand{\bfx}{\ensuremath{\mathbf{x}}}

\newcommand{\bfX}{\ensuremath{\mathbf{X}}}

\newcommand{\Sx}{S_{1\textbf{x}_1}}

\newcommand{\Sxinv}{S_{1\textbf{x}_1}^{-1}}
\newcommand{\Syinv}{S_{1\textbf{y}_1}^{-1}}
\newcommand{\Qx}{Q_{1\textbf{x}_1}}
\newcommand{\Qy}{Q_{1\textbf{y}_1}}

\newcommand{\bftheta}{\ensuremath{\mathbf{\theta}}}

\renewcommand{\strut}{\rule[-2mm]{0cm}{6mm}}

\def\be{\begin{equation}}
\def\ee{\end{equation}}
\def\bi{\begin{itemize}}
\def\ei{\end{itemize}}

\def\bdesc{\begin{description}}
\def\edesc{\end{description}}
\def\benum{\begin{enumerate}}
\def\eenum{\end{enumerate}}

\renewcommand{\Xi}{\ensuremath{X_i}}

\newcommand{\bfP}{\textbf{P}}
\newcommand{\bfQ}{\textbf{Q}}
\newcommand{\bfR}{\textbf{R}}
\newcommand{\bfS}{\textbf{S}}

\newtheorem{theorem}{Theorem}

\newtheorem{example}{Example}

\newtheorem{remarks}{Remark}

\renewcommand{\paragraph}[1]{\vspace{.5cm}{\normalfont\normalsize\bfseries\emph{#1}}}

\def\l({\left(}
\def\r){\right)}
\def\l{\left}

\def\bi{\begin{itemize}}
\def\ei{\end{itemize}}

\renewcommand{\l}{\left}
\renewcommand{\r}{\right}

\newcommand{\bdm}{\begin{displaymath}}
\newcommand{\edm}{\end{displaymath}}

\newcommand{\bI}{\mathbf{I}}

\newcommand{\bzero}{\mathbf{0}}

\newenvironment{myenv}{\begin{adjustwidth}{0.75cm}{0.25cm}}{\end{adjustwidth}}

\newcommand{\bexmp}{\vspace{0.25cm}\begin{myenv}\begin{exmp}}
\newcommand{\eexmp}{\end{exmp}\end{myenv}\vspace{0.25cm}}

\endlocaldefs

\begin{document}

\begin{frontmatter}
\title{Dimensional Analysis in Statistical Modelling}
\runtitle{Scientific and statistical modelling}

\begin{aug}
\author{\fnms{Tae Yoon} \snm{Lee} 
\ead[label=e1]{dlxodbs@student.ubc.ca}
}
\author{\fnms{James V} \snm{Zidek}\thanksref{t1}
\ead[label=e2]{jim@stat.ubc.ca}
\ead[label=u1,url]{www.stat.ubc.ca/~jim/}
}
\author{\fnms{Nancy} \snm{Heckman}
 \thanksref{t1}
\ead[label=e3]{nancy@stat.ubc.ca}
\ead[label=u2,url]{www.stat.ubc.ca/~nancy/}
}
\thankstext{t1}{The research reported in this paper was supported by  grants from the Natural Science and Engineering Research Council of Canada.}
\affiliation{University of British Columbia}

\address{Department of Statistics\\
University of British Columbia \\
2207 Main Mall\\
Vancouver, BC \\
Canada V6T 1Z4\\
\printead{e1,e2,e3}.}
\runauthor{Lee et al.}

\end{aug}

\begin{abstract}
Building on recent work in statistical science, the paper presents a theory for
 modelling natural phenomena that unifies physical and statistical 
paradigms based on the underlying principle that a model must be 
non-dimensionalizable.  After all, such phenomena cannot depend on 
how the experimenter chooses to assess them.  Yet the model itself 
must be comprised of quantities that can be determined  theoretically or 
empirically. Hence, the underlying principle requires that the model 
represents these natural processes correctly no matter what scales and units of measurement are  selected.  This goal was realized for physical modelling
through the celebrated theories of Buckingham and Bridgman and for 
statistical modellers through the invariance principle of Hunt and Stein. Building 
on recent research in statistical science, the paper shows how the latter 
can embrace and extend the former.  The invariance principle is extended 
to encompass the Bayesian paradigm, thereby enabling an assessment of model uncertainty. 
The paper covers topics not ordinarily seen in statistical science regarding 
dimensions, scales, and units of quantities in statistical modelling. 
It shows the special difficulties that can arise when models involve
 transcendental functions, such as the logarithm which is used e.g. in likelihood analysis and is a singularity in the family of Box-Cox family of transformations. Further, it demonstrates the 
importance of the scale of measurement, in particular how differently modellers
must handle ratio- and interval-scales. 
\end{abstract}



\end{frontmatter}


\section{Introduction}\label{sec:introduction}
The  importance of dimension, scale and units of measurement in modelling is largely ignored in statistical modelling.  In fact, an anonymous reviewer stated:
\begin{quotation}
	``Generally speaking, statisticians treat data as already dimensionless by taking the numeric part, which is equivalent to dividing them by their own units...''
\end{quotation}
\noindent
Others have long recognized the role   of dimensions, and hence of their scales and units of measurement, in modelling; dimensions can be used to  simplify  model fitting by reducing the number of quantities involved to a few non-dimensionalized ones.  A principal contribution of this paper makes clear to statisticians the importance of dimensions, scales and units in statistical modelling. We also develop a framework that incorporates these important ideas via a statistical invariance approach.  This allows us to extend the existing work's focus on physical quantities, which by their nature must lie on ratio-scales, to other scales and to vector spaces for multivariate responses. Yet another contribution is addressing a number of issues that are crucial in laying the foundation for the extension. These include: 
adopting  a sampling distribution supported on an interval-scale when the actual support for the sampling distribution is a subset of a ratio-scale; the meaninglessness of applying transcendental transformations such as the logarithm to quantities with units.

This paper, which is partly expository, describes and contributes to the unification of two overlapping paradigms for modelling natural phenomena.  For simplicity we will refer to these  as statistical and what \citet{meinsma2019dimensional} calls physical. Physical models are usually deterministic and developed for a specific phenomenon.

In this approach, model development cannot ignore the dimensions, scales and units of measurement on which empirical implementation and assessment will be based. Indeed, \citet{buckingham1914physically} believed that a valid model cannot depend on how an investigator chooses to measure the quantities that describe the phenomena of interest. After all, nature cannot know what the measuring methods of science are. Consequently,  Buckingham suggested that any valid model must be nondimensionalizable, leading him to his celebrated Pi-theorem.  In contrast, \citet{bridgman1931dimensions} argued that models must depend on the measurements but must be invariant under rescaling. Based on the latter premise he was able to derive the invariant \lq\lq$\pi$-functions\rq\rq~of the measurements that were central to Buckingham's theory. The work of these pioneers spurred the development of what is now known as dimensional analysis (DA), its notions of dimensional homogeneity (DH) and its quantity calculus, explored  in  depth in Section \ref{sec:dimension}.

The following example renders the ideas above in a more concrete form.
\begin{example}\label{ex:newton}
 Newton's second law of motion exemplifies the physical approach to modelling: 
\begin{equation}\label{eq:newton}
		 a  = F~ M^{-1}.
\end{equation}
Here $a$ denotes acceleration, the second derivative with respect to time of the location  of an object computed with respect to the starting point of its motion. $F$ and $M$ are, respectively, the force acting on the object and its mass.  The model in Equation (\ref{eq:newton})  satisfies the fundamental requirement of DH --  the units on the left hand side are the same as the units on the right hand side. Moreover, all three of the quantities in the model lie on a ratio-scale i.e. they are inherently positive, having a structural $0$ for an origin when and where the motion began.

The work of Buckingham and Bridgman cited above implies the quantities in a valid model have to be  transformable to dimensionless alternatives called $\pi$-functions. In the case of Newton's law, we can use $M$ and $F$  to transform $a$ into a dimensionless quantity to get the simpler but mathematically equivalent model involving a single $\pi$-function:
\begin{equation}\label{eq:newtonpi}
\pi(a, F,  M) \equiv a F^{-1} M=1.
\end{equation}
\end{example}


In contrast to physical modelling that commonly takes a bottom-up approach,  that of statistics as a discipline was top-down \citep{magnello2009karl} when, in the latter part of the nineteenth century, Karl Pearson established mathematical statistics with its focus on abstract classes of statistical models.   Pearson's view freed the statistician from dealing with the demanding contextual complexities of specific applications. In his abstract formulation, Pearson was able to:  incorporate model uncertainty expressed probabilistically; define desirable model properties;  determine conditions under which these occur; and develop practical tools to implement models that possess those properties for a specific application.  The emphasis on mathematical foundations led inevitably to an almost total disregard of the constraints brought on by dimension, scale and units of measurement.
Thus statisticians often use a symbol like $ X $ to mean a number to be manipulated in a formal analysis in equations, models and transformations.  On the other hand,  scientists use such a symbol to represent some specific aspect of a natural phenomenon or process. The scientist's goal is to  characterize $X$ through a combination of basic principles and empirical analysis. This  leads to the need to specify  one or more ``dimensions'' of $X$, e.g. length $L$. That in turn leads to the need to specify an appropriate  ``scale'' for $X$, e.g. categorical, ordinal, interval or ratio. For interval- and ratio-scales, $ X $ would have some associated units of measurement depending on the nature and resolution of the device making the measurement.  
How all of these parts of $ X $ fit together is the subject addressed in the realms of measurement theory and DA.

In recent years, importance of dimensions, scales, and units of measurements has been progressively recognized in statistics. At a deep theoretical level, \citet{hand1996statistics} considers different interpretations of measurement, studying what things can be measured and how numbers can be assigned to measurements. On a more practical side, inspired by applied scientists, \citet{finney1977dimensions} demonstrates how the principle of DH can be used to assess the validity of a statistical model. The first application of DA in statistical modeling appears in the work of \citet{vignaux1999theory}, who develop a framework for applying DA to linear regression models. 
The practicality of DA is illustrated to a great extent in
design of experiments by \citet{albrecht2013experimental}.
While much has been written in this area by nonstatisticians, such as \citet{luce1964generalization}, surprisingly little has been written by statisticians. 

A natural statistical approach to these ideas is via the statistical invariance principle due to Hunt and Stein in unpublished work \citep[Chapter 6]{lehmann2010testing}. Despite the abstraction of Pearson's approach they articulated an important principal of modelling- that when a test of a hypothesis about a natural phenomenon based on a sample of measurements rejects the null hypothesis, that decision should not change if the data were transformed to a different scale, e.g. from Celsius (interval scale) to Fahrenheit (ratio scale). That led them to the statistical invariance principle: methods and models must transform coherently under measurement scale transformations.

However, the link between DA and the statistical invariance does not seem to have been recognized until the work of Shen, Lin, and his co-investigators \citep{lin2013comment,shen2014dimensional,shen2015dimensional,shen2018conjugate,shen2019statistical}. They develop a theoretical framework for applying DA to statistical analysis and the design of experiments while showcasing the practical benefits through numerous examples involving physical quantities. 
In their framework, Buckingham's Pi-theorem plays the key role in unifying DA and the statistical invariance. In our paper, we extend their work in two ways: (1) elucidating the link between DA and the statistical invariance by removing the dependency on Buckingham's Pi-theorem and (2) in doing so, freeing ourselves from restricting to modelling physical quantities and ultimately embedding scientific modeling within a stochastic modeling framework in a general setting. 

This paper considers issues that arise when $ X $ lies on an interval-scale with values on the entire real line and when $ X $ lies on a ratio-scale, that is,  with non-negative values and a real origin $ 0 $ having a meaning of ``nothingness\rq\rq.  A good example to keep in mind is the possible scales of temperature; it can  be measured on the ratio-scale in units of degrees Kelvin ($^{\circ}K$), where $ 0^{\circ}K$  means all molecular momentum is lost. Alternatively, it can be measured on the interval-scale in units of degrees Celsius ($^{\circ}C$) where $0 ^{\circ}C$ means water freezes. The same probabilistic model cannot be used for both scales although often, in practice, the Gaussian model for the interval-scale is used as an approximation for the ratio-scale.
%
%
%

We conclude this Introduction with a summary of the paper's contents and findings. Section \ref{sec:uncstatistician} introduces us to the Unconscious Statistician through examples that  illustrate the importance of units of measurement. 
That demonstrated importance of units then leads us into Section \ref{sec:dimension}, which is largely a review of the basic elements of DA, a subject taught in the physical sciences but rarely in statistics.  
We describe a key tool, quantity calculus, 
which is the algebra of units of measurement and DH.
%
%

In Section  \ref{sec:new_transforming}, we discuss problems that might arise when statisticians transform variables.    
Sometimes the transformation leads to an unintended change of scale, e.g. when a Gaussian distribution on $ (-\infty, \infty) $ is adopted as an approximation to a distribution on $ (0, \infty) $. This can matter a lot when populations with responses on a ratio-scale are being compared. We discuss when such an approximation and hence transformation may be justified. 
Even using the famous family of Box-Cox transformations can cause problems, in particular with the limiting case logarithmic transformation.

Having investigated units and scales, we then turn to the models themselves. It turns out that when restricted by the need for DH, the class of models 
is also restricted; that topic is explored in Section \ref{sec:relationships} where we review the famous Buckingham Pi-theorem. We also see for the first time the ratio-scale's cousin, the interval-scale. All this points to the need for a general approach to scientific modelling that was foreseen in Hunt and Stein's unpublished work on the famous invariance principle.

Section \ref{sec:invariance} gets us to the major
contribution of the paper, namely extending the work of recent vintage by statistical scientists (Shen and his co-investigators) on the invariance principle as developed in its classical setting, the frequentist paradigm, and applied principally to variables on a 
ratio-scale. 
Our work constitutes a major extension of that scientific modelling paradigm. 

In particular, Section  \ref{sec:extendedinvariance} extends the invariance principle and moves our new approach to scientific modelling to the Bayesian paradigm.   This  major extension of both the scientific and statistical modelling approaches allows for quantities that could represent uncertain parameters, thereby embedding uncertainty quantification directly into the modeling paradigm. Indeed our new paradigm enables model uncertainty itself to be incorporated
.  

The paper wraps up with discussion in Section \ref{sec:discussion} and concluding remarks in Section \ref{sec:conclusions}. The supplementary material includes additional discussion, in particular a summary of the controversy of ancient vintage about whether or not taking the logarithm of measurements with units is valid, how Buckingham's theory leads us to the famous Reynolds number, a general theory for handling data on an interval-scale, and finally a brief review of statistical decision analysis for Section \ref{sec:extendedinvariance}.

\section{The unconscious statistician}\label{sec:uncstatistician} 

We start with critically examining key issues surrounding the topics of dimension and measurement scales through the Unconscious Statistician. We present three examples that illustrate some of the issues we'll be exploring.  

\vskip 10pt
\begin{example}\label{ex:poiss}
Ignoring scale and units of measurement when creating models can lead to difficulties; we cannot ignore the distinction between numbers and measurements. Consider  the Poisson random variable $ X$.  The claim is often made that the expected value and variance of $X$ are equal. But if  $ X $ has units, as it did when the distribution was first introduced in 1898 as the number of horse kick deaths in a year in the Prussian army \citep{hardle2017bortkiewicz}, then clearly, the expectation and variance will have different units and therefore cannot be equated.  
\end{example}
   
\begin{example}\label{ex:mle}
Consider a random variable representing length in millimetres, $Y$ is normally distributed with  mean $\mu$ and variance $\sigma^2$, independently measured $ n $ times to yield data $ y_1,\ldots, y_n $.  
Assume, as is common, that $\mu$ is so large that there is a negligible chance that any of the $y_i$'s are negative (we return to this common assumption in Section \ref{sec:new_transforming}).

Then the maximum likelihood estimate (MLE) of $\mu$ is easily shown to be the sample average $\bar{y}$ and the MLE of $\sigma^2$ is then maximizer of the profile likelihood
  \begin{equation}\label{eq:sigma2mle}
  	  L(\sigma^2) 
 =   (\sigma^2)^{-n/2}  ~ \exp \left\{-n \tilde{\sigma}^2/(2 \sigma^2) \right\}
  \end{equation}
where  $\tilde{\sigma}^2 = \sum_1^n(y_i-\bar{y})^2/n$, which has units of  mm$^2$. 
That estimate for $\sigma^2$ is the  MLE of $ \sigma^2$, which is easily shown by differentiating  $ L(\sigma^2) $ with respect to $ \sigma^2 $ and setting the result equal to zero. 
We note that, by any sensible definition of unit arithmetic,  $\tilde{\sigma}^2/ \sigma^2$ is unitless and so the units of $L(\sigma^2)$ are mm$^{-n}$.  

The  Unconscious Statistician simplifies the maximization of $L$ by maximizing instead, its logarithm, believing that this alternative approach yields the same result. 
.  
The statistician finds the log of $L$ to be
\[
	l(\sigma^2)  = \ln \left[   (\sigma^2)^{-n/2}  ~ \exp \left\{-n \tilde{\sigma}^2/(2 \sigma^2) \right\} \right]
=  - \frac{n}{2} \left[   \ln(\sigma^2)  + \tilde{\sigma}^2/\sigma^2 \right].
\]
Since the second term is unitless, dimensional homogeneity implies that the first term 
$\ln(\sigma^2)$ must also be unitless. 
So where did the units go?  Analyses in Subsection \ref{sec:logarithm} suggest the units reduce to a unitless $1$ by constructive processes for the logarithm that define it. 
The result is 
$\ln(\sigma^2) = \ln(\{\sigma^2 \})$, the curly brackets demarcating the numerical part of  $\sigma^2$, gotten by dropping the units of measurement.
But $\sigma^2$ itself has units mm$^2$ and it seems unsettling to have them disappear simply by taking the logarithm. 

However, the Unconscious Statistician ultimately gets the correct answer by failing to recognize 
the distinction between the scales of $\{\sigma^2\}$ and $\sigma^2$ 
So the derivative, which represents the relative rate of change between quantities on different scales, is computed as 
$d~\ln(\{\sigma^2 \})/d \sigma^2 = $ mm$^{-2} d~\ln(\{\sigma^2 \})/d \{\sigma^2\}$ rather than 
$d~\ln(\{\sigma^2 \})/d \{\sigma^2\}$. This then restores the missing units in the final result. As a fringe benefit, the second derivative $d^2 \ln(\{\sigma^2\}) / d^2 \sigma^2$, whose inverse 
defines Fisher's information, also turns out to have the appropriate units.  
 \end{example}
 
However, the story does not end there.
The problem of logarithms and their units warrants further discussion such as that in  Subsection {\ref{sec:logarithm}}. That discussion indicates that calculating the logarithm of the likelihood is, in general, not sensible.  
 
\begin{remarks}\label{rem:likelihood}
In the frequentist paradigm for statistical modelling, the likelihood is defined by the sampling distribution, which depends on the stopping rule employed in collecting the sample. The likelihood function then becomes an equivalence class. The likelihood ratio can then be used to specify a member of that class. In Example  \ref{ex:mle} a reference normal likelihood could be used with $\sigma^2$ set to  a substantively meaningful $\sigma_0^2$. The MLE of $\mu$ and $\sigma^2$ would then maximize this relative likelihood. This leads again to $\hat{\mu}=\bar{y}$, but now the MLE of $\sigma^2$ is found by maximizing the unitless $L(\sigma^2)/L(\sigma_0^2)$:
 \[
  \frac{L(\sigma^2)}{L(\sigma_0^2)} = 
   \left(     \frac{\sigma^2}{\sigma_0^2} \right)^{-n/2}
 \exp \left\{     -\frac{ n \tilde{\sigma}^2 }{ 2 \sigma_0^2}  \left[    \left( \frac{ \sigma_0^2}{\sigma^2}  \right)  - 1
  \right]
 \right\}.
 \]
 We can now maximize this ratio as a function of  the unitless $t = \sigma^2/\sigma_0^2$, by taking logarithms, differentiating with respect to $t$, setting it equal to 0 and solving for $\hat{t} = \tilde{\sigma}^2/\sigma_0^2$, and so finding that 
 $\hat{\sigma}^2 = 44.2$ mm$^2$. 
\end{remarks}
 \vskip 10pt
 
Two complimentary, unconscious choices in Example \ref{ex:mle} lead ultimately to a correct MLE. Things don't go so well for two unconscious statisticians seen in the next example.

 \vskip 10pt
\begin{example}\label{ex:unconsciousstatistician}
Here,  the data are assumed to follow the model that relates $Y_i$, a length, to $t_i$, a time:
\[
Y_i = 1 + \theta t_i + \epsilon_i, i = 1,~ \ldots, 2n.
\]
Here the $ \epsilon_i$'s are independent and identically distributed as  a $ N(0, \sigma^2) $ for a known $ \sigma $.  Suppose that $t_1 = \cdots = t_n = 1$ hour while  $ t_{n+1} = \cdots = t_{2n} = 2$ hours. Let $ \bar Y_1 = \sum_{i = 1}^n Y_i/n$, and $\bar Y_2 = \sum_{i = n+1}^{2n} Y_i/n$.  An analysis might go as follows when two statisticians A and B get involved.

First they both compute the likelihood and learn that the MLE is found by minimizing the sum of squared residuals $SSR(\theta)$:
\[
SSR(\theta) = \sum_{i = 1}^{2n} [Y_i - 1 - \theta t_i]^2
\]
which gives the MLE of $\theta$,
\[
\hat \theta = \frac{\sum_{i = 1}^{2n} t_i(Y_i-1)}{\sum_{i = 1}^{2n}t_i^2} = \frac{n\bar Y_1 + 2n\bar Y_2 - 3n}{5n} = \frac{\bar Y_1 + 2  \bar Y_2 - 3}{5}.
\]
Then for prediction at time $t=1$ hour, they get
\[
\widehat Y = 1 + \hat \theta \times 1 = 1 + \frac{\bar Y_1 + 2 \bar Y_2 - 3}{5}.
\]
Suppose that $\bar Y_1 = 1$ foot, or 12 inches, and $\bar Y_2 = 3$ feet, or 36 inches. 
Statistician A uses feet and predicts $Y$ at time $ t=1 $ hour  to be 
\[
\widehat Y_A = 1 + \frac{1 + 2 \times 3 - 3}{5} = 1.8~ {\rm{feet}} = 21.6~{\rm{inches}}.
\]
But Statistician B  uses inches and predicts $Y$ at $t=1 $  hour to be 
\[
\widehat Y_B = 1 + \frac{12 + 2 \times 36 - 3}{5} = 17.2~ {\rm{inches}}.
\]
\end{example}
\vskip .1in

 What has gone wrong here?   The problem is that the stated model implicitly depends on the units of measurement.  For instance, the numerical value of the expectation of $Y_i$ when $t_i=0$ is equal to 1, no matter what the units of $Y_i$.  When $t_i=0$, Statistician A expects $Y_i$ to equal 1 foot and Statistician B expects $Y_i$ to equal 1 inch. We can see that the problem arises because the equation defining the model does not satisfy DH, since the ``1'' is unitless. In technical terms, we would say that this model is not invariant under scalar transformations.  Invariance is important when defining a model that involves units.  However, one could simply avoid the whole problem of units in model formulation by constructing the relationship between $Y_i$ and $t_i$ so that there are no units.  This is exactly the goal of the Buckingham Pi-theorem, presented in Subsection \ref{subsec:buckingham}.

\section{Dimensional analysis}\label{sec:dimension}

Key to unifying the work on scales of measurement and the statistical invariance principle is DA.
DA has a long history, beginning with the discussion of dimension and measurement \citep{fourier1822theorie}. Since DA is key to the description of a natural phenomenon, DA lies at the root of physical modelling. A phenomenon's description begins with the phenomenon's features, each of which has a dimension, e.g. `mass' ($ M $) in physics or `utility' ($ U $)  in economics. Each dimension is assigned a scale e.g. `categorical', `ordinal', `ratio', or `interval', a choice that might be dictated by practical as well as intrinsic considerations. Once the scales are chosen, each feature is mapped into a point on its scale. For a quantitative scale, the mapping will be made by measurement or counting, for a qualitative scale, by assignment of classification. Units of measurement may be assigned as appropriate for quantitative scales, depending on the metric chosen. For example, temperature might be measured on the Fahrenheit scale, the Kelvin scale or the Celsius scale.  This paper will be restricted to quantitative features, more specifically those features on ratio- and interval-scales.  

\subsection{Foundations}
\label{subsec:dh}

One tenet of working with measured quantities is that  units in an expression or equation must \lq\lq match up\rq\rq;  relationships among measurable quantities require  {\em{\DHo}}. To check the validity of comparative statements about two quantities, $X_1$ and $X_2$, such as $X_1 = X_2$,  $X_1 < X_2$ or $X_1>X_2$, $X_1$ and $X_2$ must be the same dimension, such as time. To add $X_1$ to $X_2$,  $X_1$ and $X_2$ must not only be the same dimension but must also be on the same scale and expressed in the same units of measurement. 

To discuss this explicitly, we use a standard notation \citep{jcgm2012200} and write a measured quantity $X$ as $X=\{X\} [X]$, where $\{X\}$ is the numerical part of $X$.  $[X]$ may be read as the dimension of $X$ e.g. $[X] = L$ for length, for instance, or the units of $X$ on the chosen scale of measurement, e.g. $[X] = cm$. The latter by its nature means that the dimension is $L$. If $[X]=[1]$, then we say that $X$ is unitless or dimensionless. We define 1 or any number to be unitless, i.e., $1={1}[1]$, unless stated explicitly.

To develop an algebra for measured quantities,  for a function $ f $ we must say what we mean by  $\{f( X )\}$ (usually easy) and $[f( X )]$ (sometimes challenging). The path is clear for $f$ a simple function. For example, consider  $f(X)=X^2$.  Clearly we must have $X^2= \{X\}^2  [X]^2$, yielding, say, (3 inches)$^2 = 9$ inches$^2$. But what if $ f $ is a more complex function? This issue will be discussed in general in Subsection \ref{sec:fcts} and in detail for $f(x) = \ln (x)$ in Subsection \ref{sec:logarithm}.

For simple functions, the manipulation of both numbers and units is governed by an algebra of rules referred to as {\em{quantity calculus}} \citep{taylor2018quantity}. This set of rules states that $x$ and $y$
\begin{itemize}
	\item can be added, subtracted or compared if and only if $[x]=[y]$;
	\item can always be multiplied to get $xy =\{xy\}[xy]$ where $\{xy\}=\{x\}\{y\}$ and $[xy]=[x][y]$;
	\item can always be divided when $\{x\} \neq 0$ to get $y/x = 
	\{y/x\}[y/x]$ where $\{y/x\}=\{y\}/\{x\}$
	and $[y/x] = [y]/[x]$;
	\end{itemize}
	and that 
	\begin{itemize}
	\item $x$ can be raised to a power that is a rational fraction $\gamma$, provided that the result is not an imaginary number,  to get $x^\gamma = \{x\}^\gamma [x]^\gamma$.
\end{itemize}
Thus it makes sense to square-root transform ozone $O_3 = \{O_3\} $ parts per million (ppm) as 
$\{O_3\}^{1/2}$ ppm$^{1/2}$ since ozone is measured on a ratio-scale with a true origin of $0$  and hence must be non--negative \citep{dou2007dynamic}. These rules can be applied iteratively a finite number of times to get expressions that are combinations of products of quantities raised to powers,  along with sums and rational functions of such expressions.  
  
This subsection concludes with an example that  demonstrates the use of DH ~and quantity calculus.

\begin{example}\label{ex:canadianmodel}
This example concerns a structural engineering model for 
lumber strength now called the ``Canadian model''\citep{foschiyao1986dol}. Here $\alpha(t)$ is dimensionless and represents the somewhat abstract quantity of the damage accumulated to a piece of lumber by time $t$.  When $\alpha(t)=1$, the piece of lumber breaks. This is the only time when $\alpha(t)$ is observed. The Canadian model posits that  $\dot{\alpha}$, the derivative of $\alpha$ with respect to time,  satisfies
\begin{equation}\label{eq:canadianmodel}
  \dot{\alpha}(t) = a [\tau(t) - \sigma_0 \tau_s]_+^b
  ~+ ~c [\tau(t) - \sigma_0 \tau_s]_+^n ~\alpha(t),
\end{equation}
where $a$, $b$, $c$, $n$ and $\sigma_0$ are log-normally distributed random effects for an individual specimen of lumber, $\tau(t)$, measured in  pounds per square inch (psi), is the stress applied  to the specimen cumulative to time $ t $, $ \tau_s $  (in psi) is the specimen's short term breaking strength if it had experienced the stress pattern $ \tau(t)=kt $ for a fixed known $ k $ (in psi per unit of time), and $ \sigma_0 $  is the unitless stress ratio threshold. The expression $[ t ]_+$ is equal to $t$ if $t$ is non-negative and is equal to 0 otherwise. Let $ T_F $ denote the random time to failure for the specimen, under the specified
stress history curve, meaning $ \alpha(T_F) = 1 $. 

As has been noted \citep{kohler2002probabilistic,hoffmeyer2007duration,zhai2012a,wong2018dimensional}, this model is not  dimensionally homogeneous.  In particular, the units associated with both terms on the right hand side of Equation (\ref{eq:canadianmodel})  involve random powers, $ b $ and $ n $,  leading to random units, respectively (psi)$^b$ and (psi)$^n$. As noted by  \citet{wong2018dimensional},  the  coefficients $a$ and $c$ in Equation (\ref{eq:canadianmodel}) cannot involve these random powers and so cannot compensate to make the model dimensionally homogeneous. 

Rescaling  is a formal way of addressing this problem. \cite{zhai2012a} rescale by setting $\pi(t) = \tau(t)/\tau_s$. They let $\mu $ denote the population mean of $\tau_s$ and write a modified equation of (\ref{eq:canadianmodel}) as the dimensionally homogeneous model
\[
 \mu \dot{\alpha}(t) = a^*  [\pi(t) - \sigma_0 ]_+^b
  ~+ ~c^* [\pi(t) - \sigma_0 ]_+^n ~\alpha(t).
\]
In contrast, \citet{wong2018dimensional} propose another dimensionally homogeneous model
\begin{eqnarray}\nonumber
 	\mu \dot{\alpha}(t)  
  	&=& [(\tilde{a} \tau_s )(\pi(t) - \sigma_0)_+]^b ~
  + ~ [(\tilde{c} \tau_s )(\pi(t) - \sigma_0)_+]^n ~ \alpha(t), \nonumber
\end{eqnarray}
where $\tilde{a}$ and $\tilde{c}$ are now random effects with dimensions Force$^{-1} \cdot $ Length$^2$.  
\end{example}

We see that there may be several ways to nondimensionalize a model. Another method, widely used in the physical sciences, involves always normalizing by  the standard units specified by the Syst\`eme International d'Unit\'es (SIU), units such as meters or kilograms. So when the dimensions of a non--negative quantity $ X $ like absolute temperature have an associated SIU of $Q_0 = \{ 1 \}[Q_0]$, $X$ can be converted to a unitless quantity by first expressing $X$ in SIUs and then by using quantity calculus to rescale it as $X/Q_0$.
The next example provides an important illustration of the application of 
the standardized unit approach.

\begin{example}\label{ex:phindex}
Liquids contain both hydrogen  and hydroxide ions. In pure water these ions appear in equal numbers. But the water becomes acidic when there are more hydrogen ions and basic when there are more hydroxide ions. Thus acidity is measured by the concentration of these ions. The customary measurement is in terms of the hydrogen ion concentration, denoted $H^+$ and measured in the SIU of one mole of ions per litre of liquid. These units are denoted c$^o$ and thus, in our notation, $[H^+]=$ c$^o$. However for substantive reasons, the pH index for the acidity of a liquid is now used to characterize acidity. The index is defined by   $pH = - \log_{10} ( H^+/$c$^o )$. Distilled water  has a $ pH = 7 $ while lemon juice has a $ pH  $ level of about $ 3 $. Note that  $ \{H^+\} \in (0,\infty) $ lies on a ratio-scale while $pH$ lies on an interval-scale $ (-\infty,\infty) $ -- the transformation has changed the scale of measurement.

\end{example}


\subsection{The problem of scales} \label{subsec:scales}
The choice of scale restricts the choice of units of measurement, and these units dictate the type of model that may be used. However, comparing the size of two quantities on a ratio-scale must be made using their ratio, not their difference, whereas the opposite is true on an interval-scale where differences must be used.


Thus we need to study scales in the context of model building and  hence in the context of quantity calculus. In  his celebrated paper, \citet{stevens1946theory} starts by proposing four major scales for measurements or observations: categorical, ordinal,  interval and ratio. This taxonomy is based on the notion of permissible transformations as is the work of our Section \ref{sec:relationships}. However, our work is aimed at modelling while Stevens's work is aimed at statistical analysis. Stevens defines permissible transformations as follows. He allows permutations as the transformations of data on all four scales, allows strictly increasing transformations of data on the ordinal, ratio and interval-scales, allows scalar transformations ($f(x)=ax$) of data on the ratio and interval-scales and allows linear transformations ($f(x)=ax+b$) of data on the interval-scale.

Stevens created his taxonomy as a basis for classifying the family of all statistical procedures for their applicability in any given situation \citep{stevens1951mathematics}. 
And for instance, \citet{luce1959possible} points out that, for measurements made on a ratio-scale, the geometric mean would be appropriate for estimating the central tendency of a population distribution \citep{velleman1993nominal}. In contrast, when measurements are made on an interval-scale the arithmetic mean would be appropriate. The work of Stevens seems to be well-accepted in the social sciences, with \citet{ward2017stevens} calling his work monumental. But Steven's work is not widely recognized in statistics. \citet{velleman1993nominal} reviews the work of Stevens with an eye on potential applications in the then emerging area in statistics of artificial intelligence (AI), hoping to automate data analysis. They claim that ``Unfortunately, the use of Stevens's categories in selecting or recommending statistical analysis methods is inappropriate and can often be wrong.'' They describe alternative scale taxonomies for statistics that have been proposed, notably by \citet{mosteller1977data}. A common concern centers on the inadequacies of an automaton to select the statistical method for an AI application. Even the choice of scale itself will depend on the nature of the inquiry and thus is something to be determined by humans. For example, length might be observed on the relatively uninformative ordinal scale $\{ short,~  medium~ , long \} $, were it  sufficient for the intended goal of a scientific inquiry, rather than on the seemingly more natural  ratio-scale $(0,\infty) $. 

\subsection{Why the origin matters}\label{subsec:why}
The  interval-scale of real numbers  allows for the taking of sums, differences, products, ratios, and integer powers of values observed on that scale. Rational powers of nonnegative values are also allowed although irrational powers lead into the domain of transcendental functions and difficulties of interpretation. The same operations are allowed for a ratio-scale of real numbers provided that the differences are non-negative. So superficially, these two scales seem nearly the same. 

But there is a substantial qualitative difference between ratio- and interval-scales, so ignoring the importance of scale  when building models can result in challenges in interpretation.
The issue has to do with the meaning of the $0$ on a ratio-scale.  The next hypothetical example illustrates the point.

\begin{example}\label{eg.whyzeromeans}
When the material in the storage cabinet  at a manufacturing facility has been depleted, the amount left is $0$.  
To understand the usefulness of this origin, consider if the facility's inventory monitoring program recorded a drop of $100 kg$ during the past month. Without the knowledge of the origin,  
of where the amount of inventory lies on the scale, the implications of this drop are unclear.  If the amount left in the facility is $99,900 kg$ the drop means one thing, while if the amount left is $50 kg$, the interpretation would be completely different. 

Since the amount of inventory lies on the ratio-scale, these changes should instead be reported using ratios. The recorder in the facility in the first case would report a decline in inventory of  $100/100,000 = 0.001$ or  $0.1\%$.  In the second case, the recorder  would report a decline of $ 100/150 = 2/3$ or $66.7 \%$, the same drop but with a totally different meaning.  This example explains why stock price changes are reported on a ratio-scale, as a percentage, and not on an interval-scale.
\end{example}



\section{Transforming quantities}\label{sec:new_transforming}

The scale of the  measurement $ X $ may be transformed in a variety of ways. No change of scale occurs when the transformation is a rescaling, where we know how to transform both the numerical part of $ X $ and $ X $'s units of measurement.  When the transformation is complex,  the scale itself might change.  For instance, if $ X $ is measured on a ratio-scale, then the logarithm of $ X $ will be on an interval-scale. 
Observe that in Example \ref{ex:phindex}, the units of measurement in $ H^+ $ were eliminated before transforming by the    transcendental function $\log_{10}$. That raises the question: do we need to eliminate units before applying the logarithm?     This question and the logarithmic transformation in science have led to vigorous debate for over six decades \citep{matta2010can}.  We highlight and resolve some of that debate below in Section \ref{sec:logarithm}.

However, we begin with an even simpler situation seen in the next subsection, where we study the issues that may arise when interval-scales are superimposed on ratio-scales.

\subsection{Switching scales}\label{sec:scalesinscales}

This subsection concerns a perhaps unconscious switch in a statistical analysis from a ratio  scale, which lies on $[0,\infty)$, to an interval-scale, which lies on $(-\infty, \infty)$.  
\subsubsection*{The bell curve approximation.}\label{subsec:bellcurve}
Despite the fundamental difference between the ratio-scale and interval-scales, the normal approximation is often used to approximate the sampling distribution for a ratio-valued response quantity.  This in effect replaces the ratio-scale with an interval-scale. 
In this situation, what should be used is the truncated normal distribution approximation, although this introduces undesired complexity. For example, if the approximation for the cumulative distribution function CDF were  
$ P ( X \leq x \mid X > 0) $ where 
$X \sim N(\mu, \sigma^2$, we would have 
\begin{equation}
\label{eq:trunnormalmean}
	E(Z \mid X> 0) = 
\frac{\phi(\mu/\sigma)}{\Phi(\mu/\sigma)} 
\end{equation}
where $ Z = (X -\mu)/\sigma$ 
,  while  $\phi$ and $\Phi$ denote respectively, the standardized Gaussian distribution's probability density function and CDF. Observe that in Equation (\ref{eq:trunnormalmean}) is invariant under changes of the units in which $X$ is measured, as it should be.  Furthermore the expectation of $Z$ would be approximately $0$ if $\mu/\sigma$ were large compared to $0$ as it would be were the non-truncated Gaussian distribution for an interval-scale imposed on this ratio-scale.  This would occur if the mean  $\mu$  were much larger than the standard distribution $\sigma$. 
That suggests the bell curve approximation would not work well as an approximation were the population under investigation widely dispersed. For example it might be satisfactory if $X$  represented the height of a randomly selected adult woman, but not if it were the the height of a randomly selected human female.

As mentioned at the beginning of this section, this switch occurs when approximating a distribution. 
This switch is ubiquitous and seen in most elementary statistics textbooks.  The assumed Gaussian distribution model leads to the sample average as a measurement of the population average instead of the geometric mean, which should have been used \citep{luce1959possible}. That same switch is made in such things as regression analysis and the design of experiments. The seductive simplicity has also led to the widespread use of the Gaussian process in spatial statistics and machine learning.

The justification of the widespread use of the Gaussian approximation may well lie in the belief that the natural origin $ 0 $ of the ratio-scale lies well below the range of values of  $ X $ likely to be found in a scientific study.   
This may well  be the explanation of the reliance on interval-scales for temperature in Celsius and Fahrenheit on planet Earth at least since one would not expect to see temperatures anywhere near the true origin of temperature  $ 0^{\circ}K$ on the Kelvin scale (ratio) that corresponds to $-273^{\circ} C$ on the Celsius scale (interval). We would note in passing that these two interval-scales for temperature also illustrate the statistical invariance principle (see Subsection \ref{subsec:interval});
each scale is an affine transformation of the other.

We illustrate the difficulties that can arise when an interval-scale is misused in a hypothetical experiment where measurements are made on a ratio-scale, with serious consequences.

 \begin{example}\label{ex:switch}
The justification above for the switch from a ratio to an interval-scale can be turned into a simple approximation that may help with the interpretation of the data.
To elaborate, suppose interest lies in comparing two values of $ X $, $x_1$ and $x_2$, that lie in a ratio-scale with $0 < a <   x_1 < x_2$ for a known $ a $. Interest lies in the relative size 
of these quantities i.e. on $ r = x_2 / x_1$. An approximation to $ r $   through the first order Taylor expansion of $x_2/x_1$ at $(x_1,x_2)=(a,a)$ yields $r \approx 1+ (x_2 - x_1)/a $, thus providing an approximation to $r$ on an interval-scale. 
 For instance, with $a=120cm, ~x_1=150 cm,$ and $x_2=180 cm$, the ratio is $r=1.20$ and the approximation, $1.25$. Both are unitless. This points to the potential value of rescaled ratio data when a Gaussian approximation is to be used for the sampling distribution of a quantity on a ratio-scale.
 \end{example}
 
 \vskip .1in

\subsection{Algebraic versus transcendental functions}\label{sec:fcts} 
 A function $u$, which describes the relationship among quantities $X_1,\dots,X_p$ as 
\[
u(X_1,\dots, X_p)=0,
\]
may be  a sequence of transformations or operations involving the $X_i$'s, possibly combined with  parameters. We know how to calculate the resulting units of measurement when $u$  consists of a finite sequence of permissible algebraic operations.  The function consisting of the concatenation of such a sequence may formally be defined as a root  of a polynomial equation that must satisfy the requirement of \DHo ~(other desirable properties of $u$ along with methods for determining an allowable $u$ are discussed in Section \ref{sec:relationships}). Such a function is called algebraic.  

But $u$ may also involve non-algebraic operations leading to non-algebraic functions called transcendental (because they ``transcend'' an algebraic construction).  Examples in the univariate case ($p=1) $ are $ \sin(X) $ and  $ \cosh(X) $ and, for a given nonnegative constant $\alpha $,  $ \alpha^X $ and $ \log_\alpha{(X)} $. The formal definition of  a non-algebraic function does not explicitly say whether or not such a function can be applied to quantities with units of measurement. \citet{bridgman1931dimensions} sidesteps this issue by arguing that it is mute since valid representations of natural phenomena can always be nondimensionalized (see Subsection \ref{subsec:buckingham}).   But the current Wikipedia entry on the subject states ``transcendental functions are notable because they make sense only when their argument is dimensionless'' \citep{wiki:tf}.  The next subsection explores the  much used class of Box-Cox family \citep{box1964analysis} that includes transcendental functions.

\subsection{The Box-Cox tranformation}\label{sec:boxcox} 
 Frequently in statistical modelling,  a transformation is used to extend the domain of applicability of a procedure that assumes normally distributed measurements \citep{de1997bayesian}.  That transformation may also be seen as a formal part of statistical model building that facilitates maximum likelihood estimation of a single parameter 
\citep{draper1969distributions}.  The Box-Cox (BC)  transformations constitute an important class of such transformation and therefore a standard tool in the statistical scientist's toolbox. 

In its simplest form, a member of this family of transformations has the form of a function $ bc(X) = X^{\lambda} $ for a real-valued parameter $ \lambda \in (-\infty, \infty)$.  Here $ X $ would need to lie in $[0,\infty) $, a ratio-scale, to avoid potential imaginary numbers. However, in practice interval-scales are sometimes allowed,  a positive constant being added to avoid negative realizations of $ X$. This ad hoc procedure thus validates the use of a Gaussian distribution to approximate the sampling distribution for $ X $. 

Since $ X $ is  measured on a ratio--scale, for any two points on that scale,  $ bc(x_2/x_1)~=~ bc(x_2) / bc(x_1) $ while the scale is equivariant under multiplicative transformations, i.e., $ bc( a x) ~=~ a^{\lambda}~ bc( x) $ for any point on that scale.  Finally $ bc(X) > 0$ so that the result of the transformation also lies on a ratio--scale, even when its intended goal is an approximately Gaussian distribution for the (transformed) response.  

 \citet{box1964analysis} actually state their transformation as 
\begin{equation}\label{eq:boxcox}
bc_{\lambda}(X) = 
\begin{array}{cc}
\frac{X^\lambda - 1}{\lambda},~ & (\lambda \neq 0)\\
\end{array}
,
\end{equation} 
 that moves the origin of the ratio--scale from $ 0 $ to $ -1$ .  It is readily seen that unless $\lambda$ is a rational number, $bc_{\lambda},~\lambda \neq 0$ will be transcendental.  That fact would be inconsequential in practice in as much as a modeller would only ever use a rational number for $\lambda$ . Or at least that would be the case except that the BC class has been extended to include $\lambda = 0$ by admitting for membership $ \lim_{\lambda\rightarrow 0}bc_{\lambda} (X) = \ln(X),~X > 0$ .

 On closer inspection, we see that for validity, in Equation (\ref{eq:boxcox}), the $1$  needs to be replaced by $(1[ X ])^\lambda$ to include units of measurement. Then for $\lambda \neq 0$ the transformation becomes
\begin{equation}
\label{eq:boxcoxrevised}
	X_\lambda \doteq bc(X)=\frac{X^\lambda - 1}{\lambda} =
	\frac{\{ X \}^\lambda - 1^{\lambda}}{\lambda}[ X ]^\lambda.
\end{equation}
As $\lambda \rightarrow 0$, the only tenable limit would seem to be 
\begin{equation}
	X_0 = \ln{(\{ X \})}
	\nonumber
\end{equation}
not $\ln{(X})$. In other words, in taking logarithms in the above example, the authors may have unconsciously nondimensionalized the measurements. Taking antilogarithms would then return $\{X\}$, not $X$.

 Equation \ref{eq:boxcoxrevised} thus tells us the Box--Cox transformation may have the unintended consequence of transforming not only the numerical value of the measurements but also its units of measurement, which become $[ X ]^\lambda$. This  makes a statistical model difficult to interpret.  For example imagine the challenge of a model with a random response in units of $mm^{1/100}$.  And since the transformation is nonlinear, returning to the original scales of the data  would  be difficult. For example, unbiased estimates would become biased and their standard errors could be hard to approximate.

\begin{remarks}\label{rem:boxcox2}
\citet{box1964analysis} do not discuss the issue of scales in relation to the transformation it introduces. In that paper's second example the logarithmic transformation is applied to $X$, the number of number of cycles to failure, which may be regarded as unitless.	
\end{remarks}

 \begin{remarks}\label{rem:boxcox}
 The mathematical foundation of the Box--Cox family is quite complicated. Observe that in Equation  (\ref{eq:boxcox}),  if  $ \lambda $ is a rational number $ m/n $ for some nonnegative integers $ m $ and $ n $,  $ bc_{\lambda} $ will be an algebraic function of $ X $. So as $\lambda $ varies over its domain, $(-\infty,\infty )$, the  function flips back and forth from an algebraic to a transcendental function. For any fixed $ m $ and point $ x $, as $n$ approaches infinity, the trajectory of 
\[
\{bc_{m/n}(x): n= 1,2, \dots, \}
\]
converges to $\ln {x} $, so the family now includes the logarithmic transformation as another transformation in the statistical
analyst's toolbox, which is used when the response distribution appears to have a long right tail.  Thus a transcendental function has been added to the family of algebraic transformations obtained when $\lambda $ is chosen to be a positive rational number. It does not seem to be known if all transcendental transformations lie in the closure of the class algebraic functions under the topology of pointwise convergence. However, when this family is shifted from the domain of mathematical statistics in the human brain to that of computational statistics in the computer's processor, this complexity disappears. In the computational process, all functions are algebraic and neither the logarithm nor infinity exist.
\end{remarks}


The importance of the logarithmic transformation in statistical and scientific modelling, and issues that have arisen about its lack of units, leads next to a special subsection devoted to it.

\subsection{The logarithm: a transcendental function}{\label{sec:logarithm}} 

\subsubsection*{Does the logarithm have units?}
We have argued (see Example \ref{ex:mle}) that the answer is ``no.''  First consider applying the logarithm to a unitless quantity $ x $.  It is sensible to think that its value  will have no units, and so we take this as fact.  But what happens if we apply the logarithm to a quantity with units? For instance, is log(12 inches) = log(12) + log(inches)? This issue has been debated for decades across different scientific disciplines; we summarize recent debates in Appendix \ref{app.debate}.
 
We now discuss this issue in more detail and argue that the result must be a unitless quantity. We use the definition of the natural logarithm of $ x $ as the area under the curve of the function 
$ f(u) = 1 / u $ \citep{molyneux1991dimensions}. We follow the notation defined in Section 3.1, and for clarification, we write ``1'' as $ y \equiv 1[x]$ and $ u=\{u\} [x]$.   We then make  the change of variables $ v = u /y   $ so that $ v $ is unitless and get 
\begin{equation}
	 \ln(x) = \int_y^x  \frac{1}{ u }~ d (u) = \int_1^{x/y} \frac{1}{ y~v } ~d (y~ v) 
	= \int_1^{x/y} \frac{1}{ v } d (v), 
	\label{eq:integraldefn}
\end{equation}
which is a unitless quantity, as claimed. 
	
We now derive the more specific result, that $\ln (x ) = \ln(\{x\})$.   In other words,   applying this transcendental function to a dimensional quantity $ x $ simply causes the units to be lost. We show below, from first principles, that for $v$ unitless,
	 \[
	\frac{d\ln( v )}{dv} = \frac{1}{v}.
	\]
This implies that
\[
		 \ln(w) = \int_{1}^{w} \frac{1}{v}~ dv ,
\]
which, combined with Equation (\ref{eq:integraldefn}), implies that $\ln (x) = \ln(\{x\})$.

To show that the derivative of $\ln ( v) $ is $1/v$, 
we turn to the  original definition of the natural logarithm as the inverse of another transcendental function $\exp ( v) $, at least if $v  > 0 $. In other words 
$
 v = \exp~ ( \ln v ),~v > 0.
$
The chain rule now tells us that 
\[
1 = \frac{d\ln( v )}{dv}  \exp~ ( \ln v ).
\]
Thus
\[
\frac{d\ln( v )}{dv} =  \exp~ (-\ln v ) =  
\frac{1}{v}
\]
for any real $ v > 0 $.

\subsubsection*{Can we take the logarithm of a dimensional quantity with units?}{\label{sec:log.arg.units}}

We argue that the answer is ``no.''  We reason along the lines of  \cite{molyneux1991dimensions}, who sensibly argues that, since $\ln x$ has no units even when $x$ has units,  the result is meaningless.  In other words, since it is disturbing that the value of a function is unitless, no matter what the argument, we should not take the logarithm of a dimensional quantity with units. 

This view  agrees with that of 
\citet{meinsma2019dimensional}. He notes that the Shannon entropy of a probability density function $f$, which has units, is defined in terms of the $ln~f(x)$. He notes that Shannon found this to be an objectionable feature of his entropy but rationalized its use nevertheless. But not Meinsma, who concludes 
``To me it still does feel right....''

To consider the ramifications of ignoring this in a  statistical model, suppose that $z$ is some measure of particulate air pollution in the logarithmic scale with $ z = \ln x $ where $x$ is a measurement with units.  
This measurement appears as $\beta ~z $ in a scientific model of the impact of particulate air pollution on health \citep{cohen2004urban}.
In this model, even though $z$ is unitless, its numerical value depends on the numerical value of $x$, via $\{z\} = \ln\{x\}$. 
Thus the numerical value of $z$ depends on the units of measurement of $x$.
But, since $z$ itself is unitless, we cannot adjust $\beta$ to reflect changes in the units of $x$.
To make this point explicit, suppose that experimental data pointed to the value
$ \beta = 1,101,231.52$.
We have no idea if air pollution was a serious health problem.   
Thus, we see the problem that arises with a model that involves the logarithm of a measurement with units. 
This property of the logarithm points to the need to nondimensionalize $ x $ before applying the logarithmic transformation in scientific and statistical modelling, in keeping with the theories of Buckingham, Bridgman
and Luce.

One of the major routes taken in debates about the validity of applying the natural logarithm to a dimensional quantity involves arguments based one way or another on a Taylor expansion. A key feature of these debates involves the claim that the expansion is impossible, since the terms in the expansion have different units and so cannot be summed \citep{mayumi2010dimensions}. We show below that this claim is incorrect by showing that all of the terms in the expansion have no units (see Appendix \ref{app.debate} for more details).

Key to the Taylor expansion argument of validity is how to take the derivative of $\ln x$ when $ x $ has units. Recall that above, we calculated the derivative of 
$ \ln v $  for $ v $ unitless. 
To define the derivative of $
\ln x $ when $ x $ has units, we proceed from first principles.
Suppose we have a function $ f $ with argument $x = \{x\}[x]$.  We define the derivative of $f$ with respect to $x$ as follows.  Let $\Delta = \{ \Delta\} [\Delta]$ and $x = \{x\} [x]$ and suppose that   $\Delta$ and $x$ have the same units, that is, that $[\Delta] = [x]$.  Otherwise, we would not be able to add $x$ and $\Delta$ in what follows. 
Then we define
\begin{eqnarray}
\label{eq:log.derivative}
f'(x) &\equiv&  \lim_{ \{\Delta\} \to 0}  \frac{ f(x + \Delta) - f(x) } { \Delta} 
 \\
 &~&
  \nonumber \\
&=&   \lim_{ \{\Delta\} \to 0}  \frac{ f( \{x + \Delta\} [x + \Delta] ) - f(\{x\}[x]) } {  \{\Delta\}  [\Delta]}
   \nonumber \\
  &~&
  \nonumber \\
  &=&
  \lim_{ \{\Delta\} \to 0}  \frac{ f( \{x + \Delta\} [x ] ) - f(\{x\}[x]) } {  \{\Delta\}  [x]}. 
 \nonumber
\end{eqnarray}
 For instance, for $f(x)=x^2$
\begin{eqnarray}
\frac{d}{dx} x^2 &=&  \lim_{ \{\Delta\} \to 0}  \frac{  \{ x+\Delta\}  ^2  [x]^2 -  \{x\}^2 [x]^2 } { \{\Delta\} [x]}
 =
  \lim_{ \{\Delta\} \to 0}  \frac{  \{x + \Delta\}^2  - \{x\}^2 } { \{\Delta\} } \times [x]
  \nonumber \\
 &=&  2\{x\} [x] = 2 x.
 \nonumber \end{eqnarray}
 
To use Equation (\ref{eq:log.derivative}) to differentiate $f(x) = \ln (x)$,  
we first write  
\[
\ln (x + \Delta )  - \ln x = \ln \{x + \Delta \} - \ln \{ x \}.
\]
So
\[
\frac{d}{dx}  \ln x  = 
  \lim_{ \{\Delta\} \to 0}  \frac{ \ln \{x + \Delta\}  - \ln \{x\} } {  \{\Delta\}  [x]}  
  = \frac{d\ln \{ x \}}{ d \{ x \} } \times 
  \frac{1}{ [x] }
  =
      \frac{ 1  } { \{ x\} } \frac{ 1}{  [x]} = \frac{1}{x}.
\]

Using this definition of the derivative we can carry out a Taylor series expansion about $ x=a>0 $ to obtain
\begin{equation}\label{eq:logexpansion}
	\log(x) = \log(a) +  \sum_{k=1}^{\infty} g^{(k)}(a) \frac{(x-a)^k}{k!},
\end{equation}
where 
\begin{equation*}
g^{(k)}(a)=\bigg[ d^k \log(x) /dx^k \bigg|_{x=a}.
\end{equation*}
As $ g'(x) = 1/x $, the first term, $g'(x) (x-a)$, in the infinite summation is unitless.  Differentiating $g'(x)$ yields $g''(x) = 1/x^2$ and once again, we see that the term $g''(x) (x-a)^2/2$ is unitless.  Continuing in this way, we see that the summation on the right side of equation (\ref{eq:logexpansion}) is unitless, and so the equation satisfies dimensional homogeneity. This reasoning differs from the  incorrect reasoning of \citet{mayumi2010dimensions} in their argument that the logarithm  cannot be applied to quantities with units because the terms in the Taylor expansion would have different units. Our reasoning also differs from that of \citet{baiocchi2012dimensions} who uses a different expansion to show that the logarithm cannot be applied to measurements with units, albeit without explicitly recognizing the need for 
$ \ln x $ to be unitless. 
The expansion in Equation 
(\ref{eq:logexpansion}) is the same as  that given in \citet{matta2010can}, albeit  not in an 
explicit form for $ \ln x $. Like us, they do discredit the Taylor expansion argument  against applying $\ln x$ to quantities with units. 

\section{Allowable relationships among quantities}\label{sec:relationships}

Having explored dimensional analysis and the kinds of difficulties that can arise when scales or units are ignored, we turn to a key step in unifying physical and statistical modelling. 
We now determine how to relate quantities and hence how to specify the `law' that characterizes the phenomenon which is being modelled.  

But what models may be considered legitimate? Answers for the sciences, given long ago, were based on the principle that for a model to completely describe a natural phenomenon, it cannot depend on the units of measurement that might be chosen to implement it.  This answer was interpreted in two different ways.  In the first interpretation, the model must be nondimensionalizable, i.e., it cannot have scales of measurement and hence cannot depend on units.  In the second interpretation, the model must be  invariant under all allowable transformations of scales. 
Both of these interpretations reduce the class of allowable relationships that describe the phenomenon being modelled and place restrictions on  the complexity of any experiment that might be needed to implement that relationship.



\subsection{Buckingham's Pi-theorem}\label{subsec:buckingham}

The section begins with Buckingham's simple motivating example.
\setcounter{example}{8}
\begin{example}\label{ex:buckinghams}
This example is a characterization of properties of gas in a container, namely, a characterization of the relationship amongst the pressure ($ p $), the volume ($ v $), the number of moles of gas ($ N $) and the absolute temperature ($ \theta $) of the gas.  The absolute temperature reflects the kinetic energy of the system and is measured in degrees Kelvin ($ ^{\circ} K $), the SIUs for temperature.  A fundamental relationship amongst these quantities is given by 
\begin{equation}\label{eq:buckinghamexample}
\frac{pv}{\theta N} - D = 0 
\end{equation}
for some  constant $ D $ that doesn't depend on the gas. Since the dimension of $p v /(N\theta)$ is
(force $\times$ length$^3$)/(\# moles $\times$ temperature), as expressed, the relationship in Equation (\ref{eq:buckinghamexample}) depends on the units associated with $p,v$ and $\theta$, whereas the physical phenomenon underlying the relationship does not.  Buckingham gets around this by invoking a parameter $ R$ ($\equiv D$) with units  (\# moles $\times$ temperature)/(force $\times$ length$^3$). He rewrites Equation (\ref{eq:buckinghamexample}) as 
\[
\frac{pv}{ R\theta N} - 1 = 0. 
\]
Thus $\pi = pv/(R\theta N) $ has no units. Buckingham calls this equation  complete and hence nondimensionalizable. 
This equation is known as the Ideal Gas Law, with $R$ denoting the ideal gas constant \citep{idealgaslaw}.   
\end{example}

This example of nondimensionalizing by finding one expression, $\pi$, as in Equation (\ref{eq:buckinghamexample}) can be extended to cases where we must nondimensionalize by finding  several $\pi$ quantities. 
This extension is formalized in Buckingham's Pi-theorem.
Here is a formal statement (in slightly simplified form) as stated by 
\citet{buckingham1914physically} and discussed in a modern style in \cite{bluman.cole}.
\begin{theorem}\label{thm:buckingham1}
Suppose   $ X_1,\dots, X_p $ are $ p $ measurable quantities satisfying a defining 
 relation
\begin{equation}\label{eq:buckingham}
u(X_1,\dots, X_p)=0
\end{equation}
that is  dimensionally homogeneous.
In addition, suppose that there are $ m $ dimensions appearing in this equation, denoted $ L_1,\ldots, L_m $, and that the dimension of $ u $ can be expressed $ [u]= L_1^{\alpha_1} \times \cdots \times L_m^{\alpha_m} $ and the 
dimension of each $ X_j $ can be expressed as $[X_j]= L_1^{\alpha_{j1}} \times \cdots \times L_m^{\alpha_{jm}} $.
Then Equation (\ref{eq:buckingham}) implies the existence of  $ q $ fundamental quantities,  $ q \ge p-m $ dimensionless quantities $ \pi_1,\dots,\pi_q $ with $ \pi_i = \Pi_{j=1}^p X_j^{a_{ji}} ,~i=1,\dots, q$, and a function $ U $ such that 
\[
U(\pi_1,\dots, \pi_q)=0.
\]
\end{theorem}
\noindent In this way $u$ has been nondimensionalized.   The choice of $\pi_1, \ldots, \pi_q$ in general is not unique.

The theorem is proven constructively, so we can find $\pi_1, \ldots, \pi_q$ and $U$.
We first determine the $m$ fundamental dimensions used in $X_1,\ldots, X_p$.  We then use the quantities $X_1,\ldots, X_p$ to construct two sets of variables:  a set of $m$  primary variables also called repeating variables and a set of $q$ secondary variables, which are nondimensional.
For example, if  $X_1$ is the length of a box and $X_2$ is the height and $X_3$ is the width, then there is $m=1$ fundamental dimension, the generic length denoted $L$. 
We can choose $X_1$ as the primary variable and use $X_1$ to define two new variables $\pi_1=X_2/X_1$ and $\pi_2=X_3/X_1$.  These new variables, called secondary variables, are dimensionless.
Buckingham's theorem states the algebraic equation relating $X_1, X_2$ and $X_3$ can be re-written as an equation involving only $\pi_1$ and $\pi_2$.  Note that we could have also chosen either $X_2$ or $X_3$ as the primary variable.

A famous application of Buckingham's theorem concerns the  discovery of the Reynold's number in fluid dynamics, which is discussed in \citet{gibbons2011dimensional}.  For brevity we include that example in  Appendix \ref{sec:reynolds}. 

A link between Buckingham's approach and statistical modelling was recognized in the paper of 
\citet{albrecht2013experimental} and commented on in \cite{lin2013comment}. But its link with the statistical invariance principal seems to have been first identified in the thesis of \cite{shen2015dimensional}.
This connection provides a valuable approach for the statistical modelling of scientific phenomena.
Shen builds a regression model starting with Buckingham's approach and thereby a nondimensionalized relationship amongst the variables of interest. We propose a different approach in Section \ref{sec:invariance}. We present Shen's illustrative example next. 

\begin{example} \label{ex:shen}
This example, from  \citet{shen2015dimensional}, concerns a model for the predictive relationship between 
the volume $ X_3 $ of wood in a pine tree and its height $ X_1 $ and diameter $X_2 $.
The dimensions are  $[X_1] = L $, $ [X_2] = L$ 
and $ [X_3] = L^3 $.  Shen chooses $ X_1 $ as the repeating variable and calculates the $\pi$-functions $ \pi_1 =X_2 X_1^{-1} $ 
and $ \pi_2 = X_3 X_1^{-3}$. He then applies the Pi-theorem to get the dimensionless version of the relationship amongst the variables:
\[
\pi_2 = g (\pi_1)
\]
for some function $ g $.  He correctly recognizes that $ ( \pi_1,\pi_2) $ is the maximal invariant under the scale transformation group, although the connection to the ratio-scale of Stevens is not made explicitly.  He somewhat arbitrarily chooses the class of relationships given by 
\begin{equation}\label{eq:shenlinearizes}
 \pi_2 = k \pi_1^\gamma.
 \end{equation}
  He linearizes the model in Equation (\ref{eq:shenlinearizes}) by taking the logarithm 
  and adds a residual to get a standard regression model, susceptible to standard methods of analysis.  In particular the least squares estimate $ \hat{\gamma} = 1.942 $ turns out to provide a good fit judging by a scatterplot.  
  
  Note that application of the logarithmic transformation is justified since the $\pi$-functions are dimensionless. 
\end{example}

\subsection{Bridgman's alternative}\label{subsec;bridgman}

We now describe an alternative to the approach of \citet{buckingham1914physically} due to \citet{bridgman1931dimensions}. At around the same time that Edgar Buckingham was working on his Pi-theorem, Percy William Bridgman was giving lectures at Harvard on the topic of nondimensionalization that were incorporated in a book whose first edition was published by Yale University Press in 1922. The second edition came out in 1931 \citep{gibbons2011dimensional}.  Bridgman thanks Buckingham for his papers but notes their approaches differ. And so they do.  For a start, Bridgman asserts his disagreement with the position that seems to underlie Buckingham's work that ``a dimensional formula has some esoteric significance connected with the `ultimate nature' of things.''  Thus those that espouse that point of view  must ``find the true dimensions and when they are found, it is expected that something new will be suggested about the physical properties of the system.'' Instead, Bridgman takes measurement itself as the starting point in modelling and even the collection of data:
``Having obtained a sufficient area of numbers by which the different quantities are measured, we search for relations between these numbers, and if we are skillful and fortunate, we find relations which can be expressed in mathematical form.'' He then seeks to characterize a measured quantity as either primary, that is, the product of direct measurement, or secondary, that is, computed from the measurements of primary quantities, as, for instance, velocity is computed from the primary quantities of length and time. Finally he sees the basic scientific issue as that of characterizing one quantity in terms of the others as in our explication of Buckingham's work above in terms of the function $u$.  

Bridgman considers the functional relationship between secondary and primary measurements under what statistical scientists might call ``equivariance'' under multiplicative changes of scale in the primary units. He proves that the functional relationship must be based on monomials with possible fractional exponents, not unlike the form of the $\pi$-functions above.  Thus   Bridgman  is able to re-derive Buckingham's $\pi$ formula, albeit with the added assumption that $u$ is differentiable with respect to its arguments.  

\subsection{Beyond ratio-scales}\label{subsec:alternatives}

Nondimensionalization seems more difficult outside of the domain of the physical sciences. For example, the dimensions of quantities such as utility cannot be characterized by a ratio-scale.  And the choice of the primary dimensions is not generally so clear, although \citet{baiocchi2012dimensions} does provide an example in macroeconomics where time $[T]$, money $[\$]$, goods $[R]$ and utility $[U]$ may together be sufficient to characterize all other quantities.  

Bridgman's results on allowable laws were limited to laws involving quantities measured on ratio-scales.  A substantial body of work has been devoted to extending these results to laws involving quantities measured on nonratio-scales,  beginning with the seminal paper of  \citet{luce1959possible}.  To quote the paper by \citet{aczel1986scientific}, which contains an extensive review of that work, 
``Luce shows that the general form of a `scientific law' is greatly restricted by knowledge of the `admissible transformations' of the dependent and independent variables.'' 
It seems puzzling that this principle has been recognized little if at all in statistical science. This may be due to the fact fact that little attention is paid to such things as dimensions and units of measurement.

The substantial body of research that followed Luce's publication covers a variety of scales, e.g. ordinal, among other things.   Curiously that body of work largely ignores the work of Buckingham in favor of Bridgman even though the former preceded the latter. Also ignored is the work on statistical invariance described, which goes back to G. Hunt and C. Stein in 1946 in unpublished  but well-known work that led to optimum statistical tests of hypotheses.  

To describe this important work by Luce, we re-express Equation (\ref{eq:buckingham}) as
\begin{equation}\label{eq:buckingham2}
X_p = u^*(X_1,\dots, X_{p-1})    
\end{equation}
for some function $u^*$ and thereby define a class of all possible laws that could relate  $X_p$ to  the  predictors  $X_1,\dots, X_{p-1}$, before turning to a purely data-based empirical assessment of the possible $u^*$'s.   
Luce requires that $u^*$ satisfy an invariance condition. Specifically, he makes the strong assumption that the scale of each $X_i,~i=1,~\dots,~ p-1$, is  susceptible to a transformation $T_i \in {\cal F}_i$, i.e. $X_i \rightarrow T_i(X_i)$ for  some 
sets of possible transformations ${\cal F}_1,\dots, {\cal F}_{p-1}$.
Furthermore he assumes that the $X_i$'s are transformed independently of one another; no structural constraints are imposed.  Luce assumes the existence of a function $D$ such that 
\begin{equation}\label{eq:ustar}
 u^*\left(T_1(X_1),\dots, T_{p-1}(X_{p-1})\right)  = D(T_1,~\dots,~T_{p-1}) ~ 
u^*(X_1,\dots, X_{p-1})	
\end{equation}
for all possible transformations and choices of $X_i,~i=1,~\dots,~ p-1$. He 
determines that under these conditions, if each $X_i,~i=1,~\dots,~ p$, lies on a ratio-scale then
\[
u^*(X_1,\dots, X_{p-1}) \propto \Pi_{i=1}^{p-1}~X_i^{\alpha_i},
\]
where the $\alpha_i$'s are nondimensional constants. This is Bridgman's result, albeit proved by Luce without assuming differentiability of $u^*$.  If on the other hand some of the $X_i$'s, $i=1,~\dots,~ p-1$,  are on a ratio-scale while others are on an interval-scale and $X_p$ is on an interval-scale, then Luce proves $u^*$ cannot exist except in the case where $p=2$ and $X_1$ is on an interval-scale. 

However, as noted by \citet{aczel1986scientific}, the assumption in Equation (\ref{eq:ustar}) of the independence of the transformations $T_i$ seems unduly strong for many situations, and weakening that assumption expands the number of possibilities for the form of $u^*$.  Further work culminated in that of \citet{paganoni1987funct}. 
%
While this work was for $X_i$'s in a general vector space, for simplicity we present it here in our context, where $X_i \in \Real$, $i=1,\ldots, p$. 
Let ${\cal X}$ and ${\cal{P}}$ be nonempty subsets of $\Real^{p-1}$ and ${\cal R}$ a  set of $(p-1)$ by $(p-1)$ real-valued matrices.  Suppose that
\begin{enumerate}
	\item $\textbf{x}  + \textbf{p}  \in {\cal X}$ for all $\textbf{x} \in {\cal X}$ and $\textbf{p} \in {\cal P}$;
	\item  the identity matrix is in ${\cal R}$ and, for all $R \in {\cal R}$ and all  ${\textbf x} \in {\cal X}$,  $R {\textbf x} \in {\cal X}$;
	\item if ${\cal P} \neq \{ 0 \}$,   then $\lambda R \in {\cal R}$ for all $R \in {\cal R}$ and all $\lambda > 0$.
\end{enumerate}
Suppose also that the function $u^*$ in Equation (\ref{eq:buckingham2}) satisfies
\[
u^*(\textbf{R}~\textbf{x}+ \textbf{p}) = \alpha(\textbf{R},\textbf{p}) ~u^*(\textbf{x}) ~+ ~\beta(\textbf{R},\textbf{p})
 \label{eq:functeqn}  
\]
for all  $\textbf{R} \in \mathcal{R}$, $\textbf{x} \in {\cal  X} $ and $\textbf{p} \in {\cal P}$  
for some 
positive-valued function $\alpha$ and real-valued function $\beta$.  
Paganoni then determines the possible forms of $\alpha$ and $\beta$.

\section{Statistical invariance}\label{sec:invariance}

Having covered some important issues at the foundations of modelling in previous sections, we now turn to the modelling itself. It usually starts with a study question   on the relationship among a set of 
 specified, observable or measurable attributes of members of a population, $ \omega \in \Omega$.  A random sample of its members is to be collected to address the study question.   

A fundamental principle (Principle 1) for constructing a model for natural phenomena, which is embraced by Buckingham's Pi-theory, asserts that the model cannot depend on the scales and consequent units in which the attributes are to be measured. That principle can be extended to cover other features deemed to be irrelevant.  

Principle 2 for constructing a model calls for judicious attribute choices and transformations to reduce  the sample size needed to fit the model. The specified attributes could be design variables selected in advance of sampling to maximize the value of the study.  A notable example comes from computationally expensive computer simulators. These are run at a selected set of input attributes to develop computationally cheap emulators. 
  This in turn leads to a need to reduce the number of inputs in the predictive model that would need to be fitted empirically. 
A classical example is given in Appendix \ref{sec:reynolds} where Buckingham's theorem leads to a model with just a single predictand and a single predictor, the latter being derived from the original set 
of five predictors.

Finally we dichotomize general approaches to modelling. Approach 1 i.e. scientific modelling, leads to  what 
\citet{meinsma2019dimensional} calls  ``physical models.'' The models are generally 
deterministic with attributes measured on a ratio-scale. The $u$ in Equation \eqref{eq:buckingham}
is commonly known at least up to unknown parameters (e.g., Example 
\ref{ex:newton}),  
before any sampling is done.

Approach 2 leads to a second type of models commonly seen in the social and medical sciences. There, we have a sample of $n$ attribute-vectors, each of dimension $p$ to which the invariance principle is applied. 
And that application can lead to a nondimensionalization of the data with a consequent reduction in the number of attributes, all based on the aggregated sample of $n$ attribute-vectors. But beyond eliminating irrelevant units of measurement, applying the principle can eliminate other irrelevant features of the data, such as angle of rotation. In our approach to be described, the entire sample is holistically incorporated into model development and implementation. Now a single maximal invariant is used to summarize the sample.

In keeping with the goal of generalizing Buckingham's theory, our approach will focus on the construction of a predictive distribution model. Model uncertainties can then be characterized through such things as conditional variances and residual analysis. Furthermore principled empirical assessments of the validity of $u$ can be made given the replicate samples.

Scales play a prominent role in modelling as well. So for categorical attributes, e.g. R red, Y yellow, G green, the model should be invariant under permutations of the code by which the attributes are recorded. Models with ordinal attributes e.g. small, medium, large should be invariant under positive monotone transformation. But, as noted in Section \ref{sec:introduction}, this paper will focus mainly on ratio-scales and interval-scales. For all scales, and both approaches to modelling, transformation groups to which we now turn play a key role.

\subsection{Transformation groups}
\label{subsec:transformation}

This subsection reviews the theory of transformation groups and the statistical invariance principle, a topic that has a rich history \citep{eaton1983multivariate}. These are needed for extending the Buckingham Pi-theory. That need is recognized by \cite{shen2019statistical}, although their applications concern ratio-scales and physical models. To introduce these
groups, for simplicity, in both this section and the next, we will focus on how groups transform the sample space. 
Later, in Sections \ref{subsec:statinvariance} and \ref{sec:extendedinvariance}, we will use 
the same concepts 
for the full general theory of statistical invariance and generalized statistical invariance.

Each $\omega\in\Omega$ has a vector of measurable attributes 
$X=X(\omega)\in {\cal X}$: 
 \begin{equation}
X = \{X\}[X]=\left(
\begin{array}{c}
X_1\\
\vdots\\
X_p\\
\end{array}
\right).
\nonumber
\end{equation}
A sample of $\omega$'s is to be drawn according to a probability 
distribution $P$ on $\Omega$,
with $P$ inducing a probability 
distribution on $X$. Buckingham's theory (see Subsection \ref{subsec:buckingham}) aims at relating $\omega$'s attributes through 
a model like that in Equation (\ref{eq:buckingham2}). Our extension of that theory below will be stochastic in nature and assign $X_p$ the special role of predictand.

A sample of size $n$ yields a sample of observations $X_{ij}$ represented by 
${\bf X}^{p \times n}\in {\cal X}^n$. 
The statistical invariance principle posits that randomized statistical decision rules 
that determine  actions should be invariant under $1:1$ transformations by members $g$ of an algebraic group of transformations $G$.  That is, any pair of points ${\bf x},{\bf x}^\prime\in {\cal X}^n$ are considered equivalent for statistical inference if and only if 
${\bf x} =g({\bf x}^\prime)$ for some $g \in G$.  This equivalence is denoted ${\bf x}\sim {\bf x}^\prime$. By definition, the equivalence classes formed by $G$ are disjoint and exhaustive so we can index them by a
parameter $\gamma \in \Gamma$ and let 
${\cal X}^n_\gamma,~\gamma\in\Gamma$, 
represent an equivalence class.
The ${\cal X}^n_\gamma,~\gamma\in\Gamma$,  are referred to as orbits, which could be indexed by a set of points, $\{{\bf x}_\gamma,~\gamma\in\Gamma\}$. If the set of points satisfies some regularity conditions, then it is called a cross-section, denoted  ${\cal X}^n/G$, and its existence is studied by
\cite{wijsman1967cross}. 
Assuming a cross-section does exist, we may write
\begin{equation}
	{\cal X}^n= G \times {\cal X}^n/G.
	\nonumber
\end{equation}
In other words, any point ${\bf x} \in {\cal X}^n$ is 
represented by $(g,{\bf x}_\gamma)$ for appropriately chosen $g$ and ${\bf x}_\gamma$.



The statistical invariance principle states that a statistical decision rule must be invariant, that is, the rule must take the same value for all points in a single orbit $\gamma$.
Maximal invariant functions play a special role in statistics.  The function $M$ is invariant if its value is constant on each orbit. Further, $M$ is a maximal invariant if it takes different values on each orbit.


The following example shows the statistical invariance principle in action.

\begin{example}\label{ex:simpleexample.group}
A hard-to-make, short-lived product has an exponentially distributed random time $X~ hr$ to failure. A process-capability-analysis led to a published value of $\lambda_0$ for that product's average time-to-failure. The need to assure that standard is valid has led to a periodic sample of size $n=2$ resulting in a sample vector ${\bf X} = (X_1,X_2)$.  To make inference about $\lambda$, the expected value of $X$, following Remark \ref{rem:likelihood},  the analyst relies on the relative likelihood, i.e., ignoring irrelevant quantities  
\begin{equation}\label{eq:rellik}
\tilde{L}(\lambda) 
\doteq \frac{L(\lambda)}
{L(\lambda_0)} 	= 
-2 [\ln{\tilde{\lambda}} + \tilde{\bar{x}}\tilde{\lambda}^{-1} ] 
\end{equation}
where in general for any quantity $u$ with the same units as $\lambda_0$,  $\tilde{u} = 
u/\lambda_0$ is unitless. Differentiating the relative likelihood in Equation (\ref{eq:rellik}) yields the maximum likelihood (MLE)
\begin{equation}
	\hat{\tilde{\lambda}}_{MLE} = \tilde{\bar{x}}.
	\nonumber
\end{equation}
Using the relative likelihood thus leads to any change in $\lambda$ relative to the published value being expressed by their ratio as mandated by their lying on a ratio scale. The same is true of relative change estimated by the MLE. 

The group $G = \{g_c,~c>0\}$ transforms any realization ${\bf X}={\bf x}$ as follows:
\[
g_c ({\bf x}) = (c~x_1,c~x_2).
\]
As a maximal invariant i.e. $\pi$-function we may take, 
\begin{equation}
	\pi=M(x_1,x_2) = (x_1/x_{\cdot},x_2/x_{\cdot} )	
	= M(\tilde{x}_1, \tilde{x}_2),
	\nonumber
\end{equation}
where $x_{\cdot} = x_1 + x_2$.
The range of $M$ in $(-\infty, \infty)^2$ is given by 
\[
{\cal M} = \{(m_1, m_2): m_2 = 1-m_1,~m_1,~m_2 > 0\}. 
\]
Points in 
${\cal M}$ index the orbits of the group $G$. To locate a point ${\bf x}$ on the orbit, it entails taking
${\bf{m}}=(x_1/x_{\cdot}, x_2/x_{\cdot})$ and applying the transformation $g_c,~ c = x_{\cdot}$, to 
${\bf m}$. 
Thus, the orbits created by $G$ are rays in the positive quadrant, emanating from, but not including, the point $(0~hr,0~hr)$. Thus ${\cal X}^2$ is the union of these rays. Finally, we may let ${\cal X}^2/G = {\cal M}$. 

$M$, as a $\pi$ function, plays a key role in developing the randomized (if necessary) statistical procedures that are invariant under transformations of ${\cal X}^2$.
For example
\begin{equation}
	\hat{\tilde{\lambda}}_{MLE} = \tilde{\bar{x}}
	= c \times \upsilon[M(\tilde{x}_1, \tilde{x}_2)]
	\nonumber
\end{equation}
where $ c = \tilde{x}_{\cdot}$ and 
$\upsilon[M(\tilde{x}_1, \tilde{x}_2)] \equiv 1/2$.
But 
better choices of $\upsilon$ may be  dictated by the 
manufacturer's loss function. 

Note that here $c$, $\upsilon$, and $M$ are all unitless. 
\end{example}

Invariance of statistical procedures under the action of transformation groups may be a necessity of modelling. For instance, consider the extension of Newton's second law (Example \ref{ex:newton})  to the case of vector fields where velocity  replaces speed and direction now plays a role.  The statistical model for this extension may need to be invariant under changes of direction. In other cases, invariance may be required under permutations and monotone transformations. So in summary transformation groups may play an important role in both scientific and statistical modelling.

\subsection{Nondimensionalization}
\label{subsec:motivatingexample} 
This section presents a novel feature of this paper, the need for dimensional consistency combined with the nondimensionalization 
principle, that no model should depend on the units in which the data have been measured. 
Of particular note is the comparison  of the strict application of Buckingham's $\pi$-theory 
as described in \citet{shen2019statistical} (Approach 1) and the one we are proposing (Approach 2).  The comparison is best described in terms of an hypothetical example.


\begin{example}\label{ex:simpleexample.rain}
In a study of the 
magnitude of rainfall,  
the primary (repeating) variables  are ${ X}_1$ and ${X}_2$, denoting the depth of the rain collected in a standardized cylinder and the duration of the rainfall, respectively.  
The third quantity $X_3$ represents the magnitude of the rainfall as measured by an electronic sensor that computes a weighted average of $X_1$ as a process over the continuous time period ending at time $ X_2 $.   The dimensions of the three measurable quantities are $ [X_1] = L $, $[X_2] = T $ and $ [X_3] = LT^{-1} $, which is secondary.
Thus the attribute-vector is the column  vector
\begin{equation}
X =\left(
\begin{array}{c}
X_1\\
X_2\\
X_3
\end{array}
\right).
\nonumber
\end{equation}
The attribute-space is all possible values of $X\in {\cal X}$.

The scales and units of measurement are selected by
the investigators.  These could be changed by an arbitrary linear transformation
\begin{equation}\label{eq:transformation}
{X} \rightarrow {\bf C}^{3 \times 3}{X}   =  [\diag\{c_1,c_2,c_3\}]~ X, \quad c_i >0.
\end{equation}
 But the dimension of $X_3$  $([L]/[T])$ is related to those of $X_1$ ($[L]$), and $X_2$ ($[T]$). This relationship must be taken into account when the scales of these dimensions are selected with their associated units of measurement.

To begin, an experiment is to be performed twice and all three attributes measured each time.  
 The result will be a $3 \times 2$ dimensional matrix
 
\begin{equation}
	{\bf X} = 
\left( 
\begin{array}{ccc} 
  	X_{11} & X_{12} \\
 	 X_{21} &X_{22} \\
 	 X_{31} &X_{32} \\
 \end{array}
 \right)
 .
 \nonumber
\end{equation}
Thus, the sample space ${\cal X}^2$ will be the set of all possible realizations of ${\bf X}$. 

%
%


Now a predictive model is to be constructed in accordance with Buckingham's desideratum that the model should not depend on the measurement system the experimenters happen to choose. Furthermore, dimensional consistency dictates that any changes in the measurements must be consistently applied to all the attributes. More precisely the scales of measurement would require that $c_3 = c_1/c_2$ in the transformation matrix of Equation (\ref{eq:transformation}). 

Approach 1 focuses on $X$, not ${\bf{X}}$, and is based on considering length and time as fundamental quantities, the primary attributes are  $X_1$ and $X_2$, with respective dimensions, $L$ and $T$.
The predictand, $X_3$, 
must 
be nondimensionalized as a Buckingham  $\pi$-function. Thus,  we get 
\begin{equation}
\pi_3 = \frac{X_2 X_3}{ X_1}
\nonumber
\end{equation} 
for each one of the two sampled vectors.
The primary variables are labelled $\pi_1=\pi_2=1$. 
We define the nondimensionalized
attributes vector as 
\begin{equation}
\pi({\bf{X}}) =\left[
\begin{array}{c c c}
1  & ~~~ &1\\
1  & ~~~& 1\\
X_{21} X_{31}/ X_{11}  & & X_{22} X_{32}/{ X_{12} }
\end{array}
\right].
\nonumber
\end{equation}
 In other words, in Buckingham's theory, the function that expresses the relationship among the variables for each of the two samples is
 \[
 \pi_{3j} = u^*( \pi_{1j}, \pi_{2j}) \equiv K
 \]
 But the right hand side of this equation is a constant, which is not unreasonable since Buckingham's model was intended to be deterministic. To deal with that issue we might adopt the ad hoc solution proposed by Shen in a different example  \citep[p.~17]{shen2015dimensional} by introducing a further variable, namely a model error $\epsilon$:
\begin{equation}\label{eq:predrelation1}
\pi_{3j} = K \exp{\epsilon_{3j}}.
\end{equation}
Taking logarithms and fitting the resulting 
model yields an estimate $\hat{K}$ for 
$K$. 


In predictive form, for a 
future $\omega$, without an electronic sensor for measuring rainfall,
 Equation (\ref{eq:predrelation1}) yields, after estimating $K$,
\begin{equation}
X_{3f} = \hat{K} \frac{X_{1f}}{ X_{2f}}
\nonumber
\end{equation}
where $X_{1f}$ and 
$X_{2f}$ are the depth and duration measurements.
On the other hand, there are technical advantages to ignoring units of measurements as is commonly done in 
developing and validating statistical models, as noted by the anonymous reviewer quoted in Section \ref{sec:introduction}. In that case we would obtain 
\begin{equation}
\{X_{3f}\} = \hat{K} \frac{
\{X_{1f}\} }
{
\{X_{2f}\}}.
\nonumber
\end{equation}

\begin{remarks}\label{rem:adderror}
A more formal approach would   include the model error $X_0$ in a nondimensional form so 
that $\pi_0 = X_0$. Going through the steps above  yields 
\[
\pi_{3j} = u^*(\pi_{0j},\pi_{1j},\pi_{2j}),~j=1,2.
\]
\end{remarks}

In contrast to Approach 1, our approach, 
Approach 2, treats the sample holistically, considering the whole data matrix ${\bf{X}}$, not just $X$. We nondimensionalize
the problem by choosing as the primary  variables $X_{1j}$ and $X_{2j}$,  $j=1,2$, although 
 other  choices are available.
Let $\hat{X}_{i} =  (X_{i1}X_{i2})^{1/2}$. We then form the $\pi$-functions 
\begin{equation}
	\pi_{1j}=\frac{X_{1j}}{\hat{X}_1},\quad
\pi_{2j}=\frac{X_{2j}}{\hat{X}_2}, 
\quad
\pi_{3j}=\frac{X_{3j}\hat{X}_2}{\hat{X}_1}, \quad j=1,2.
\nonumber
\end{equation}
Then for each of the two samples we obtain
\begin{equation}
\pi_{3j}=u^*(\pi_{1j}, \pi_{2j}),~j=1,2.
\nonumber
\end{equation}
In predictive form this result becomes
\begin{equation}
X_{3j}= \frac{\hat{X}_1}{\hat{X}_2}
u^*(\pi_{1j}, \pi_{2j}),~j=1,2.
\nonumber\end{equation}
Suppose we take 
\begin{equation}\label{eq:ustar2}
u^*(\pi_{1j}, \pi_{2j})=K \frac{\pi_{1j}}{\pi_{2j}},~j=1,2
\nonumber
\end{equation}
for some positive $K$. 
Then  
\begin{equation}
X_{3j}	= \frac{\hat{X}_1}{\hat{X}_2}	 
u^*(\pi_{1j}, \pi_{2}j)
=  \frac{\hat{X}_1}{\hat{X}_2}	 
\left(K \frac{\pi_{1j}}{\pi_{2j}}\right)
=	K \frac{X_{1j}}{ X_{2j}}.\\
\nonumber \end{equation}
From the last result we obtain the model of \cite{shen2019statistical}
\begin{equation}
\pi_{3j} = K,~j=1,2.
\nonumber
\end{equation}
However, the final choice for $u^*$ could be dictated by an analysis of the data, an advantage of our 
holistic approach.

Finally we summarize our choice of $\pi$-functions as 
a maximal invariant 
\begin{equation}
	M({\bf X}) = 
\left( 
\begin{array}{cc} 
  	\pi_{11} & \pi_{12} \\
 	 \pi_{21} & \pi_{22} \\
 	 \pi_{31} & \pi_{32} \\
 \end{array}
 \right).
 \nonumber
\end{equation}


\begin{remarks}\label{rem:comparisons}
This example shows that  Approach 2 can yield the same model as Approach 1 even though Approach 1 is designed for a single $3 \times 1$ dimensional attribute vector unlike Approach 2, which starts with the entire 
$3 \times n$ sample matrix (with $n=2$). This phenomenon will be investigated in future work.
\end{remarks}

Following Example \ref{ex:simpleexample.group}, we can formalize the creation of orbits, 
$\pi$-functions and so on in terms of a transformation group
$G = \{g=g_{c_1,c_2,c_1/c_2}:~c_i > 0~[c_i],~i=1,2\}$ 
acting on the attribute-vector $X$. 
 A subgroup $G^*$ obtains by 
 restricting $c_3=c_1/c_2$.

We are now prepared to move to the general case and a generalization of the concepts seen in this example.
 
  \end{example}
  
%

\subsection{Invariant statistical models}\label{subsec:statinvariance}  
This subsection builds on Subsections 
\ref{subsec:transformation} and \ref{subsec:motivatingexample} to obtain a generalized
version of Buckingham's Pi-theory. This means transforming the scales of each of the so-called primary attributes $X_1,\dots,X_q$, which leads ineluctably to transforming the scales of the remaining,  secondary attributes $X_{q+1},\dots,X_p$. Models like that in Equation (\ref{eq:buckingham}) must reflect that link. In this subsection, we  extend Buckingham's idea beyond changes of scale  by considering the application of a 
transformation of the attribute-scales. Our approach assumes a sample of attribute vectors. When prediction is the ultimate goal of inference, as in Example (\ref{ex:newton}), our inferential aim is to construct a model as expressed in Equation (\ref{eq:buckingham2}). 



\noindent
\subsubsection*{Sample space.} We consider a sample of $n$ possibly dependent attribute-vectors collected from a sample of $\omega$'s from the population $\Omega$.  
The sample matrix, ${\bf X}^{p \times n}$, is partitioned to reflect the primary and secondary attributes as follows.  
\begin{equation}
{\bf X} = 
\left( 
\begin{array}{ccc} 
  	X_{11} & \dots &X_{1n} \\
 	 \vdots &\vdots &\vdots \\
 	 X_{p1} & \dots &X_{pn} \\
 \end{array} 
 \right)  \label{eq:sample}
\equiv
 \left( \begin{array}{c}
{\bf X}_{1}^{q \times n} \\
~ ~{\bf X}_{2}^{ (p-q)\times n}
\end{array}\right)
 \equiv
 \left( \begin{array}{c}
{\bf X}_{1} \\
{\bf X}_{2}
\end{array}\right) .
\end{equation}

Let ${\cal X}_j$  denote all possible values of  ${\bf{ X}}_j$, $j=1, 2$.  We define a group $G^*$ of  transformations on ${\cal X}_1 \times {\cal X}_2$ through the following theorem.
Each transformation is first defined on ${\cal X}_1$, with an extension to ${\cal X}_2$ that yields unit consistency. 
\begin{theorem}
Let $G_1$ be a group of transformations on ${\cal X}_1$ with identity element $e_1$. 
Assume the following.
\begin{enumerate}
    \item 
There exists a function $H$ defined on ${\cal X}_1 \times {\cal X}_2$  so that
$H({\bf X}_1, {\bf X}_2)$ is always unitless.
\item For each $g \in G_1$, there exists a  $\tilde{g}_g:  {\bf X}_2 \to {\bf X}_2$  with $ H (g ( {\bf X}_1), \tilde{g}_g({\bf X}_2)) = H({\bf X}_1, {\bf X}_2)$ for all ${\bf X}_1 \in {\cal X}_1$ and ${\bf X}_2 \in {\cal X}_2$.
\item
{\label{assump:gtildee}}
For all ${\bf X}_2 \in {\cal X}_2$,  $\tilde{g}_{e_1}({\bf X}_2)= {\bf X_2}$. 
\item
{\label{assump:composition}}
For all $g_1, g_2 \in G_1$, 
   $\tilde{g}_{g_1 \circ g_2} = \tilde{g}_{g_1} \circ \tilde{g}_{g_2}$.
 \end{enumerate}
Let $G^*$ be the set of all transformations from ${\cal X}_1 \times {\cal X}_2 $ to  ${\cal X}_1 \times {\cal X}_2 $ of the form 
\[
g^*({\bf X}_1, {\bf X}_2) = \left( g( {\bf X}_1), \tilde{g}_g( {\bf X}_2) \right),  \quad g \in G_1
.\]
  Then $G^*$ is a group under composition.
\end{theorem}

\begin{proof}
To show that $G^*$ is closed under composition, let $g_1^*$ and $g_2^*$, both in $G^*$, be associated with, respectively, $g_1 $and $g_2$, both in $ G_1$.  Then
\begin{eqnarray}
(g_1^* \circ g_2^*) ( {\bf X}_1, {\bf X}_2) &=&
g_1^* \left( g_2^* ( {\bf X}_1, {\bf X}_2)
\right) =
g_1^*\left( g_2(X_1), \tilde{g}_{g_2}(X_2)\right)
\nonumber
\\
&=&
\left(g_1 ( g_2(X_1)), \tilde{g}_{g_1}(\tilde{g}_{g_2}(X_2))\right)
= \left(  (g_1 \circ g_2) ({\bf X}_1),  (\tilde{g}_{g_1} \circ \tilde{g}_{g_2}) ( {\bf X}_2) \right)
\nonumber
\\
&=&
\left(  (g_1 \circ g_2) ({\bf X}_1),  \tilde{g}_{g_1 \circ g_2} ( {\bf X}_2) \right)
\nonumber
\end{eqnarray}
by Assumption \ref{assump:composition}. 
So $g_1^* \circ g_2^*$ is associated with $g_1 \circ g_2 \in G_1$.
We easily see that $H((g_1 \circ g_2)(X_1),\tilde{g}_{g_1 \circ g_2}(X_2) = H(X_1, X_2)$, and so
 $g_1^* \circ g_2^*$ is in $G^*$.
Clearly, the identity element of $G^*$ is given by
$e({\bf X}_1, {\bf X}_2) \equiv (e_1({\bf X}_1),  \tilde{g}_{e_1}( {\bf X}_2)) $, which equals  $({\bf X}_1,   {\bf X}_2)$ 
 by the definition of $e_1$ and by Assumption \ref{assump:gtildee}.
The inverse of $g \in G^*$ is easily found:   
if $g^*({\bf X}_1,{\bf X}_2) = (g({\bf X}_1), \tilde{g}_{g}({\bf X}_2))$, then
$(g^{*})^{-1}({\bf X}_1,{\bf X}_2) = (g^{-1}({\bf X}_1), \tilde{g}_{g^{-1}}({\bf X}_2))$.
\end{proof}

Illustrating the Theorem via Example \ref{ex:simpleexample.group}, we have 
\[
{\bf X}_1 =  \left(\begin{array}{cc} 
   X_{11} ~&~  X_{12} \\
 	 X_{21} ~&~  X_{22}
 \end{array}
 	 \right)
 	 \quad {\rm{and}} \quad
{\bf X}_2 =  \left(\begin{array}{cc} 
   X_{31} ~&~  X_{32} \\
 \end{array}
 	 \right).
\]
$G_1$ has members  
\[
g_{c_1,c_2}({\bf X}_1)  
 =
\left( 
\begin{array}{cc} 
  	c_1 X_{11} ~&~ c_1 X_{12} \\
 	c_2  X_{21} ~&~ c_2 X_{22} 
 \end{array}
 \right).
\]
One choice for the function $H$ is  
\[
H({\bf X}_1, {\bf X}_2)   
=
\left( 
 	 X_{31} X_{21}/ X_{11} ,~ X_{32} X_{22}/X_{12} 
 \right) .
\]
For each $g_{c_1, c_2} \in G_1$, we see that
$\tilde{g}_{g_{c_1,c_2}} ({\bf X}_2) = (c_1/c_2) {\bf X_2}$.
We also see that $\tilde{g}_{e_1}({\bf X}_2) = {\bf X}_2.$
The set $G^*$ consists of transformations of the form
\[
 g^*_{c_1,c_2}({ \bf X}_1, {\bf X}_2)
=
\left( 
\begin{array}{cc} 
  	c_1 X_{11} ~&~ c_1 X_{12} \\
 	c_2  X_{21} ~&~ c_2 X_{22} 
 	\\
 	(c_1/c_2)~ X_{31} ~&~ (c_1/c_2)~ X_{32}
 \end{array}
 \right).
 \]
We easily see that $\tilde{g}_{g_1 \circ g_2} = \tilde{g}_{g_1} \circ \tilde{g}_{g_2}$.  
Therefore, $G^*$ is a group.

\vskip 20pt

  


 
 Thus, by the Theorem, given a group of transformations $G_1$ on the primary attributes in the sample, we can construct a group $G^*$ of transformations on all attributes and we can write 
  \[
  {\cal X}_1 \times {\cal X}_2 = G^* \times 
  ( {\cal X}_1 \times {\cal X}_2)/G^*.
\]
  Orbits will be indexed by $\gamma\in\Gamma$ and 
  $\pi({\bf X})$ will denote a maximal invariant under the action of $G^*$.
  Let 
  \begin{equation}
  \pi({\bf X}) =  (\pi_{ij}({\bf X}))^{p \times n}.
  \nonumber
  \end{equation}

 Therefore, by the statistical invariance principle,  acceptable randomized decision rules, which include equivariant estimators as a special case, depend on ${\bf X}$ only through the maximal invariant. 
  We obtain as a special case the Buckingham $\pi$-functions as a special case where, in particular, the attributes are assessed on ratio-scales. Note that the 
$ \pi$-functions obtained in this way are not unique.    


\subsubsection*{The maximal invariant's distribution.}


Suppose that the distribution of ${\bf X}$ is  in the collection of probability distributions,  
${\cal P} = \{P_{n,\boldsymbol{\lambda}},~\boldsymbol{\lambda} \in {\Lambda}\}$. Assume, for all $g \in G^*$, the distribution of $g({\bf X})$ is also contained in ${\cal P}$.  More precisely assume that for each $g \in G^*$,  there is a one-to-one transformation $\bar{g}$ 
of ${\Lambda} $ onto ${\Lambda}$ such that ${\bf X} $ has distribution $P_{n,\boldsymbol{\lambda}}$ if and only if $g({\bf X})$ has distribution $P_{n,\bar{g}(\boldsymbol{\lambda})}$. 
 Assume further that the set $ \bar{G}^*$ of all $\bar{g}$ is a transformation group under composition, with identity  $\bar{e}$.  Assume also that $\bar{G}^*$ is homomorphic to $G^*$, i.e.~that there exists a one-to-one mapping $h $ from $ G^* $ onto $\bar{G}^*$ such that, for all $g, g^* \in G^*$, $h(g \circ g^*) = h(g) \circ h(g^*)$; $h(e) =\bar{e}$,  and $h(g^{-1}) = \{h(g)\}^{-1}$.

Let $ \pi^{-1} $ denote the set inverse, that is, $\pi^{-1}(C) = \{{\textbf X} \in {\cal X}^n$ : $\pi({\textbf{X}}) \in C\}$.  Then since $\pi({\textbf{X}})=\pi(g({\textbf{X}}))$ for any $g \in G^*$, for all  $g \in G^*$ and 
$\boldsymbol{ \lambda} \in { \Lambda} $, 
\begin{eqnarray*}
	 P_{n,\boldsymbol{ \lambda}} [ \pi({\textbf{X}})\in B ] & = &  P_{n,\boldsymbol{ \lambda}} [ \pi(g({\textbf{X}}) ) \in B ] \\
	 & = & P_{n,\boldsymbol{ \lambda}} [ g(\textbf{X}) \in \pi^{-1}(B) ] \\
	 & = &  P_{n,\bar{g} (\boldsymbol{ \lambda})}[\textbf{X} \in \pi^{-1}(B)  ] \\
	 & = & P_{n,\bar{g}( \boldsymbol{ \lambda}) } [ \pi({\textbf{X}})\in B  ].
\end{eqnarray*}
Thus, any $\boldsymbol{\lambda}^*$ ``connected to'' ${\boldsymbol{ \lambda}}$ via some $\bar{g} \in G_0$ 
induces the same distribution on $\pi({\textbf{X}})$.
This implies that $\upsilon (\boldsymbol{ \lambda} ) 
 \doteq  P_{\boldsymbol{ \lambda}} [ \pi({\textbf{X} }) \in B ] $ is invariant under transformations in $\bar{G}^*$ and hence that 
  $\upsilon (\boldsymbol{ \lambda} ) $ 
    depends on $\boldsymbol{ \lambda }$ 
    only through a maximal invariant on
     $\Lambda$.  We denote that maximal invariant by
     $\boldsymbol{\pi}_{\boldsymbol{ \lambda}} $. 
         Finally we relabel the distribution of  $ \pi({\textbf{X} })$ under ${\boldsymbol{ \lambda}}$ (and  under all of the associated $ \boldsymbol{ \lambda}^*$'s by  $P_{\boldsymbol{\pi}_{\boldsymbol{ \lambda}}}$.

The actions of the group $ G^* $  have nondimensionalized $ {\bf X}$ as $ {\bf X} \rightarrow  \pi(\textbf{X})$. Thus we obtain a stochastic version of the Pi-theorem. More precisely using the general notation $ [ {\bf U }] $  to represent ``the distribution of'' for any random object $ {\bf U }$ we have a nondimensionalized conditional distribution of the nondimensionalized predictand from sample $j$
given the transformed predictors of all samples
as
\begin{equation}\label{eq:shenassum2}
[ \pi_{pj} \mid  \pi_{1:(p-1),1:n},~ \boldsymbol{\pi}_{\boldsymbol{\lambda}}].
\end{equation}
More specifically, we have derived the result seen in Equation (\ref{eq:shenassum2}), which is the conditional distribution 
assumed in a special case by \citet{shen2015dimensional} in his Assumption 2.
Furthermore we predict for $ X_p $ by its conditional expectation, using the distribution in equation (\ref{eq:shenassum2}), which
 can be derived once the joint distribution of the attributes has been 
specified.
  The conditional variance would express the predictor's uncertainty. 
Hence statistical invariance implies that information about the variables can be summarized by maximal invariants in the sample space and in the parameter space.

\subsection{interval-scales}\label{subsec:interval}

Returning to Equation (\ref{eq:buckingham}), recall that underlying the Buckingham Pi-theorem are $ p $ variables that together describe a natural phenomenon through the relationship expressed in that equation.  The Pi-theorem assumes that $q$ of these variables are designated as the repeating or primary variables, while the remainder, which are secondary, have scales of measurement that involve the dimensions of the primary variables.  It is the latter that are converted to the $\pi$-functions in the Buckingham theorem.  But as we have seen in Subsection \ref{subsec:statinvariance}, it is these same variables that together yield the maximal invariant under the actions of a suitably chosen group, which  was fairly easily identified in the case of ratio-scales.

Subsection \ref{subsec:statinvariance} provides the bridge between the statistical invariance principle and the deterministic modeling theories described in Section \ref{sec:relationships} (i.e., the deterministic modeling frameworks developed in the physical sciences where ratio--scales are appropriate). 
Appendix \ref{app.interval} develops a similar bridge with such models in the social sciences where quantities on interval-scales are involved. For such quantities, allowable transformations extend from simple scale transformations to affine transformations. Examples of such quantities can be found in \citet{kovera2010encyclopedia}: the Intelligence Quotient (IQ), Scholastic Assessment Test (SAT), Graduate Record Examination (GRE), Graduate Management Admission Test (GMAT), and Miller Analogies Test (MAT).  Models for such quantities might involve variables measured on a ratio--scale as well. Since much of the development parallels that in Subsection \ref{subsec:statinvariance}, we omit a lot of the details for brevity and those that we do provide are in Appendix \ref{app.interval}.

\subsection{Extending invariance to random effects and the Bayesian paradigm}\label{sec:extendedinvariance}

This section extends the previous sections to
incorporate random effects and the Bayesian paradigm. Its foundations 
lie in statistical decision theory as sketched in Appendix 
\ref{sec:modellingfoundations}. Here 
the action that is a component of decision theory is 
 prediction based on a prediction model as in Equation 
(\ref{eq:buckingham2}). A training sample  of $n$ attribute vectors of length $p$ provides data for building the prediction model. Thus the predictors and predictand 
are observed for each of $n$ sampled $\omega$'s to yield the random sample's $p \times n$ matrix ${\bf X}$ seen in Equation (\ref{eq:sample}), which we denote ${\bf X}^{training}$.
Given a future ($n+1)$st attribute $p$-vector $X^{\rm{\it{future}}}$, the goal is the prediction 
of its $p$th component, $X^{\rm{\it{future}}}_p$, based on observations of its first $p-1$ components $X^{\rm{\it{future}}}_{1:(p-1)}$, all within a Bayesian framework with an appropriate extension of the framework presented in earlier sections. The situation is the one confronting the analyst who must fit a regression model based on $n$ data points and then predict a response  given only the future predictors.
We let ${\bf{X}}^{sample}$ denote the current data matrix ${\bf{X}}^{training}$ and the future data vector $X^{\rm{\it{future}}}$, with ${\bf{X}}^{sample}$, a $p \times (n+1)$ matrix in ${\cal{X}}$.


The sampling distribution of ${\bf X}^{sample}$ is determined conditional on the random parameters $\boldsymbol{\lambda}\in \Lambda$. That means specifying $\boldsymbol{\lambda}$'s prior distribution, which in turn is conditional on the set of (specified) hyperparameters $\boldsymbol{\phi}\in\Phi$.

To extend the invariance principle requires, in addition to the structures described above, 
an action space, ${\cal A}$, that is the space of possible future 
predictions of the missing observation, a prior distribution on the parameter space and a loss function, 
which remain to be specified.  We also require the specified transformation groups for ${\cal X}$ and $\Lambda$, in addition to the   transformation groups for ${\cal A}$ and $\Phi$. 
In summary, we have  the homomorphically related transformation groups $G^*,\bar{G},\hat{G}$ and $\tilde{G}$ acting on, respectively, 
${\cal X}, \Lambda, {\cal A}$ and $ \Phi$. The extended invariance principle then reduces points in these four spaces to their maximal invariants i.e. $\pi$-functions, that can be used to index the orbits induced by their respective groups. Assuming a convex loss function, the Bayes predictor in this reduced problem is a nonrandomized decision rule leading to an action in ${\cal A}$.  Each of the spaces ${\cal X}, \Lambda, {\cal A}, \Phi$ can (subject to regularity conditions)  be represented in the form 
\[
W = H_g \times W/H_g
\]
for the appropriate transformation group $H_g$
\citep{zidek1969representation}.  The corresponding maximal invariants can be expressed as matrices:
\begin{equation}
	 \pi^{{sample}}, ~ \pi^{{parameter}}, ~
	\pi^{{action}} ~
{\rm{and}}~	\pi^{{hyperparameter}}.
\nonumber
\end{equation} 
Finally using square brackets to represent  the distributions involved, we get the predictive distribution of interest conditional on quantities we know:
\begin{equation}\label{eq:bayespredictor}
[\pi_{p ,n+1}^{sample}  ~
\mid ~
\pi_{1:(p-1) , n+1}^{sample},  ~~
\pi_{1:p ,1:n}^{ sample},  ~~
\pi^{hyperparameter}
].
\end{equation}

To fix ideas we sketch an application in the following example, where we take advantage of the sufficiency and ancillarity principles to simplify the construction of the principle.


\begin{example}\label{ex:bda}
Assume the vector of observable attributes,  
$ X^{5 \times 1}$,
is  normally distributed, conditional on the mean $\mu$ and covariance matrix $\Sigma$.  We will sometimes parameterize $\Sigma$ in terms of  the diagonal matrix of standard deviations 
${\bf{\sigma}} = \diag\{\sigma_1,\dots,\sigma_5\}$ and the correlation matrix ${\bf{\rho}}$, with
 $\Sigma = {\bf{\sigma}}{\bf{\rho}} {\bf{\sigma}}$.  
 Therefore, the parameters are $\boldsymbol{\lambda}=\{ \mu, \rho, \sigma\}$ and, conditional on ${\boldsymbol{\lambda}}$,  
$X^{5 \times 1}\sim 
N_5(\lambda, \sigma {\rho} \sigma)$.
In practice, $X_5$, the fifth of these observable attributes is difficult to assess, leading to the idea of making $X_5$ a predictand and the remaining four attributes  predictors. All attributes lie on an interval-scale so a conventional approach would be seem to be straightforward:  multivariate regression analysis. Simply collect a training sample of $n$ independent vectors and fit a regression model for the intended purpose  

Complications arise due to the varying dimensions on which these attributes are to be measured. That in turn leads to different scales and different units of measurement, depending on how they are to be measured. That would not pose a problem for the unconscious statistician who might simply ignore the units. A better approach would be that suggested in Farawell (\citet{faraway2015linear}, p.103), namely rescale the measurements in a thoughtful way to eliminate those units (see also Section \ref{sec:discussion}). However, neither of those approaches deals with the rigid structural issue imposed by the need for dimensional consistency. That is, the units of measurement for the $X_i$'s are respectively   
$u_i,~i=1,...,5$, with $u_4$ and $u_5$ constrained to be $ u_2 u_1^{-1} $ and 
$u_3 u_2^{-1}$, respectively. 
To overcome the problem, Buckingham's Pi-theorem suggests itself. Thus we might use $X_{1:3}$ as primary variables to nondimensionalize $X_{4:5}$. But that does not work either since our variables lie on interval scales with $0$'s as conceptually possible values in the appropriate units. That is, $\pi$-functions as simple ratios of these variables, cannot be constructed directly to be used to nondimensionalize the attribute measurements. So ultimately we turn to the statistical invariance principle to solve the problem. The relevant transformation groups are described in what follows.

The first step creates the training set of $n$ random vectors of attribute measurements, recorded in the $5 \times n$ matrix ${\bf{X}}^{training}$.  
Letting $x_{\cdot j}$ denote the $j$th column of a realization of ${\bf{X}}^{training}$, $\bar{x}= \sum_{j=1}^n x_{\cdot j}$ and ${\bf S} = 
\sum_{j=1}^n 	(x_{\cdot j} - \bar{x})
				(x_{\cdot j} - \bar{x})^T$, the sample sum of squares,  
the likelihood function is   
\begin{eqnarray}\nonumber
L(\boldsymbol{\lambda}) &\propto &   
\mid \Sigma \mid^{-n/2}
\exp 
\left(
-\sum_{j=1}^n 
(x_{\cdot,j} - \mu)^T \Sigma^{-1} (x_{\cdot,j} - \mu)/2
\right)	\\\nonumber
\label{eq:likelihood}
& = &
\mid \Sigma \mid^{-n/2}
\exp {
- 
tr~ (\Sigma^{-1} )~{\bf S}/2
} \\ \nonumber
 & & \times 
\exp{ - 
(\mu - \bar{x})^T(\Sigma/n)^{-1}( \mu - \bar{x})^T/2).
}	 \nonumber
\end{eqnarray}
Conditional on $\mu$ and $\Sigma$, we may invoke the sufficiency principle and replace the training matrix ${\bf{X}}^{training}$ with its sufficient statistics
\begin{equation}\label{eq:suffstats}
{\bf X}^{\rm{\it{suff}}}= ( \bar{X}^{5 \times 1},~~ {\bf S}^{5 \times 5}  ),
\end{equation}
i.e., the matrix whose first column consists of the sample (row) means of ${\bf X}^{training}$
and the last five columns contain ${\bf S}$. Thus we may estimate the covariance by $\hat{\Sigma}$,  factored as 
\begin{equation}
\hat{\Sigma} = \hat{\sigma} \hat{\rho} \hat{\sigma}.
\nonumber
\end{equation}
Here $\hat{\sigma}$ denotes the diagonal matrix of estimates of the population standard deviations of the 
five attributes. Furthermore 
$\hat{\rho}$ denotes the estimate of the matrix of 
correlations between the random attributes. It is invariant under changes of scale and transformations of their origins.  
Furthermore, these quantities would be independent given the parameters. 

 Turning to the Bayesian layer, we will  adopt a conjugate prior \citep{gelman2014bayesian} for the illustrative purpose of this example, with  hyperparameters $\boldsymbol{\phi} =\{ \mu_H, \Sigma_H, B\}$:
 \begin{eqnarray}
\label{eq:priormean}
[~\mu ~ \mid ~ \Sigma~] & = & N_p (\mu_H, \Sigma/B)
\\ \label{eq:priorcovariance}
[~\Sigma ~] & = &  
Inv{\text{-}}Wishart_{B}
((B \Sigma_H)^{-1}) .\\\nonumber
\end{eqnarray}
We specify the hyperparameters 
by equating prior knowledge with a hypothetical sample of $\omega$'s and their associated attribute vectors 
 We will add a superscript $h$ on quantities below to indicate their hypothetical nature. Thus the hypothetical sample is of size $n^h$ with a likelihood  
 derived from a prior sample with $p\times n^h$ matrix
${\bf X}^h$, with sample mean $\bar{x}^h$ and sample sum of squares ${\bf{S}}^h$. 
Thus 
we obtain the 
hypothetical likelihood for $\mu$ and $\Sigma$ given the independence 
of ${\bf S}^h$ and $\bar{x}^h$:
\begin{eqnarray*}
L_{n^h}(\mu, \Sigma) &\propto &   
\mid \Sigma \mid^{-n^h/2}
\exp {
- tr~ \bigg( \Sigma^{-1} 
 {\bf S}^{h}/2 \bigg)
 }\\ \label{eq:quasilikelihood}
& & \times 
\exp{ 
- 
(\mu - \bar{x}^h)^T
(\Sigma/n_h)^{-1}
(\mu - \bar{x}^h)/2
}.	\\ \nonumber
\end{eqnarray*}

Finally complement the hypothetical 
likelihood with a noninformative improper prior on $\Sigma$
with density $\mid \Sigma\mid^{(d+1)/2}$ where $d > p$ and obtain the specification of the prior distribution.
We take the hyperparameters for the prior distributions in 
Equations (\ref{eq:priormean}) and (\ref{eq:priorcovariance}) to get  
$\mu_H=  \bar{x}^h$, $\Sigma_H =  {\bf{S}}^h/B$  and $B=n^h$.  That completes the construction of the prior.

For our prediction problem involving a future $X^{\rm{\it{future}}}$, we
 will use the posterior distribution of $\mu$ and $\Sigma$ based on the training data ${\bf{X}}^{training}$ via the sufficient statistics ${\bfX}^{ \rm{\it{suff}}}$.
To get the posterior distributions for 
$\mu$ and $\Sigma$ entails taking the product of 
the prior density as determined above based on the hypothetical sample, with the actual likelihood.

This determines a posterior density 
of ${\boldsymbol{\lambda}}$. 
Thus the covariance matrix associated with the posterior distribution of $\mu$ would depend on both the sample's 
sum of squares and the hypothetical sample's sum of squares matrix as 
\begin{equation}
	{\bf S}^{posterior} =
{\bf S}  + {\bf S}^{h} + 
\frac{
(\bar{x} - \bar{x}^h)	
(\bar{x} - \bar{x}^h)^T
}
{
1/n + 1/n^h
}.
\nonumber
\end{equation}
In other words, $S^{posterior}$ would 
replace $S_h$ to get us from its inverted Wishart prior distribution of $\Sigma$ to its posteriori distribution. Moreover, its degrees of freedom would move to $n + n_h$ to reflect the 
larger sample size. 
We omit further details since our primary interest lies in the reduced model obtained by applying the invariance principle, a reduction to which we now turn.

 We now turn to the construction of the transformation groups 
 needed to implement the invariance principle.     
 All transformations of data derive from transforming an attribute vector $X$ as follows:
 \begin{equation}
x \rightarrow    {\bf C}~(x + b) 
\nonumber
 \end{equation}
 for any possible realization $x$ of 
 $X$, where the coordinates of $b^{5 \times 1}$ have the same units as the corresponding coordinates of $x$, 
${\bf C} = \diag\{c_1,c_2,c_3, c_4^*,c_5^*\}$
where $c_4^* = c_2 c_1^{-1}$  and 
$c_5^* = c_3 c_2^{-1}$ to ensure dimensional consistency.
 The same ${\bf{C}}$ and $b$ are used to transform all data vectors and so the transformation 
of the sufficiency-reduced  matrix ${\bf{X}}^{\rm{\it{suff}}}$ in Equation (\ref{eq:suffstats}) and  $X^{future}$
is as follows:
\begin{equation}\label{eq:g}
g_{b,{\bf C}}(\bar{x}, {\bf s}, x^{\rm{\it{future}}})
= \left({\bf C}~(\bar{x} + b), ~{\bf C}~{\bf s}~{\bf C},~{\bf C}~(x^{\rm{\it{future}}} + b) \right).
\end{equation} 
We can study orbits and maximal invariants by considering the transformation of ${\bf{X}}^{\rm{\it{suff}}}$ separately from the transformation of $X^{\rm{\it{future}}}$.

Consider the decomposition of ${\bf{X}}^{\rm{\it{suff}}}$ into orbits. To index those orbits we first determine a maximal invariant.
We do this in two ways, to make a point:  first, we ignore dimensional consistency and then we include it.
 To begin, for both ways we set $b = -\bar{x} $ to transform  $\bar{x}$ to 
$ {\bf 0}~[\bar{x}_n]$, where $[\bar{x}_n]$ denotes the vector of units of $\bar{x}_n$. This means that all the points in an orbit can be reached from the new origin by choosing $b$ to be the sample average. Next we observe that we may estimate the population 
covariance for the attributes 
vector by 
${\bf s}/n = \hat{\Sigma} = \hat{\sigma} \hat{\rho} \hat{\sigma}$. 
So the matrix ${\bf C}$ in 
Equation (\ref{eq:g}), in effect, acts on the diagonal matrix $\hat{\sigma}$. 
If no restrictions were placed on ${\bf C}$'s fourth and fifth diagonal elements, then the orbits could simply be indexed by $\hat{\rho}$. In short the maximal invariant could be defined by $M (\bar{x}, {\bf s})
= ({\bf 0}, \hat{\rho})$.
But dimensional 
consistency does not allow that choice of ${\bf C}$. Instead it demands the structural requirement that we use the $c^*_4$ and $c^*_5$ specified above. The modified transformation would then act on 
$\hat{\sigma}$ as follows:
\begin{equation}
\hat{\sigma} \rightarrow 
\hat{\sigma}^*=\diag\{1,1,1,
\hat{\sigma}_1
\hat{\sigma}_4\hat{\sigma}_2^{-1},\hat{\sigma}_2\hat{\sigma}_5\hat{\sigma}_3^{-1} \}.
\nonumber
\end{equation}
The result would mean the changes in $\sigma_4$ and $\sigma_5$ would be cancelled out by the changes in the transformations of the first three $\sigma$'s. The maximal invariant 
would then be   dimensionless as required. It would make the maximal invariant for the  sufficiency-reduced sample space
\begin{equation}\label{eq:maxinvEx14sample}
	\pi^{\rm{\it{suff}}}
	=({\bf 0}, \hat{\sigma}^*
	\hat{\rho}\hat{\sigma}^*),
\end{equation}
not $({\bf 0}, \hat{\rho})$.

\begin{remarks}\label{eq:buckinghamrevisited}
The Buckingham theory concerned attributes measured on a ratio-scale. Were that the case in this example, we could have used the primary and secondary attributes differently. More precisely we could have let  
$c_i = X_i^{-1},~i=1,2,3$,  
$c_4= X_1 X_4 X_2^{-1}$ and $c_5=X_2 X_5 X_3^{-1}$.
The result would eliminate the first three variables while achieving the primary objective of non-dimensionalizing the model. We have achieved this goal using the standard deviations instead. But the method suggested here could also be used for scales other than interval, a subject of 
current research.
\end{remarks}

The corresponding maximal invariant in the parameter space for 
$\Lambda = \{ (\mu,\Sigma) \}$ would be identical to that in 
Equation (\ref{eq:maxinvEx14sample}), albeit with the hats removed, to get 
\begin{equation}
	\pi^{parameter}
	=({\bf 0}, \sigma^*
	\rho\sigma^*).
	\nonumber
\end{equation}

Observe that the ratios of the $ \sigma $'s with and without hats would be unitless and hence ancillary quantities, thus independent of the sufficient statistics. Hence the maximal invariants can be constructed from them by Basu's theorem \citep{basu1958statistics}. Finally for the hyperparameter space we would obtain the analogous result: 
\begin{equation}
	\pi^{hyperparameter}
	=({\bf 0}, \sigma_0^*
	\rho_0\sigma_0^*).
	\nonumber
\end{equation}

We can now 
compute the posterior distribution 
\begin{equation}
[~\pi^{parameter}~ \mid ~\pi^{\rm{\it{suff}}},\pi^{hyperparameter}]	
\nonumber\end{equation}
but we skip the details for brevity.

We now come to the principal objective of this example, namely a model for predicting  
a future but as yet unobserved value of the predictand, $X^{future}_{5}$, based on the future covariate 
vector $X^{future}_{1:4}$.  The data in 
${\bf X}^{\rm{\it{suff}}}$ is used as the sufficiency-reduced training sample. 
As well we assume that, given the parameters of the sampling distribution, a future attribute vector $X^{future}$ is normal with mean $\mu$ and covariance matrix $\Sigma$, and is conditionally independent of ${\bf X}^{\rm{\it{suff}}}$, given $\mu$ and $\Sigma$. 


In conformity with the modelling above, which through application of the invariance principle led to the $\pi$-functions required to nondimensionalize the problem, we transform the predictand using statistics computed from the data in ${\bf X}^{\rm{\it{suff}}}$:
$(X^{\rm{\it{future}}}_{5}-\bar{X}_5)
/\hat{\sigma}_5$. Furthermore normalized in this way, the predictand becomes invariant of those population parameters. But one more step is necessary to 
ensure that we have nondimensionalized the predictand in its $\pi$-function, namely
to align the dimensions 
of the predictand with those in the predictors to ensure dimensional consistency. The result is
\begin{equation}
\label{eq:predpifunction}
\pi^{\rm{\it{future}}}_{5}
= {\hat{\sigma}_2
\hat{\sigma}_5 \over 
\hat{\sigma}_3}
{(X^{\rm{\it{future}}}_{5}-\bar{X}_5)
\over
\hat{\sigma}_5}.
\end{equation}
We would also need to 
convert the four predictor's into their
$\pi$-functions and that would be done as in 
Equation (\ref{eq:predpifunction}).
The result will be 
$\pi^{future}_{1:4,(n+1)}$.

That predictor is found using Equation 
(\ref{eq:bayespredictor}), with modified notation. It is given by 
\begin{eqnarray*}
\label{eq:predictorEx14}
& & E[~\pi^{\rm{\it{future}}}_{5}~ 
\mid
\pi^{future}_{1:4}, 
\pi^{\rm{\it{suff}}},
\pi^{hyperparameter}
]\\
&=& E\{ E[~\pi^{future}_{5}~ 
\mid
\pi^{\rm{\it{future}}}_{1:4},\pi^{parameter}] 
\mid 
\pi^{\rm{\it{future}}}_{1:4}, 
\pi^{\rm{\it{suff}}},
\pi^{hyperparameter}
\}.
\end{eqnarray*}

\end{example}

\section{Discussion}\label{sec:discussion}


Its roots in mathematical statistics along with its formalisms made dimensional analysis (DA) seem unnecessary in statistical science. In fact, \citet{shen2019statistical} seem to have written the first paper in statistical science to recognize the need to incorporate units. For example, the authors propose what they call ``physical Lebesgue measure'' that integrates the units of measurement into the usual Lebesgue measure. Yet application of Buckingham's desideratum eliminates those units. Paradoxically it does so by exploiting the units it eliminates. That is, it exploits the intrinsic structural relationships among those units that dictate how the model must be set up. This vital implicit connection is recognized in this paper and earlier, in other papers, in more specialized contexts  \citep{shen2015dimensional,shen2019statistical}. 

\begin{remarks}
	The linear model with a Gaussian stochastic structure implicitly assumes data are measured on an interval-scale. But for physical quantities on a ratio-scale, that model would at best be an approximation. \citet{shen2018conjugate} in their Section 3.3 argue in favor of using a power (meaning product rather than sum) model in this context. We give arguments in our Subsection \ref{sec:scalesinscales} to show why additive linear models are inappropriate in this context other than as an approximation to a ratio in some sense.

That said, the simplest way of nondimensionalizing a model is by dividing each coordinate of ${\bf x}$ by a known constant with the same units of measurement as the coordinate itself, thereby removing its units of measurement.   Then $k=p$ in Buckingham's notation, and $\pi_1({\bf x}) = (x_1/c_1, x_2/c_2, \ldots, x_p/c_p)$ where $c_i= \{1\}[x_i]$. This is in effect the approach used by \citet{zhai2012a} and \citet*{zhai2012b} to resolve dimensional inconsistencies in models.
It is also the approach generally implicit in regression analysis where 
e.g. 
\[
X_5 = \beta_1 X_1+ \beta_2 X_2 +  \beta_3 X_3 + X_4
\]
with $ X_1 = 1 $  being unitless and $ X_4 $ representing a combination of measurement and modelling error. The
$ \beta_i $'s play the key role of forcing the model to adhere to the principle of \DHo~when the 
$ X_i  $'s have different units of measurement. A preferable approach would be to nondimensionalize the   $ X_i  $'s themselves in some meaningful way. For example if $ X_2 $ were the air pollution level at a specific site at a specific time, it might be divided by its long-term average over the population of sites and times.  The relative sizes of the now dimensionless $ \beta_i $'s are then readily interpretable -- a  large $ \beta $ would mean the associated $ X $ contributes a lot to the overall mean effect. 
\end{remarks}

Buckingham's theory does not specify the partition of attributes into the primary and secondary categories as is needed when deriving Buckingham's $\pi$-functions. 
That topic is a current area of active research. Recently, two approaches have been 
proposed in the context of the optimal design of computer experiments. \cite{arelis2020improve} suggests using functional analysis of variance to choose the quantities that contribute most to the variation in the output of interest as the base quantities. \cite{yang2021note}, on the other hand, propose a criterion based on measurement errors and choose the quantities that best minimize the measurement errors.

That optimal design is an ideal application for the Buckingham theory. There, the true model, called a simulator, is a deterministic numerical computer model of some phenomenon. But it is computationally intensive. So a stochastic alternative called an emulator is  fitted to a sample of outputs from the simulator at a judiciously selected set of input vectors called design points, although they represent covariates in the reality the simulator purports to represent. The Buckingham Pi-theorem simplifies the model by reducing the dimension of the input vector and hence the number of design points. It also simplifies the form of the emulator in the process. That kind of application is discussed in \citet{shen2015dimensional}, \citet{arelis2020improve} and 
\citet{adragni2009sufficient}.

  
%
%
%
%

Our approach to extending Buckingham's work differs 
from that in \citet{shen2015dimensional}. Shen 
restricts the quantities to lie on ratio-scales so
he can base his theory directly on Buckingham's
Theorem. His starting point is the application of 
that theorem and the dimensionless $ \pi $- 
functions it generates. In contrast, our theory 
allows a fully general group of transformations and
arbitrary scales.  Like Buckingham, we designate 
certain dimensions such as length $ L $ as primary (or
fundamental) while the others are secondary.  We
require that a transformation of any primary 
scale must be made simultaneously to all 
scales involving that primary scale including 
secondary scales.
That requirement ensures consistency of change
across all the quantities and leads to our 
version of the $ \pi $-functions.
However, that leaves open the issue of
which variables to serve as the primary and which
 the secondary variables, a topic under active investigation.


The paper has explored the nature and possible application of DA with the aim of integrating physical and statistical modelling. The result has been an extension of  the statistical invariance principle as a way of embracing the principles that lay behind Buckingham's development of his Pi-theory.  The result is a restriction on the class of allowable models and resulting optimal statistical procedures based on those models. How does the performance of these procedures compare with the general class of unrestricted procedures? Would a minimax or Bayesian procedure in the restricted class of allowable procedures have these same performance properties if they were thrown in with the entire set of decision rules? Under certain conditions, the answer is affirmative in the minimax case  \citep{kiefer1957invariance} and in the Bayesian case
\citep{zidek1969representation}. 

  \section{Concluding remarks}\label{sec:conclusions}

This paper has given a comprehensive overview of DA and its importance in statistical modelling. 
Dimensions have long been known to lie at the foundations of deterministic modelling, with each dimension requiring the specification of a scale and each scale requiring the specification of units of measurement.  Dimensions, their scales, and the associated units of measurement lie at the heart of empirical science.  

However,  statistical scientists regularly ignore their importance. We have demonstrated with the examples presented in Section \ref{sec:uncstatistician} that ignoring scales and units of measurement can lead to results that are either wrong or meaningless. This points to the need for statistics education to incorporate some basic training on quantity calculus and the importance of scales, along with the impact at a fundamental level of transforming data. Statistics textbooks should reflect these needs. Going beyond training is the whole process of disseminating statistical research. There again the importance of these concepts should be recognized by authors and reviewers to ensure quality. 

But does any of this really matter? We assert that it does. First we have described in Example 
\ref{ex:canadianmodel}, the genesis of 
this paper, an important example of dimensional inconsistency. An application of Buckingham's theory showed that an additional quantity needs to be added to complete the model in the sense of Buckingham to make it nondimensionalizable
\citep{wong2018dimensional}. The importance of this famous model, a model which was 
subsequently revised, lay in its use in assessing the reliability of lumber as quantified in terms of its return period,  an important component in the calculation of quality standards for lumber when placed in service as a building material. Papers flowing from that discovery soon followed 
\citep{wong2018dimensional,yang2018bayesian}.
The work reported in this paper led us to a deeper level than mere dimensional consistency, namely the discovery that the units impose important intrinsic structural links among the various quantities involved in the model. These links lead in Example \ref{ex:bda} to a new version of transformation groups usually adopted in invariant Bayesian multivariate regression models. This new version requires use of a subgroup dictated by those links. 

At a still deeper level, we are confronted by structural constraints imposed by the scales. For example, the artificial origin $0[u]$,  where $[u]$ denotes the units of $0$, in the interval scale rule out use of Buckingham's Pi-theory. Furthermore it leads to a new calculus for ratio-scales, a topic under active investigation.

Further, we have shown that, surprisingly,  not all functions are candidates for use in formulating relationship among attribute variables.  Thus functions like $ g(x) = \ln(x) $ are transcendental and hence inadmissible for that role.
This eliminates from consideration in relationships not only the natural logarithm but also, for example, the hyperbolic trigonometric functions.    This knowledge would be useful to statistical scientists in developing statistical models.


On the other hand the paper reveals  an important deficiency of deterministic physical models of natural phenomena in their failure to reflect their uncertainty about these phenomena.  
An approach to doing so 
 is presented in the paper along with an extension of the classical theory to incorporate the Bayesian approach.  
  That approach to this union of the different frameworks is reached via the statistical invariance principle, yielding a generalization  of the famous  theories of Buckingham, Bridgman  and Luce. 
  

In summary, each of the two approaches to modelling, physical and statistical, has valuable aspects that can inform the other.  This paper provides the groundwork for the unification of these approaches, setting the stage for future research.


\begin{appendix}
\section{Validity of using $ \ln( x ) $ when $ x $ has units of measurement: The debate goes on.}
\label{app.debate}

Whether as a transcendental function, the function $\ln{( x )}$ may be applied to
measurements $ x $ with units of measurement has been much discussed in other scientific disciplines and we now present some of that discussion for illustrative purposes.  \citet{molyneux1991dimensions} points out that both affirmative and  negative views had been expressed on this issue. He argues in favor of a compromise, namely defining the logarithm by exploiting one of its most fundamental properties as $\ln{(X)} = \ln{(\{X\}[X])}=  \ln{(\{X\}}) + \ln{([X])}$. He finds support for his proposal by noting that the derivative of the  constant  term, $\ln([X])$, would be zero. It follows that 
\[
\frac{d\ln{(X)}}{d X} = \frac{d\ln{(\{X\})}}{d \{X\}}.
\]
To see this under his definition of the logarithm
\begin{eqnarray*}
\frac{
\ln { ( x + \Delta x ) } - \ln { ( x ) }
}
{
\Delta \{ x \} 
} 
&=& \frac{
\ln  \{x + \Delta x \} + \ln [ x] - \ln \{  x \} -  \ln [ x]  
}
{
\Delta \{ x \} 
}  \\
&=& 
\ln  \left( 1 + {\frac{\Delta  \{ x \}}{\{ x \}} } \right)^{1/ \Delta \{ x \}} \\
&\rightarrow & \frac{1}{\{x\}}
\end{eqnarray*}
where $1/\Delta \{x\}\rightarrow \infty$.
Note that his definition of the derivative differs from ours, given in Equation (\ref{eq:log.derivative}); we include the units of $x$ in the denominatorm as the derivative is taken with respect to $x$, units and all.
Furthermore \citet{molyneux1991dimensions} argues that the proposal makes explicit the  units that are sometimes hidden, pointing to the same example, Example \ref{ex:phindex}, that we have used to make the point. It is unitless because the logarithm is applied to a count, not a measurement, that count being the number of SIUs.   \citet{molyneux1991dimensions} gives other such examples. The proposal not only makes the units explicit, but on taking the antilog of the result, you get the original value of $ X $ on the raw scale with the units $ [X] $ correctly attached. 

But, in a letter to the Journal Editor \citep{mills1995dimensions},
Ian Mills quotes \citet{molyneux1991dimensions}, in which Molyneux himself  says that his  proposal ``has no meaning.'' Furthermore, Mills says he is ``inclined to agree with him.'' Furthermore Mills argues, like Bridgman, that the issue is mute since in practice since the logarithm is applied in the context of  the difference of two logarithms, leading to
$\ln{(u/v)} = \ln{ (\{u\} / \{v\} } ) $, a unitless quantity.  In the same issue of the journal, Molyneux publishes a lengthy rejoinder saying amongst other things that Mills misquoted him.

However, in so far as the authors of this paper are aware, Molyneux's proposal was not accepted by the scientific community, leaving unresolved the issue of applying the natural logarithm to a dimensional quantity.  In particular \citet{matta2010can} also rejects it in a totally different context.  	
\citet{mayumi2010dimensions} pick up on this 
discussion in a recent paper regarding dimensional analysis in economics and the frequent application of logarithmic specifications. Their approach is based on Taylor expansion arguments that show that application of  the logarithm to dimensional quantities $ X $ is fallacious since in the expansion 
\begin{equation}\label{eq:logexpansion_1}
	\ln~(1 +  X) = X  + \frac{ X^2 }{2} + \dots .
\end{equation}
the terms on the right hand side would then have different units of measurement.

\citeauthor{mayumi2010dimensions}  then go on to describe a number of findings that are erroneous due to the misapplication of the logarithm. They also cite a ``famous controversy'' between A.C. Pigou and Milton Friedman that according to the authors, revolved around dimensional homogeneity 
\citep{pigou1936marginal,arrow1961capital}, although not specifically involving the logarithm.   One of the findings criticized in  \citet{mayumi2010dimensions} is subsequently defended by 
\citet{chilarescu2012dimensions}.  But \citet{mayumi2012response} publish a rejoinder 
in which they reassert the violation of the principle of dimensional homogeneity  
in that finding and declare that the claim in \citet{chilarescu2012dimensions}  ``is completely wrong. 
So, contrary to Chilarescu and Viasu's claim, log(V/ L) or log W in \citet{arrow1961capital}
 can never be used as a scientific representation.''

Although agreeing with the conclusion that the logarithm cannot be applied to a dimensional $ X $, \citet{matta2010can} states that Taylor expansion argument above, which the authors attribute to a Wikipedia article in September 2010, is fallacious. [The Wikipedia article actually misstates the Taylor expansion as
\begin{equation}\label{eq:logexpwrong}
\ln~ X =  X + X^2/2 + \dots. 
\end{equation}
but that does not negate the thrust of their argument.]
 They argue that the Taylor expansion should be 
\begin{equation}\label{eq:Taylorexp}
g(X) = g(X_o) + (X-X_o)\frac{dg}{dX}\mid_{X_o} +\dots
\end{equation}
so that if $ X $   had units of measurement, they would cancel out.  But then the authors don't state that expansion for the logarithm. If they did, they would have had  to deal with the issue of the units of $g(x_0)=\ln(x_0)$,
term, while the remainder of the expansion is unitless (see our comments on this issue in 
Subsection \ref{sec:log.arg.units}).

\citet{baiocchi2012dimensions} points out that if the claims of \citet{mayumi2010dimensions} were valid, they would make most ``applications of statistics, economics, ... unacceptable'' for statistical inference based on the use the exponential and logarithmic transformations. Baiocchi then tries a rescue operation by arguing that the views of \citet{mayumi2010dimensions} go ``against well established theory and practice of many disciplines including ...statistics,...and that it rests on an inadequate understanding of dimensional homogeneity and the nature of empirical modeling.'' The paper invokes the dominant theory of measurement in the social sciences that the author claims makes a numerical statement meaningful if it is invariant under legitimate scale transformations of the underlying variables. That idea of meaningfulness can then be applied to the logarithmic transformation of dimensional quantities in some situations.

To explain this idea, Baiocchi first gives the following analysis involving quantities on the ratio-scale. 
Start with the model 
$\ln~X_2  = \alpha X_1 $. Let us rescale $ X_2 $ as 
$ k~ X_2 $. That leads to the need to appropriately rescale $ X_1 $ say as $ m~ X_1 $. Consequently our model becomes 	$\ln~(k~X_2) = \alpha m~X_1$ or $\ln~X_2= \hat{\alpha} x- \ln~k$ with $\hat{\alpha}= m\alpha$. But then this model for $X_2$ cannot be reduced to its original form because of its second log term. Thus the model would be considered empirically meaningless. 
		
On the other hand if $ X_2 $ were unique up to a power transformation we would get $\ln~X_2^k = \alpha(m~X_1)$ or 
$\ln~X_2 = \hat{\alpha}X_1$ with   $ \hat{\alpha}= m\alpha/k$. Therefore the model would be invariant under admissible transformations and hence empirically meaningful. So the situation is more complex than the paper of \citet{mayumi2010dimensions} would suggest.

\citet{baiocchi2012dimensions} also addresses other arguments given by \citet{mayumi2010dimensions}. In particular, he is concerned with their Taylor expansion argument 
$\ln{(1 + X)} =  X -X^2/2 +\cdots$.  They point out that for  $ 1 + X $ to make sense, the $1$ would have to have the same units as  $ X $. They use the expansion  $\ln~(X_0+ X) = \ln~X_0 + x/a - (x/a)^2/2 + \cdots $ to make the point that when $ a = 1 $ has the same units as $ X $, the expansion is valid.  However this  argument ignores the fact that in $ \ln~X_0 $,  $ X_0 $ has units so this argument seems tenuous and therefore leaves doubt about their success in discrediting the arguments in \citet{mayumi2010dimensions}. For brevity we will terminate our review of \citet{baiocchi2012dimensions} on that note. It is a lengthy paper with further discussion of the \citet{mayumi2010dimensions} arguments  
and a very lengthy bibliography of relatively recent relevant work on this issue.  

\section{Application of Buckingham's theorem and the discovery of Reynold's number}\label{sec:reynolds}

This section provides a well-known example from fluid dynamics.

\begin{example}	\label{ex:reynolds} 
The example is a model for fluid flow around a sphere for the calculation of the drag force $F$. It turns out that the model depends only on something called the coefficient of drag and on a complicated, single dimensionless number called the Reynolds number that incorporates all the relevant quantities. 

To begin, we list all the relevant quantities, the drag force ($F$), velocity ($V$), fluid density ($\rho$), viscosity ($\mu$) and sphere diameter ($D$).  We see that we have $p=5$ $X$s in the notation of Buckingham's theorem. 
We first note that the dimensions of these five quantities can be expressed in terms of the three dimensions, length ($ L $), mass ($ M$) and time ($ T $).    We treat these as the three primary dimensions  and this tells us that we need at most 
$ 5-3 = 2 $ dimensionless $\pi$ functions to define for our model. 
 
We first write down the dimensions of each of the five quantities in terms of $L$, $M$ and $T$:  
\begin{equation}
\label{eq:drag.step1}
 [ F ] = MLT^{-2};  ~~~ [ V ] = LT^{-1}; ~~~  [\rho] =ML^{-3}; ~~~  [\mu] =ML^{-1}T^{-1}; ~~~ [ D ] = L. 
 \end{equation}
 We now proceed to sequentially eliminate the dimensions $L$, $M$ and $T$ in all five equations.
 First we use $[D]=L$ to eliminate $L$.  The first four equations become
\[
[FD^{-1}] =MT^{-2} ;  ~~~       [VD^{-1}] = T^{-1}   ;  ~~~   [D^3 \rho] = M    ;  ~~~    [D\mu] = MT^{-1}.      
\]
We next eliminate $M$ via $D^3 \rho$, yielding
\[
[FD^{-1} D^{-3}\rho^{-1}] = T^{-2}    ;  ~~~  [V D]^{-1} = T^{-1}   ;  ~~~   
[  D \mu  D^{-3} \rho^{-1}] = T^{-1} ,     
\]
that is
\[
[FD^{-4}\rho^{-1}] =T^{-2}   ;  ~~~    [VD^{-1}] = T^{-1}  ;  ~~~  [\mu D^{-2}\rho^{-1}]= T^{-1} .
\]
To eliminate $T$, we could use $[VD^{-1}]$ or 
$[\mu D^{-2} \rho^{-1}]$ or even, with a bit more work, 
$[F D^{-4}\rho^{-1}]$.   We use $[VD^{-1}]$, yielding
\[
[FD^{-4}\rho^{-1}V^{-2} D^2] =1   ~~~     {\rm{and}} ~~~~      [\mu D^{-2} \rho^{-1} V^{-1} D ]= 1,
\]
that is
\[
[F D^{-2} \rho^{-1}V^{-2}] =1         ~~~      {\rm{and}}     ~~~    [\mu D^{-1}  \rho^{-1}V^{-1}] = 1.
\]
All the dimensions are now gone so we have nondimensionalized the problem and in the process found $\pi_1$ and $\pi_2$ as implied by Buckingham's theorem:
\[
\pi_1(F,V, \mu,\rho, D) =\frac{F}{\rho D^2 V^2} \hskip 20pt
{\rm{and}} \hskip 20pt
\pi_2(F,V, \mu,\rho, D) =\frac{\mu}{\rho D V} .
\]
Therefore, for some function $U$, 
\[
U \left( \frac{F}{\rho D^2 V^2} ,\frac{\mu}{\rho D V}\right) = 0.
\]
Remarkably we have also found the famous Reynolds number, $\rho DV /\mu$ \citep{friedmann1968laminar}. The Reynolds number determines the coefficient of drag, $\rho D^2 V^2/F$, according to a fundamental law of fluid mechanics.   

\vskip 10pt
If we knew $u$ in equation (\ref{eq:buckingham}) to begin with, we could track the series of  transformations starting at (\ref{eq:drag.step1}) to find $U$.  If, however, we had no specified $u$ to begin with, we could use $\pi_1$ and $\pi_2$  to determine a model, that is, to find $U$.  For instance, we could carry out  experiments, make measurements  and determine $U$ from the data.  In either case, we can use $U$ to determine the coefficient of drag from the Reynolds number and in turn calculate the drag force.
 
\end{example}

Invoking the principle of invariance enables us to embed Example \ref{ex:reynolds} in a stochastic framework using Approach 1 as follows.

\begin{example}[continues=ex:reynolds]


We continue to use the same notation. In this example, the random variable to be replicated in independent experiments is ${\bf X} = (V,\rho,\mu,D,N)\in {\cal X}=\Real^5_+$.

\subsubsection*{The sample space.} 
The creation of the transformation group and relevant subgroup follow the lines of Example \ref{ex:simpleexample.rain}. We choose $L,M,T$ as the primary dimensions. Then with ${\bf c} = (c_1,c_2,c_3)$ the corresponding group of transformations is
\[
g_{\bf c}(V,\rho,D,\mu, N) = \bigg( \frac{c_1}{c_3} V , \frac{c_2}{c_1^3}\rho, c_1 D,\frac{c_2}{c_1c_3}\mu,\frac{c_1c_2}{c_3^2}N \bigg).
\]
For indexing the cross sections of ${\cal X}$ we have the maximal invariant 
\begin{equation}\label{eq:crosssection}
M({\bf X}) = (\frac{V}{V},\frac{\rho}{\rho},\frac{D}{D},\frac{\mu}{\rho V D^2},\frac{N}{\rho V^2 D^2})  =  (1,1,1,\pi_\mu,\pi_N, )
 \end{equation}
  where $\pi_\mu = \mu/ (\rho V D)$ and $\pi_N= N/ (\rho V^2 D^2)$.  Let 
$\boldsymbol{\pi}_{\bf X} = (\pi_\mu, \pi_N ) $.
To show that $ M $ is a maximal invariant, first  observe that $M(X)$ is invariant since each term is dimensionless. Thus showing $M(X)$ is a maximal invariant reduces to finding a subgroup element for which $X^*=g_{{\bf c}^*}(X)$ when $M(X)=M(X^*)$. 
For $N$, 
\begin{align*}
g_{{\bf c}^*}(N) &= \frac{c_1^*c_2^*}{(c_3^{*})^2}N \\
&= \frac{D^*}{D} \frac{\rho^*(D^*)^3}{D^3 \rho}  \frac{D^2}{V^2}\frac{(V^*)^2}{(D^*)^2}N \\
&=  \frac{D^*}{D} \frac{\rho^*(D^*)^3}{D^3 \rho}  \frac{D^2}{V^2}\frac{(V^*)^2}{(D^*)^2} \rho D^2 V^2  \frac{N}{\rho D^2 V^2}\\
&=\rho^* (D^*)^2 (V^*)^2 \pi_N\\
&=\rho^* (D^*)^2 (V^*)^2 \pi_{N^*} \mbox{ using the assumption that $M(X)=M(X^*)$}\\
&=\rho^* (D^*)^2 (V^*)^2 \frac{N^*}{\rho^* (D^*)^2 (V^*)^2}\\
&= N^*.
\end{align*}
Similarly, we get that $\mu \rightarrow \mu^*$.  Relating these results to Shen's thesis, this is essentially his Lemma 5.5. However Shen does not derive the maximal invariant; he simply uses the Pi quantities derived from Buckingham's Pi-theorem as the maximal invariant. In contrast, for us the maximal invariant  emerges in $M(X)$ purely as an artifact of the need for dimensional consistency as expressed through the application of the invariance principle. 

Observe that all points in ${\cal X}$ obtain from the cross section in Equation (\ref{eq:crosssection}) by application of the appropriate element of the group of transformations.  To see this let us first choose $c_1 = D^{-1}$. Then we have
\[
g_{\bf c} (\textbf{x}) =  \bigg(\frac{1}{D c_3} V, c_2D^3\rho,1,\frac{c_2D}{c_3}\mu,\frac{c_2}{D c_3^2}N \bigg).
\]
Next let $c_2 = (D^3 \rho)^{-1}$ and get 
\[
g_{\bf c}(\textbf{x}) =  \bigg(\frac{1}{D c_3} V, 1,1, \frac{1}{D^3\rho c_3}\mu, \frac{1}{D^4 \rho c_3^2}N \bigg).
\]
Finally choose $c_3 = VD^{-1}$, which yields
\[
g_{\bf c}(\textbf{x}) =  \bigg(1, 1, 1,M(\textbf{x})\bigg).
\]
Inverting this transformation takes us from the cross section to the point ${\textbf{ x}}$.

\subsubsection*{The sampling distribution}

The analysis above naturally suggests the transformation group $\bar{G}$ and its cross section for the parameter space, 
$$\Lambda = \{(\lambda_V,\lambda_\rho,\lambda_D,\lambda_\mu,\lambda_N),~\lambda_i>0,~i=V,\dots,N \}$$
namely 
\[
M(\lambda) = \bigg(1,1,1, \boldsymbol{\pi}_{ \lambda }\bigg)
\]
where with $\pi_{\lambda_\mu} = \lambda_\mu/(\lambda_\rho\lambda_D \lambda_V)$ and 
$\pi_{\lambda_N} = \lambda_N /(\lambda_\rho \lambda_D^2 \lambda_V^2)$ and  $\boldsymbol{\pi}_{ \lambda } = (\pi_{\lambda_\mu},\pi_{\lambda_N})$ characterizes the maximal invariant over 
the parameter space. It follows that for any $\lambda\in \Lambda$, 
\[
\lambda = \bar g_{{\bf c}}(M(\lambda),)
\]
where 
\begin{align*}
&c_1 = 1/ \lambda_D\\
&c_2 = 1/ (\lambda_\rho \lambda_{D}^3 )\\
&c_3 = \lambda_V / \lambda_D. \\
\end{align*}

The statistical invariance implies that
\begin{equation}\label{eq:stochmodel2}
	F(X|\lambda)=P(X\leq x|\lambda) = P(g_c(X) \leq x | \bar g_c(\lambda))
\end{equation}
for any $c_i>0, ~ i=1,2,3$. Notice that
\begin{align*}
P(g_c(X) \leq x | \bar g_c(\lambda)) &= P(X \leq g_c^{-1}(x) | \bar g_c(\lambda))\\
&= F(g_c^{-1}(x) | \bar g_c(\lambda))\\
\end{align*}
Now by taking the partial derivatives with respect to the variables, we obtain
\begin{align*}
f(x | \lambda ) = f(g_c^{-1}(x) | \bar g_c(\lambda)) 
 \frac{c_3}{c_1} 
 \frac{c_1^3}{c_2}
 \frac{1}{c_1} \frac{c_1 c_3}{c_2} 
  \frac{c_3^2}{c_1c_2}
\end{align*}
 Since this must hold for any $c_i >0$, we may choose 
 $c_1 = \lambda_D, 
 c_2 = \lambda_D^3 \lambda_\rho, 
 c_3 = \lambda_D/ \lambda_V$. Then 
\begin{align*}
f(g_c^{-1}(x) | \bar g_c(\lambda)) &= f(g_c^{-1}(x) | 
\bar g_c^{-1}(\lambda))\\
&=f(\frac{V}{\lambda_V},
\frac{\rho}{\lambda_\rho},
\frac{D}{\lambda_D},
\frac{\mu}{\lambda_\rho \lambda_D \lambda_V},
\frac{N}{\lambda_\rho \lambda_D^2 \lambda_V^2}) 
| \boldsymbol{\pi}_{\lambda})\\
&=f(\frac{V}{\lambda_V},
\frac{\rho}{\lambda_\rho},
\frac{D}{\lambda_D},
\frac{\lambda_\mu}{\lambda_\rho\lambda_D \lambda_N}\frac{\mu}{\lambda_\mu},
\frac{\lambda_N}{\lambda_\rho \lambda_D^2 \lambda_V^2}\frac{N}{\lambda_N}
 | \boldsymbol{\pi}_{\lambda})\\
\end{align*}
Thus the joint PDF is proportional to
\begin{equation}\label{eq:pdfex2}
f(
{V \over \lambda_V},
{\rho \over \lambda_\rho},
{D \over \lambda_D}, 
{\pi_{\lambda_\mu} \mu \over \lambda_\mu},
{ \pi_{\lambda_N} N \over \lambda_N} | \boldsymbol{\pi}_{\lambda} ).
\end{equation}

Hence the statistical invariance implies that information about the variables can be summarized by maximal invariants in the sample space and in the parameter space.

\subsubsection*{The sample} 
Now suppose $ n $ independent experiments are performed and that they yield data $\textbf{x}_1,\dots, \textbf{x}_n$. 
Further suppose for this illustrative example, that the model in Equation (\ref{eq:stochmodel2}) and resulting likelihood derived from Equation (\ref{eq:pdfex2}), the sufficiency principle implies 
${\bf S} =  \Sigma_j {\bf x}_j = (S_V, S_\rho ,S_D, S_\mu, S_N)$ is a sufficient statistic.  Then a maximal invariant for transformation group is
\[
M(V,\rho,D,\mu, N) =  \bigg(
 \frac{V}{S_V},
 \frac{\rho}{S_\rho},
 \frac{D}{S_D},
 {\mu \over (S_\rho S_V S_{D^2})}, 
  {
 N \over (S_\rho S_{V^2} S_{D^2})
 }
 \bigg), ~ c_i > 0, ~\forall i.
\]
To see this, observe that each term is dimensionless, so $M$ is certainly invariant. Now suppose $M(V,\rho,D,\mu, N)=M(V^*,\rho^*,D^*,\mu^*,N^*)$. Then we need to show that there exists $\{c_i^*\}$ such that $(V,\rho,D,\mu, N)=g_{c_1^*,c_2^*,c_3^*}(V^*,\rho^*,D^*,\mu^*, N^*)$. These do exist and they are 
\begin{align*}
c_1^*&= \frac{S_{D^*}}{S_D},\\
c_2^*&=\frac{S_{\rho^*} S_{{D^*}^3}}{S_{\rho} S_{D^3}},\\
c_3^*&=\frac{\frac{S_{D^*}}{S_{V^*}}}{\frac{S_{D}}{S_{V}}}.
\end{align*}

We conclude our discussion of this example. Proceeding further would entail the specification of the sampling distribution and that in turn would depend on contextual details.  
\end{example}

\section{Invariance models for interval-scales}\label{app.interval}

This section develops the theory for the interval case, which parallels that for ratio-scales seen in Subsection \ref{subsec:statinvariance}.

\subsubsection*{The sample space}\label{subsec:samplespace3}
We first partition the response vector $\textbf{X}$ as in Equation (\ref{eq:sample}). These partitions correspond to the primary and secondary quantities as in the Buckingham theory, although that distinction was not made as far as we know in the Luce work and its successors. Of particular interest again is $X_p$ in the model of Equation  (\ref{eq:buckingham2}) .  The first step in our construction entails a choice of the transformation group $G^*$.  That choice will depend on the dimensions involved.  However, given that we are assuming in this subsection that quantities  line an affine space, we will in the sequel rely on \citet{paganoni1987funct} as described in Subsection \ref{subsec:alternatives} for an illustration in this subsection.

We begin with a setup more general than that of \citet{paganoni1987funct} and would include for example the discrete seven point Semantic Differential Scale (SDM).  So we we extend Equation (\ref{eq:functeqn})  as follows
\begin{eqnarray}\label{eq:pagmodel1}
	g_1(\textbf{x}_1)  &=& \textbf{R}_1 \textbf{x}_1 + \textbf{P}_1\\ \label{eq:pagmodel2}
	g_2(\textbf{x}_2)  &= & \textbf{R}_2 \textbf{x}_2 +  \textbf{P}_2.
	\end{eqnarray}
       where  now $x_p$ is the final coordinate of
        $ \textbf{x}_2$ when $p-k +1 > 1$. Note that in the univariate version of the model model proposed by  \citet{paganoni1987funct}, Equation \ref{eq:functeqn} has the vector $ {\bf x}_2 $ replaced with $   x_p $. Here both the rescaling matrix $  \textbf{R}_2  $ and the translation vector $ \textbf{P}_2 $  depends on the pair $  \textbf{R}_1  $ and $ \textbf{P}_1 $, i.e. $ \textbf{R}_2  = R( \textbf{R}_1, \textbf{P}_1 ) $  
and $   \textbf{P}_2 = P(\textbf{R}_1, \textbf{P}_1 ) $.  
Note that the ratio--scales are formally incorporated in this extension
simply by setting to $ 0 $, the relevant coordinates of $ P_1 $ and $ P_2 $.  

Conditions are needed to ensure that  $ G^* $ is a transformation group. 
For definiteness we choose $\textbf{R}_2 = M(\textbf{R}_1) $ and 
$\textbf{P}_2  =  \psi(\textbf{R}_1) $ where in general $M(\textbf{S}\textbf{R})
= M(\textbf{S})M(\textbf{R})$ and 
$ \psi(\textbf{S} \textbf{R}) = M(\textbf{S})\psi(\textbf{R}) + \psi(\textbf{S}) $.  The objects 
$\textbf{R}_1$ and $\textbf{P}_1$ lie in the subspaces described in Subsection 
\ref{subsec:alternatives} while $\textbf{R}_2$ and $\textbf{P}_2$ lie in multidimensional rather than one dimensional spaces as before.   We omit details for brevity.

 Finally we index the transformation group $ G_0 $ acting on $\textbf{x}$ by 
$ (\textbf{R}_1,\textbf{P}_1) $ 
and define that associated transformation by 
\begin{equation}\label{eq:grouptransformation.2}
	g_{ (\textbf{R}_1,\textbf{P}_1) }  (\textbf{x}) =  
[g_1(\textbf{x}_1) ,  g_2(\textbf{x}_2) ].
\end{equation}
It remains to show that in this case, $G_0$ is a transformation group and for this we need the conditions presented by \citet{paganoni1987funct}.

\begin{theorem}
	The set $G_0$ of transformations defined by Equation (\ref{eq:grouptransformation.2})
	is a transformation group acting on the sample space.\\
	\end{theorem}
	 
	 {\bf Proof.}
First we show that  $G_0$  possesses an identity transformation. This is found 
simply by  taking $\textbf{R}_1 = 
\textbf{I}_k $ and $\textbf{P}= \textbf{0}_k$ and invoking the definitions 
of $ M^* $, $ \psi $ and $A$: 
\begin{align*}
g_1(\bfx_1)&=\bfR_1  \bfx_1 + \bfP_1= \bI_k x_1 + \bzero_k = \bfx_1. \\
g_2(\bfx_2)&= \bfR_2 \bfx_2 + \bfP_2 \\
&= M(\bfR_1)\bfx_2 + \psi(\bfR_1) + A(\bfP_1)\\
&= \bfx_2 + 0 + 0\\
&= \bfx_2.
\end{align*}

Next  we show that the composition of two transformations indexed by
$(\textbf{S}_1, \textbf{Q}_1)$ and $(\textbf{R}_1,\textbf{P}_1)$ yield a transformation in $G_0$.  First we obtain $
	g_{(\textbf{R}_1,\textbf{P}_1)}(\textbf{x}) = (\textbf{x}_1^1,\textbf{x}_2^1) $ where 
	\begin{eqnarray}\label{eq:thm1.6}
	\textbf{x}_1^1 &=&\textbf{R}_1 \textbf{x}_1+\textbf{P}_1,~ {\rm  and}\\\label{eq:thm1.7}
	\textbf{x}_2^1 &= & \textbf{R}_2 \textbf{x}_2 +\textbf{P}_2 = M(\textbf{R}_1) \textbf{x}_2 + \psi(\textbf{R}_1) + A(\bfP_1).\\ \nonumber
	\end{eqnarray}
	
Next we compute 
\begin{eqnarray*}
g_{(\textbf{S}_1,\textbf{Q}_1)} (\textbf{x}^1) &=&   \textbf{(S}_1 \textbf{x}_1^1+\textbf{Q}_1, \textbf{S}_2 \textbf{x}_2^2 +\textbf{Q}_2)  \\
&=&    \textbf{(S}_1 \textbf{x}_1^1+\textbf{Q}_1,  M( \textbf{S}_1) \textbf{x}_2^1 + \psi( \textbf{S}_1)+A(\bfQ_1)).\\
 \end{eqnarray*}
But 
\begin{align*}
M( \textbf{S}_1) \textbf{x}_2^1 + \psi( \textbf{S}_1) + A(\bfQ_1)&= M(\textbf{S}_1\textbf{R}_1) \textbf{x}_2 +
M(\textbf{S}_1)\psi(\textbf{R}_1) + M(\bfS_1)A(\bfP_1)\\
&+\psi(\textbf{S}_1)+A(\bfQ_1) \\
&= M(\textbf{S}_1\textbf{R}_1)\textbf{x}_2 + \psi(\textbf{S}_1\textbf{R}_1)+ A(\bfS_1 \bfP_1) + A(\bfQ_1)\\
&=M(\textbf{S}_1\textbf{R}_1)\textbf{x}_2 + \psi(\textbf{S}_1\textbf{R}_1)+ A(\bfS_1 \bfP_1+\bfQ_1),
\end{align*}
which proves that the composition is an element of $G_0$.

Finally we need to show that for any member of $G_0$ indexed by 
$(\textbf{R}_1,\textbf{P}_1)$ there exists an inverse.  
Starting with the transformed quantities
in Equations (\ref{eq:thm1.6})
and (\ref{eq:thm1.7}),  let $(S_1,Q_1) = (\textbf{R}_1^{-1},-\bfR_1 \textbf{P}_1)$. 
Then we find that 
\begin{align*}
g_{(\textbf{S}_1,\textbf{Q}_1)} ( \textbf{R}_1 \textbf{x}_1+\textbf{P}_1, \textbf{R}_2 \textbf{x}_2 +\textbf{P}_2 ) &= ( \bfS_1(\bfR_1 \bfx_1 + \bfP_1)+\bfQ_1, \bfS_2(\bfR_2 \bfx_2 +\bfP_2)+ \bfQ_2)\\
&=(\bfS_1 \bfR_1 \bfx_1 + \bfS_1 \bfP_1 + \bfQ_1, \bfS_2 \bfR_2 \bfx_2 + \bfS_2 \bfP_2+ \bfQ_2).
\end{align*}
But
\begin{align*}
\bfS_1 \bfR_1 \bfx_1 + \bfS_1 \bfP_1 + \bfQ_1 &= \bfR_1^{-1} \bfR_1 \bfx_1 + \bfR_1^{-1} \bfP_1 + (-\bfR_1^{-1} \bfP_1)\\
&=\bfx_1 + 0 = \bfx_1,
\end{align*}
 and
\begin{align*}
\bfS_2 \bfR_2 \bfx_2 &+ \bfS_2 \bfP_2+ \bfQ_2 = M(\bfS_1)M(\bfR_1)\bfx_2 + M(\bfS_1)(\psi(\bfR_1)+A(\bfP_1)) + (\psi(\bfS_1)+A(\bfQ_1))\\
&=M(\bfR_1^{-1})M(\bfR_1)\bfx_2 + M(\bfR_1^{-1})\psi(\bfR_1)+ M(\bfR_1^{-1})A(\bfP_1)+ \psi(\bfR_1^{-1}) + A(-\bfR_1^{-1} \bfP_1)\\
&=M(\bfR_1^{-1}\bfR_1)\bfx_2 + \psi(\bfR_1^{-1} \bfR_1)+ A(\bfR_1^{-1}\bfP_1) + A(-\bfR_1^{-1} \bfP_1)\\
&=M(\bI_k)\bfx_2 + \psi(\bI_k)+ A(\bfR_1^{-1} \bfP_1 - \bfR_1^{-1} \bfP_1)\\
&= \bfx_2 + \bzero_{p-k+1} + A(\bzero_{k})\\
&= \bfx_2 + \bzero_{p-k+1} + \bzero_{p-k+1}\\
&= \bfx_2.
\end{align*}
Thus the transformation indexed by 
 $(\textbf{R}_1^{-1},-\textbf{P}_1)$ is the inverse of that indexed by 
 $(\textbf{R}_1,\textbf{P}_1)$. That concludes the proof that $G_0$ is a transformation group.
 
\vskip .1in

We now proceed, as outlined in Subsection \ref{subsec:statinvariance}, to find the analogues of the Pi function in Buckingham's theory, which in our extension of that theory 
are coordinates of the maximal invariant under the transformation group
 $  G_0 $.  To that end
we seek that transformation  
for which 
$ g_{1\textbf{x}_1} (\textbf{x}_1) = \textbf{P}_{10}$ 
i.e. 
$ \textbf{x}_1 =g_{1\textbf{x}_1}^{-1} ( \textbf{P}_{10}) =  \textbf{S}_{1\textbf{x}_1} \textbf{P}_{10}  + 
\textbf{Q}_{1\textbf{x}_1} $ 
for an appropriate $\textbf{S}_{1\textbf{x}_1}$ and 
$\textbf{Q}_{1\textbf{x}_1}$, 
where $\textbf{S}_{1g_{(R,P)}(\textbf{x}_1)}=R\textbf{S}_{1\textbf{x}_1}$ and $\textbf{Q}_{1g_{(R,P)}(\textbf{x}_1)}=R\textbf{Q}_{1\textbf{x}_1}+P$. It follows that $\textbf{P}_{10} = \textbf{S}_{1\textbf{x}_1}^{-1}(\textbf{x}_1 -  \textbf{Q}_{1\textbf{x}_1}  ) $
for a designated fixed origin $\textbf{P}_{10}$  in the range of $\textbf{X}_1$.  
Dimensional consistency calls for the transformation of 
$\textbf{x}_2$  by the $g_2$ that complements the 
$g_1$ found in the previous paragraph, the one indexed by $(
\textbf{S}_{1\textbf{x}_1}^{-1},~- 
\textbf{S}_{1\textbf{x}_1}^{-1 }
\textbf{Q}_{1\textbf{x}_1}
)$.
If we invoke the invariance principle, we may thus transform  $\textbf{x} = (\textbf{x}_1,\textbf{x}_2)$
to 
\begin{equation}
(\boldsymbol{\pi}_{1x},  \boldsymbol{\pi}_{2x}),
\nonumber
\end{equation}
where $\boldsymbol{\pi}_{1x} = \textbf{P}_{10}$ 
and  $\boldsymbol{\pi}_{2x} = M(\textbf{S}_{1\textbf{x}_1}^{-1}) \textbf{x}_2 + 
\psi(\textbf{S}_{1\textbf{x}_1}^{-1})$ is the maximal invariant. Certainly it is invariant. Now we need to show there exists $(S^*,Q^*)$ such that $\textbf{x}=g_{(S^*,Q^*)}(\textbf{y})$ when $(\pi_{1x},\pi{2x})=(\pi_{1y},\pi{2y})$. So suppose $(\pi_{1x},\pi{2x})=(\pi_{1y},\pi{2y})$. We claim that $\textbf{x}_1=g_{(S^*,Q^*)}\textbf{y}_1$, and hence $\textbf{x}_2=g_{(M(S^*),\psi(S^*))}\textbf{y}_2$, where $ S^*=S_{1\textbf{x}_1}S_{1\textbf{y}_1}^{-1}$ and $Q^*=-(S_{1\textbf{x}_1}S_{1\textbf{y}_1}^{-1}Q_{1\textbf{y}_1}+Q_{1\textbf{x}_1})$.

{\bf Proof.}
 Assume below that $M^{-1}(X)=M(X^{-1})$.
\begin{align*}
\pi_{1x} = \pi_{1y} &\\
\iff \Sxinv (\textbf{x}_1-\Qx)&= \Syinv (\textbf{y}_1-\Qy)  \\
\iff \textbf{x}_1&=\Sx \Syinv (\textbf{y}_1-\Qy) + \Qx  \\
&=\Sx \Syinv \textbf{y}_1 - \Sx \Syinv \Qy  + \Qx \\
&=S^*\textbf{y}_1 + Q^*\\
&=g_{(S^*,Q^*)} \textbf{y}_1\\
\pi_{2x} &= \pi_{2y} \\
\iff M(\Sxinv)\textbf{x}_2 + \psi(\Sxinv) &= M(\Syinv)\textbf{y}_2 + \psi(\Syinv)\\
\iff \textbf{x}_2 &= M^{-1}(\Sxinv)M(\Syinv)\textbf{y}_2+ M^{-1}(\Sxinv)(\psi(\Syinv)- \psi(\Sxinv))\\
&=M(\Sx)M(\Syinv)\textbf{y}_2+ M(\Sx)(\psi(\Syinv)- \psi(\Sxinv))\\
&=M(\Sx \Syinv)\textbf{y}_2 + M(\Sx)\psi(\Syinv)-M(\Sx)\psi(\Sxinv)\\
&=M(\Sx \Syinv)\textbf{y}_2+ M(\Sx)\psi(\Syinv)-(\psi(\Sx \Sxinv) + \psi(\Sx))\\
&=M(\Sx \Syinv)\textbf{y}_2+ M(\Sx)\psi(\Syinv)-(\psi(I) - \psi(\Sx))\\
&=M(\Sx \Syinv)\textbf{y}_2+ M(\Sx)\psi(\Syinv)-0 + \psi(\Sx)\\
&=M(\Sx \Syinv)\textbf{y}_2+ \psi(\Sx \Syinv)\\
&=g_{(M(S^*),\psi(S^*))}\textbf{y}_2.
\end{align*}
Thus the proof is complete.

\begin{example}
The linear regression model is one of the most famous models in statistics: 
$ y^{1 \times 1} = \beta \textbf{x}^{ (p-1) \times 1} $. \citet{shen2015dimensional} shows using dimensional analysis that this model is inappropriate when all the variables are on a ratio--scale. Instead in that case the right hand side should be the product of powers of  P--functions of the coordinates of $ \textbf{x}$. But this section shows how to handle the case where the variables are regarded as interval--valued. The Pi functions would then be combinations of the $ x $ coordinates depending on the units of measurement of $ y $ and those.  

To more specific we begin by defining
	for every $ \textbf{x}\in {\cal X} $, a
  	 $ g_{\textbf{x}} \in  G_0 $ such that 
  	$ g_{\textbf{x}} ( \textbf{x}) = (g_{1\textbf{x}}(\textbf{x}_1) , g_{2\textbf{x}}(\textbf{x}_2) )  =
  	 (\boldsymbol{\pi}_{\textbf{x}1}^{1\times k}, \boldsymbol{\pi}_{\textbf{x}2}^{1 \times (p - k + 1)})$,  where 
  	 $[\boldsymbol{\pi}_{\textbf{x}1}]= [\textbf{1}_k]$, 
    $[\boldsymbol{\pi}_{\textbf{x}2}] = [\textbf{1}_{(p - k + 1)}] $ 
  	and in general
  	   $\textbf{1}_r$ denotes the vector of dimension $r$, all of whose elements are $1$, representing generically the unit on the coordinate's scale.  
  	   For the regression example the final coordinate in $ {\bf x}_2 $ is $ x_p = Y $.  It then follows from the above analysis in the notation used there  that where $$\boldsymbol{\pi}_{1x} = \textbf{P}_{10}$$ 
and  $$\boldsymbol{\pi}_{2x} = M(\textbf{S}_{1\textbf{x}_1}^{-1}) \textbf{x}_2 + 
\psi(\textbf{S}_{1\textbf{x}_1}^{-1})$$ 
is the non--dimensionalized maximal invariant. The distribution of $ \boldsymbol{\pi}_{2X}$ then determines the nondimensionalized regression model.  We omit the details for brevity.
\end{example}

\section{Foundations of statistical modelling}\label{sec:modellingfoundations}
 
After a sample ${\bf X}^{p\times n} = {\bf x}$ is selected, a statistical inquiry is expected to lead to a decision $d({\bf x})= a$ chosen from an action space ${\cal A}$. That action may be chosen by a randomized decision rule, i.e. a probability distribution $\delta(D; {\bf x})$ for events $D\subset {\cal A}$. A nonrandomized rule $d({\bf x})$ would then correspond to a degenerate probability distribution for 
$\delta(\{d({\bf x})\};{\bf x})=1$. 
 
 The decision would be based on the loss function $L(a,\lambda)$ or rather the expected loss function called the risk function
 \[
 r(\delta, \lambda) = \int L(a,\lambda) \delta(da;{\bf x})P_\lambda(d{\bf x}).
 \]
The minimax criterion is commonly used to determine the optimal decision rule as 
\[\delta_{minimax}=\mathop{\mathrm{argmin}}_{\delta}\max_\lambda   r(\delta, \lambda).
\]
An alternative is the Bayes rule where, given a prior distribution $\Pi$ for 
$\lambda$, the Bayesian decision rule is found by minimizing 
\[
 R(\delta) = \int L(a,\lambda) \delta(da;{\bf x})P_\lambda(d{\bf x})
 \Pi(\lambda).
 \]
 That prior distribution will be usually indexed by a hyperparameter vector 
$\upsilon$ lying in a hyperparameter space 
$\Upsilon$.

\end{appendix}

\section*{Acknowledgements}
We are  indebted to Professor George Bluman of the Department of Mathematics at the University of British Columbia (UBC) for his helpful discussions on dimensional analysis. Our gratitude goes to Professor Douw Steyn of the Earth and Ocean Sciences Department, also at UBC, for introducing the second author to the Unconscious Statistician.
We also thank Yongliang (Vincent) Zhai, former Masters student of the last two authors of this paper, for contributing to Example \ref{ex:unconsciousstatistician} and for his work in his Masters thesis that inspired much of this research.
Finally we acknowledge the key role played by the Forest Products Stochastic Modelling Group centered at UBC and funded by the a combination of FPInnovations and the Natural Sciences and Engineering Research Council of Canada through a Collaborative Research and Development Grant. 
It was the work of that Group that sparked our interest in the problems addressed in this paper. 
The research of the last two authors was also supported via the Discovery Grant program of the  Natural Sciences and Engineering Research Council of Canada. 



\bibliographystyle{imsart-nameyear} 
\bibliography{bib_frontmatter_to-do-list/DAbib2021}       

\begin{thebibliography}{64}

\bibitem[\protect\citeauthoryear{Acz{\'e}l, Roberts and
  Rosenbaum}{1986}]{aczel1986scientific}
\begin{barticle}[author]
\bauthor{\bsnm{Acz{\'e}l},~\bfnm{J{\'a}nos}\binits{J.}},
  \bauthor{\bsnm{Roberts},~\bfnm{Fred~S}\binits{F.~S.}} \AND
  \bauthor{\bsnm{Rosenbaum},~\bfnm{Zangwill}\binits{Z.}}
(\byear{1986}).
\btitle{On scientific laws without dimensional constants}.
\bjournal{Journal of Mathematical Analysis and Applications}
\bvolume{119}
\bpages{389--416}.
\end{barticle}
\endbibitem

\bibitem[\protect\citeauthoryear{Adragni and
  Cook}{2009}]{adragni2009sufficient}
\begin{barticle}[author]
\bauthor{\bsnm{Adragni},~\bfnm{Kofi~P}\binits{K.~P.}} \AND
  \bauthor{\bsnm{Cook},~\bfnm{R~Dennis}\binits{R.~D.}}
(\byear{2009}).
\btitle{Sufficient dimension reduction and prediction in regression}.
\bjournal{Philosophical Transactions of the Royal Society A: Mathematical,
  Physical and Engineering Sciences}
\bvolume{367}
\bpages{4385--4405}.
\end{barticle}
\endbibitem

\bibitem[\protect\citeauthoryear{Albrecht
  et~al.}{2013}]{albrecht2013experimental}
\begin{barticle}[author]
\bauthor{\bsnm{Albrecht},~\bfnm{Mark~C}\binits{M.~C.}},
  \bauthor{\bsnm{Nachtsheim},~\bfnm{Christopher~J}\binits{C.~J.}},
  \bauthor{\bsnm{Albrecht},~\bfnm{Thomas~A}\binits{T.~A.}} \AND
  \bauthor{\bsnm{Cook},~\bfnm{R~Dennis}\binits{R.~D.}}
(\byear{2013}).
\btitle{Experimental design for engineering dimensional analysis}.
\bjournal{Technometrics}
\bvolume{55}
\bpages{257--270}.
\end{barticle}
\endbibitem

\bibitem[\protect\citeauthoryear{Arelis}{2020}]{arelis2020improve}
\begin{bphdthesis}[author]
\bauthor{\bsnm{Arelis},~\bfnm{Rodr{\'\i}guez Alexi~Gilberto}\binits{R.~A.~G.}}
(\byear{2020}).
\btitle{How to improve prediction accuracy in the analysis of computer
  experiments, exploitation of low-order effects and dimensional analysis},
\btype{PhD thesis},
\bpublisher{University of British Columbia}.
\end{bphdthesis}
\endbibitem

\bibitem[\protect\citeauthoryear{Arrow et~al.}{1961}]{arrow1961capital}
\begin{barticle}[author]
\bauthor{\bsnm{Arrow},~\bfnm{Kenneth~J}\binits{K.~J.}},
  \bauthor{\bsnm{Chenery},~\bfnm{Hollis~B}\binits{H.~B.}},
  \bauthor{\bsnm{Minhas},~\bfnm{Bagicha~S}\binits{B.~S.}} \AND
  \bauthor{\bsnm{Solow},~\bfnm{Robert~M}\binits{R.~M.}}
(\byear{1961}).
\btitle{Capital-labor substitution and economic efficiency}.
\bjournal{The Review of Economics and Statistics}
\bpages{225--250}.
\end{barticle}
\endbibitem

\bibitem[\protect\citeauthoryear{Baiocchi}{2012}]{baiocchi2012dimensions}
\begin{barticle}[author]
\bauthor{\bsnm{Baiocchi},~\bfnm{Giovanni}\binits{G.}}
(\byear{2012}).
\btitle{On dimensions of ecological economics}.
\bjournal{Ecological Economics}
\bvolume{75}
\bpages{1--9}.
\end{barticle}
\endbibitem

\bibitem[\protect\citeauthoryear{Basu}{1958}]{basu1958statistics}
\begin{barticle}[author]
\bauthor{\bsnm{Basu},~\bfnm{D}\binits{D.}}
(\byear{1958}).
\btitle{On statistics independent of sufficient statistics}.
\bjournal{Sankhy{\=a}: The Indian Journal of Statistics}
\bpages{223--226}.
\end{barticle}
\endbibitem

\bibitem[\protect\citeauthoryear{Bluman and Cole}{1974}]{bluman.cole}
\begin{bbook}[author]
\bauthor{\bsnm{Bluman},~\bfnm{G~W}\binits{G.~W.}} \AND
  \bauthor{\bsnm{Cole},~\bfnm{J~D}\binits{J.~D.}}
(\byear{1974}).
\btitle{Similarity Methods for Differential Equations}.
\bseries{Applied Mathematical Sciences}.
\bpublisher{Springer-Verlag}.
\end{bbook}
\endbibitem

\bibitem[\protect\citeauthoryear{Box and Cox}{1964}]{box1964analysis}
\begin{barticle}[author]
\bauthor{\bsnm{Box},~\bfnm{George~EP}\binits{G.~E.}} \AND
  \bauthor{\bsnm{Cox},~\bfnm{David~R}\binits{D.~R.}}
(\byear{1964}).
\btitle{An analysis of transformations}.
\bjournal{Journal of the Royal Statistical Society. Series B (Methodological)}
\bvolume{26}
\bpages{211--252}.
\end{barticle}
\endbibitem

\bibitem[\protect\citeauthoryear{Bridgman}{1931}]{bridgman1931dimensions}
\begin{bbook}[author]
\bauthor{\bsnm{Bridgman},~\bfnm{P.~W.}\binits{P.~W.}}
(\byear{1931}).
\btitle{Dimensional Analysis, Revised Edition}.
\bpublisher{Yale University Press}, \baddress{New Haven}.
\end{bbook}
\endbibitem

\bibitem[\protect\citeauthoryear{Buckingham}{1914}]{buckingham1914physically}
\begin{barticle}[author]
\bauthor{\bsnm{Buckingham},~\bfnm{Edgar}\binits{E.}}
(\byear{1914}).
\btitle{On physically similar systems; illustrations of the use of dimensional
  equations}.
\bjournal{Physical Review}
\bvolume{4}
\bpages{345--376}.
\end{barticle}
\endbibitem

\bibitem[\protect\citeauthoryear{Chilarescu and
  Viasu}{2012}]{chilarescu2012dimensions}
\begin{barticle}[author]
\bauthor{\bsnm{Chilarescu},~\bfnm{Constantin}\binits{C.}} \AND
  \bauthor{\bsnm{Viasu},~\bfnm{Ioana}\binits{I.}}
(\byear{2012}).
\btitle{Dimensions and logarithmic function in economics: A comment}.
\bjournal{Ecological Economics}
\bvolume{75}
\bpages{10--11}.
\end{barticle}
\endbibitem

\bibitem[\protect\citeauthoryear{Cohen et~al.}{2004}]{cohen2004urban}
\begin{barticle}[author]
\bauthor{\bsnm{Cohen},~\bfnm{Aaron~J}\binits{A.~J.}},
  \bauthor{\bsnm{Anderson},~\bfnm{H~Ross}\binits{H.~R.}},
  \bauthor{\bsnm{Ostro},~\bfnm{Bart}\binits{B.}},
  \bauthor{\bsnm{Pandey},~\bfnm{K~Dev}\binits{K.~D.}},
  \bauthor{\bsnm{Krzyzanowski},~\bfnm{Michal}\binits{M.}},
  \bauthor{\bsnm{K{\"u}nzli},~\bfnm{Nino}\binits{N.}},
  \bauthor{\bsnm{Gutschmidt},~\bfnm{Kersten}\binits{K.}},
  \bauthor{\bsnm{Pope~III},~\bfnm{C~Arden}\binits{C.~A.}},
  \bauthor{\bsnm{Romieu},~\bfnm{Isabelle}\binits{I.}},
  \bauthor{\bsnm{Samet},~\bfnm{Jonathan~M}\binits{J.~M.}} \betal{et~al.}
(\byear{2004}).
\btitle{Urban air pollution}.
\bjournal{Comparative quantification of health risks: global and regional
  burden of disease attributable to selected major risk factors}
\bvolume{2}
\bpages{1353--1433}.
\end{barticle}
\endbibitem

\bibitem[\protect\citeauthoryear{De~Oliveira, Kedem and
  Short}{1997}]{de1997bayesian}
\begin{barticle}[author]
\bauthor{\bsnm{De~Oliveira},~\bfnm{V.}\binits{V.}},
  \bauthor{\bsnm{Kedem},~\bfnm{B.}\binits{B.}} \AND
  \bauthor{\bsnm{Short},~\bfnm{D.~A.}\binits{D.~A.}}
(\byear{1997}).
\btitle{Bayesian prediction of transformed {G}aussian random fields}.
\bjournal{Journal of the American Statistical Association}
\bvolume{92}
\bpages{1422--1433}.
\end{barticle}
\endbibitem

\bibitem[\protect\citeauthoryear{Dou, Le and Zidek}{2007}]{dou2007dynamic}
\begin{btechreport}[author]
\bauthor{\bsnm{Dou},~\bfnm{YP}\binits{Y.}},
  \bauthor{\bsnm{Le},~\bfnm{ND}\binits{N.}} \AND
  \bauthor{\bsnm{Zidek},~\bfnm{JV}\binits{J.}}
(\byear{2007}).
\btitle{A dynamic linear model for hourly ozone concentrations}
\btype{Technical Report} No. \bnumber{228},
\bpublisher{Statistics Department, University of British Columbia}.
\end{btechreport}
\endbibitem

\bibitem[\protect\citeauthoryear{Draper and
  Cox}{1969}]{draper1969distributions}
\begin{barticle}[author]
\bauthor{\bsnm{Draper},~\bfnm{N~R}\binits{N.~R.}} \AND
  \bauthor{\bsnm{Cox},~\bfnm{D~R}\binits{D.~R.}}
(\byear{1969}).
\btitle{On distributions and their transformation to normality}.
\bjournal{Journal of the Royal Statistical Society: Series B (Methodological)}
\bvolume{31}
\bpages{472--476}.
\end{barticle}
\endbibitem

\bibitem[\protect\citeauthoryear{Eaton}{1983}]{eaton1983multivariate}
\begin{bbook}[author]
\bauthor{\bsnm{Eaton},~\bfnm{Morris~L}\binits{M.~L.}}
(\byear{1983}).
\btitle{Multivariate Statistics: a Vector Space Approach.}
\bpublisher{John Wiley \& Sons, Inc., 605 Third Ave., New York, NY 10158, USA,
  1983, 512}.
\end{bbook}
\endbibitem

\bibitem[\protect\citeauthoryear{Faraway}{2015}]{faraway2015linear}
\begin{bbook}[author]
\bauthor{\bsnm{Faraway},~\bfnm{Julian~J}\binits{J.~J.}}
(\byear{2015}).
\btitle{Linear models with R: Second Edition}.
\bpublisher{Chapman and Hall/CRC}.
\end{bbook}
\endbibitem

\bibitem[\protect\citeauthoryear{Finney}{1977}]{finney1977dimensions}
\begin{barticle}[author]
\bauthor{\bsnm{Finney},~\bfnm{DJ}\binits{D.}}
(\byear{1977}).
\btitle{Dimensions of statistics}.
\bjournal{Journal of the Royal Statistical Society, Series C (Applied
  Statistics)}
\bvolume{26}
\bpages{285--289}.
\end{barticle}
\endbibitem

\bibitem[\protect\citeauthoryear{Foschi and Yao}{1986}]{foschiyao1986dol}
\begin{binproceedings}[author]
\bauthor{\bsnm{Foschi},~\bfnm{R.~O.}\binits{R.~O.}} \AND
  \bauthor{\bsnm{Yao},~\bfnm{F.~Z.}\binits{F.~Z.}}
(\byear{1986}).
\btitle{Another look at three duration of load models}.
In \bbooktitle{Proceedings, XVII IUFRO Congress}.
\end{binproceedings}
\endbibitem

\bibitem[\protect\citeauthoryear{Fourier}{1822}]{fourier1822theorie}
\begin{bbook}[author]
\bauthor{\bsnm{Fourier},~\bfnm{Joseph}\binits{J.}}
(\byear{1822}).
\btitle{Th{\'e}orie Analytique de la Chaleur, par M. Fourier}.
\bpublisher{Chez Firmin Didot, P{\`e}re et Fils}.
\end{bbook}
\endbibitem

\bibitem[\protect\citeauthoryear{Friedmann, Gillis and
  Liron}{1968}]{friedmann1968laminar}
\begin{barticle}[author]
\bauthor{\bsnm{Friedmann},~\bfnm{M}\binits{M.}},
  \bauthor{\bsnm{Gillis},~\bfnm{J}\binits{J.}} \AND
  \bauthor{\bsnm{Liron},~\bfnm{N}\binits{N.}}
(\byear{1968}).
\btitle{Laminar flow in a pipe at low and moderate Reynolds numbers}.
\bjournal{Applied Scientific Research}
\bvolume{19}
\bpages{426--438}.
\end{barticle}
\endbibitem

\bibitem[\protect\citeauthoryear{Gelman et~al.}{2014}]{gelman2014bayesian}
\begin{bbook}[author]
\bauthor{\bsnm{Gelman},~\bfnm{Andrew}\binits{A.}},
  \bauthor{\bsnm{Carlin},~\bfnm{John~B}\binits{J.~B.}},
  \bauthor{\bsnm{Stern},~\bfnm{Hal~S}\binits{H.~S.}},
  \bauthor{\bsnm{Dunson},~\bfnm{David~B}\binits{D.~B.}},
  \bauthor{\bsnm{Vehtari},~\bfnm{Aki}\binits{A.}} \AND
  \bauthor{\bsnm{Rubin},~\bfnm{Donald~B}\binits{D.~B.}}
(\byear{2014}).
\btitle{Bayesian Data Analysis; Third edition}.
\bpublisher{CRC press}.
\end{bbook}
\endbibitem

\bibitem[\protect\citeauthoryear{Gibbings}{2011}]{gibbons2011dimensional}
\begin{bbook}[author]
\bauthor{\bsnm{Gibbings},~\bfnm{J.~C.}\binits{J.~C.}}
(\byear{2011}).
\btitle{Dimensional Analysis}.
\bpublisher{Springer Verlag}, \baddress{London}.
\end{bbook}
\endbibitem

\bibitem[\protect\citeauthoryear{Hand}{1996}]{hand1996statistics}
\begin{barticle}[author]
\bauthor{\bsnm{Hand},~\bfnm{David~J}\binits{D.~J.}}
(\byear{1996}).
\btitle{Statistics and the theory of measurement}.
\bjournal{Journal of the Royal Statistical Society. Series A (Statistics in
  Society)}
\bvolume{159}
\bpages{445--492}.
\end{barticle}
\endbibitem

\bibitem[\protect\citeauthoryear{H{\"a}rdle and
  Vogt}{2015}]{hardle2017bortkiewicz}
\begin{barticle}[author]
\bauthor{\bsnm{H{\"a}rdle},~\bfnm{Wolfgang~Karl}\binits{W.~K.}} \AND
  \bauthor{\bsnm{Vogt},~\bfnm{Annette~B.}\binits{A.~B.}}
(\byear{2015}).
\btitle{Ladislaus von {B}ortkiewicz---Statistician, Economist and a {E}uropean
  Intellectual}.
\bjournal{International Statistical Review}
\bvolume{83}
\bpages{17-35}.
\bdoi{10.1111/insr.12083}
\end{barticle}
\endbibitem

\bibitem[\protect\citeauthoryear{Hoffmeyer and
  S{\o}rensen}{2007}]{hoffmeyer2007duration}
\begin{barticle}[author]
\bauthor{\bsnm{Hoffmeyer},~\bfnm{Preben}\binits{P.}} \AND
  \bauthor{\bsnm{S{\o}rensen},~\bfnm{John~Dalsgaard}\binits{J.~D.}}
(\byear{2007}).
\btitle{Duration of load revisited}.
\bjournal{Wood Science and Technology}
\bvolume{41}
\bpages{687--711}.
\end{barticle}
\endbibitem

\bibitem[\protect\citeauthoryear{Kiefer et~al.}{1957}]{kiefer1957invariance}
\begin{barticle}[author]
\bauthor{\bsnm{Kiefer},~\bfnm{Jack}\binits{J.}} \betal{et~al.}
(\byear{1957}).
\btitle{Invariance, minimax sequential estimation, and continuous time
  processes}.
\bjournal{The Annals of Mathematical Statistics}
\bvolume{28}
\bpages{573--601}.
\end{barticle}
\endbibitem

\bibitem[\protect\citeauthoryear{K{\"o}hler and
  Svensson}{2002}]{kohler2002probabilistic}
\begin{binproceedings}[author]
\bauthor{\bsnm{K{\"o}hler},~\bfnm{Jochen}\binits{J.}} \AND
  \bauthor{\bsnm{Svensson},~\bfnm{Staffan}\binits{S.}}
(\byear{2002}).
\btitle{Probabilistic modelling of duration of load effects in timber
  structures}.
In \bbooktitle{Proceedings of the 35th Meeting, International Council for
  Research and Innovation in Building and Construction, Working Commission
  W18--Timber Structures, CIB-W18, Paper}
\bvolume{35-17}
\bpages{1}.
\end{binproceedings}
\endbibitem

\bibitem[\protect\citeauthoryear{Kovera}{2010}]{kovera2010encyclopedia}
\begin{bmisc}[author]
\bauthor{\bsnm{Kovera},~\bfnm{Margaret~Bull}\binits{M.~B.}}
(\byear{2010}).
\btitle{Encyclopedia of research design}.
\end{bmisc}
\endbibitem

\bibitem[\protect\citeauthoryear{Lehmann and Romano}{2010}]{lehmann2010testing}
\begin{bbook}[author]
\bauthor{\bsnm{Lehmann},~\bfnm{Erick~L}\binits{E.~L.}} \AND
  \bauthor{\bsnm{Romano},~\bfnm{Joseph~P}\binits{J.~P.}}
(\byear{2010}).
\btitle{Testing Statistical Hypotheses}.
\bpublisher{Springer}.
\end{bbook}
\endbibitem

\bibitem[\protect\citeauthoryear{LibreTexts}{2019}]{idealgaslaw}
\begin{bmisc}[author]
\bauthor{\bsnm{LibreTexts}}
(\byear{2019}).
\btitle{The {I}deal {G}as {L}aw}.
\bhowpublished{{\url{https://chem.libretexts.org/Bookshelves/Physical_and_Theoretical_Chemistry_Textbook_Maps/Supplemental_Modules_(Physical_and_Theoretical_Chemistry)/Physical_Properties_of_Matter/States_of_Matter/Properties_of_Gases/Gas_Laws/The_Ideal_Gas_Law}}}.
\bnote{Accessed 01/04/2020}.
\end{bmisc}
\endbibitem

\bibitem[\protect\citeauthoryear{Lin and Shen}{2013}]{lin2013comment}
\begin{barticle}[author]
\bauthor{\bsnm{Lin},~\bfnm{Dennis~KJ}\binits{D.~K.}} \AND
  \bauthor{\bsnm{Shen},~\bfnm{Weijie}\binits{W.}}
(\byear{2013}).
\btitle{Comment: some statistical concerns on dimensional analysis}.
\bjournal{Technometrics}
\bvolume{55}
\bpages{281--285}.
\end{barticle}
\endbibitem

\bibitem[\protect\citeauthoryear{Luce}{1959}]{luce1959possible}
\begin{barticle}[author]
\bauthor{\bsnm{Luce},~\bfnm{R~Duncan}\binits{R.~D.}}
(\byear{1959}).
\btitle{On the possible psychophysical laws}.
\bjournal{Psychological Review}
\bvolume{66}
\bpages{81}.
\end{barticle}
\endbibitem

\bibitem[\protect\citeauthoryear{Luce}{1964}]{luce1964generalization}
\begin{barticle}[author]
\bauthor{\bsnm{Luce},~\bfnm{R~Duncan}\binits{R.~D.}}
(\byear{1964}).
\btitle{A generalization of a theorem of dimensional analysis}.
\bjournal{Journal of Mathematical Psychology}
\bvolume{1}
\bpages{278--284}.
\end{barticle}
\endbibitem

\bibitem[\protect\citeauthoryear{Magnello}{2009}]{magnello2009karl}
\begin{barticle}[author]
\bauthor{\bsnm{Magnello},~\bfnm{M~Eileen}\binits{M.~E.}}
(\byear{2009}).
\btitle{Karl {P}earson and the establishment of mathematical statistics}.
\bjournal{International Statistical Review}
\bvolume{77}
\bpages{3--29}.
\end{barticle}
\endbibitem

\bibitem[\protect\citeauthoryear{Matta et~al.}{2010}]{matta2010can}
\begin{barticle}[author]
\bauthor{\bsnm{Matta},~\bfnm{Ch{\'e}rif~F}\binits{C.~F.}},
  \bauthor{\bsnm{Massa},~\bfnm{Lou}\binits{L.}},
  \bauthor{\bsnm{Gubskaya},~\bfnm{Anna~V}\binits{A.~V.}} \AND
  \bauthor{\bsnm{Knoll},~\bfnm{Eva}\binits{E.}}
(\byear{2010}).
\btitle{Can one take the logarithm or the sine of a dimensioned quantity or a
  unit? Dimensional analysis involving transcendental functions}.
\bjournal{Journal of Chemical Education}
\bvolume{88}
\bpages{67--70}.
\end{barticle}
\endbibitem

\bibitem[\protect\citeauthoryear{Mayumi and
  Giampietro}{2010}]{mayumi2010dimensions}
\begin{barticle}[author]
\bauthor{\bsnm{Mayumi},~\bfnm{Kozo}\binits{K.}} \AND
  \bauthor{\bsnm{Giampietro},~\bfnm{Mario}\binits{M.}}
(\byear{2010}).
\btitle{Dimensions and logarithmic function in economics: A short critical
  analysis}.
\bjournal{Ecological Economics}
\bvolume{69}
\bpages{1604--1609}.
\end{barticle}
\endbibitem

\bibitem[\protect\citeauthoryear{Mayumi and
  Giampietro}{2012}]{mayumi2012response}
\begin{barticle}[author]
\bauthor{\bsnm{Mayumi},~\bfnm{Kozo}\binits{K.}} \AND
  \bauthor{\bsnm{Giampietro},~\bfnm{Mario}\binits{M.}}
(\byear{2012}).
\btitle{Response to dimensions and logarithmic function in economics: A
  comment}.
\bjournal{Ecological Economics}
\bvolume{75}
\bpages{12--14}.
\end{barticle}
\endbibitem

\bibitem[\protect\citeauthoryear{Meinsma}{2019}]{meinsma2019dimensional}
\begin{barticle}[author]
\bauthor{\bsnm{Meinsma},~\bfnm{Gjerrit}\binits{G.}}
(\byear{2019}).
\btitle{Dimensional and scaling analysis}.
\bjournal{SIAM review}
\bvolume{61}
\bpages{159--184}.
\end{barticle}
\endbibitem

\bibitem[\protect\citeauthoryear{Mills}{1995}]{mills1995dimensions}
\begin{barticle}[author]
\bauthor{\bsnm{Mills},~\bfnm{Ian~M}\binits{I.~M.}}
(\byear{1995}).
\btitle{Dimensions of logarithmic quantities}.
\bjournal{Journal of Chemical Education}
\bvolume{72}
\bpages{954}.
\end{barticle}
\endbibitem

\bibitem[\protect\citeauthoryear{Molyneux}{1991}]{molyneux1991dimensions}
\begin{barticle}[author]
\bauthor{\bsnm{Molyneux},~\bfnm{P.}\binits{P.}}
(\byear{1991}).
\btitle{The dimensions of logarithmic quantities: implications for the hidden
  concentration and pressure units in p{H} values, acidity constants, standard
  thermodynamic functions, and standard electrode potentials}.
\bjournal{Journal of Chemical Education}
\bvolume{68}
\bpages{467}.
\end{barticle}
\endbibitem

\bibitem[\protect\citeauthoryear{Mosteller and Tukey}{1977}]{mosteller1977data}
\begin{bbook}[author]
\bauthor{\bsnm{Mosteller},~\bfnm{Frederick}\binits{F.}} \AND
  \bauthor{\bsnm{Tukey},~\bfnm{John~Wilder}\binits{J.~W.}}
(\byear{1977}).
\btitle{Data Analysis and Regression: A Second Course in Statistics.}
\bpublisher{Addison-Wesley Series in Behavioral Science: Quantitative Methods}.
\end{bbook}
\endbibitem

\bibitem[\protect\citeauthoryear{{Joint Commitee on Guides in
  Metrology}}{2012}]{jcgm2012200}
\begin{bmisc}[author]
\bauthor{\bsnm{{Joint Commitee on Guides in Metrology}}}
(\byear{2012}).
\btitle{200: 2012 --- {I}nternational {V}ocabulary of {M}etrology {B}asic and
  General Concepts and Associated Terms ({V}{I}{M})}.
\end{bmisc}
\endbibitem

\bibitem[\protect\citeauthoryear{Paganoni}{1987}]{paganoni1987funct}
\begin{barticle}[author]
\bauthor{\bsnm{Paganoni},~\bfnm{L}\binits{L.}}
(\byear{1987}).
\btitle{On a functional equation concerning affine transformations}.
\bjournal{Journal of Mathematical Analysis and Applications}
\bvolume{127}
\bpages{475--491}.
\end{barticle}
\endbibitem

\bibitem[\protect\citeauthoryear{Pigou, Friedman and
  Georgescu-Roegen}{1936}]{pigou1936marginal}
\begin{barticle}[author]
\bauthor{\bsnm{Pigou},~\bfnm{Arthur~C}\binits{A.~C.}},
  \bauthor{\bsnm{Friedman},~\bfnm{Milton}\binits{M.}} \AND
  \bauthor{\bsnm{Georgescu-Roegen},~\bfnm{Nicholas}\binits{N.}}
(\byear{1936}).
\btitle{Marginal utility of money and elasticities of demand}.
\bjournal{The Quarterly Journal of Economics}
\bvolume{50}
\bpages{532--539}.
\end{barticle}
\endbibitem

\bibitem[\protect\citeauthoryear{Shen}{2015}]{shen2015dimensional}
\begin{bphdthesis}[author]
\bauthor{\bsnm{Shen},~\bfnm{Weijie}\binits{W.}}
(\byear{2015}).
\btitle{Dimensional analysis in statistics: theories, methodologies and
  applications},
\btype{PhD thesis},
\bpublisher{The Pennsylvania State University}.
\end{bphdthesis}
\endbibitem

\bibitem[\protect\citeauthoryear{Shen and Lin}{2018}]{shen2018conjugate}
\begin{barticle}[author]
\bauthor{\bsnm{Shen},~\bfnm{Weijie}\binits{W.}} \AND
  \bauthor{\bsnm{Lin},~\bfnm{Dennis~KJ}\binits{D.~K.}}
(\byear{2018}).
\btitle{A conjugate model for dimensional analysis}.
\bjournal{Technometrics}
\bvolume{60}
\bpages{79--89}.
\end{barticle}
\endbibitem

\bibitem[\protect\citeauthoryear{Shen and Lin}{2019}]{shen2019statistical}
\begin{barticle}[author]
\bauthor{\bsnm{Shen},~\bfnm{Weijie}\binits{W.}} \AND
  \bauthor{\bsnm{Lin},~\bfnm{Dennis~KJ}\binits{D.~K.}}
(\byear{2019}).
\btitle{Statistical theories for dimensional analysis}.
\bjournal{Statistica Sinica}
\bvolume{29}
\bpages{527--550}.
\end{barticle}
\endbibitem

\bibitem[\protect\citeauthoryear{Shen et~al.}{2014}]{shen2014dimensional}
\begin{barticle}[author]
\bauthor{\bsnm{Shen},~\bfnm{Weijie}\binits{W.}},
  \bauthor{\bsnm{Davis},~\bfnm{Tim}\binits{T.}},
  \bauthor{\bsnm{Lin},~\bfnm{Dennis~KJ}\binits{D.~K.}} \AND
  \bauthor{\bsnm{Nachtsheim},~\bfnm{Christopher~J}\binits{C.~J.}}
(\byear{2014}).
\btitle{Dimensional analysis and its applications in statistics}.
\bjournal{Journal of Quality Technology}
\bvolume{46}
\bpages{185--198}.
\end{barticle}
\endbibitem

\bibitem[\protect\citeauthoryear{Stevens}{1946}]{stevens1946theory}
\begin{barticle}[author]
\bauthor{\bsnm{Stevens},~\bfnm{Stanley~Smith}\binits{S.~S.}}
(\byear{1946}).
\btitle{On the theory of scales of measurement}.
\bjournal{Science}
\bvolume{103}
\bpages{677-680}.
\end{barticle}
\endbibitem

\bibitem[\protect\citeauthoryear{Stevens}{1951}]{stevens1951mathematics}
\begin{barticle}[author]
\bauthor{\bsnm{Stevens},~\bfnm{Stanley~Smith}\binits{S.~S.}}
(\byear{1951}).
\btitle{Mathematics, measurement, and psychophysics}.
\bjournal{Handbook of Experimental Psychology}
\bpages{1-49}.
\end{barticle}
\endbibitem

\bibitem[\protect\citeauthoryear{Taylor}{2018}]{taylor2018quantity}
\begin{barticle}[author]
\bauthor{\bsnm{Taylor},~\bfnm{Barry~N}\binits{B.~N.}}
(\byear{2018}).
\btitle{Quantity calculus, fundamental constants, and SI units}.
\bjournal{Journal of Research of the National Institute of Standards and
  Technology}
\bvolume{123}
\bpages{123008}.
\end{barticle}
\endbibitem

\bibitem[\protect\citeauthoryear{Velleman and
  Wilkinson}{1993}]{velleman1993nominal}
\begin{barticle}[author]
\bauthor{\bsnm{Velleman},~\bfnm{Paul~F}\binits{P.~F.}} \AND
  \bauthor{\bsnm{Wilkinson},~\bfnm{Leland}\binits{L.}}
(\byear{1993}).
\btitle{Nominal, ordinal, interval, and ratio typologies are misleading}.
\bjournal{The American Statistician}
\bvolume{47}
\bpages{65--72}.
\end{barticle}
\endbibitem

\bibitem[\protect\citeauthoryear{Vignaux and Scott}{1999}]{vignaux1999theory}
\begin{barticle}[author]
\bauthor{\bsnm{Vignaux},~\bfnm{VA}\binits{V.}} \AND
  \bauthor{\bsnm{Scott},~\bfnm{JL}\binits{J.}}
(\byear{1999}).
\btitle{Theory \& methods: Simplifying regression models using dimensional
  analysis}.
\bjournal{Australian \& New Zealand Journal of Statistics}
\bvolume{41}
\bpages{31--41}.
\end{barticle}
\endbibitem

\bibitem[\protect\citeauthoryear{Ward}{2017}]{ward2017stevens}
\begin{barticle}[author]
\bauthor{\bsnm{Ward},~\bfnm{Lawrence~M}\binits{L.~M.}}
(\byear{2017}).
\btitle{{SS} {S}tevens's invariant legacy: scale types and the power law}.
\bjournal{American Journal of Psychology}
\bvolume{130}
\bpages{401--412}.
\end{barticle}
\endbibitem

\bibitem[\protect\citeauthoryear{Wijsman}{1967}]{wijsman1967cross}
\begin{binproceedings}[author]
\bauthor{\bsnm{Wijsman},~\bfnm{Robert~A}\binits{R.~A.}}
(\byear{1967}).
\btitle{Cross-sections of orbits and their application to densities of maximal
  invariants}.
In \bbooktitle{Proc. Fifth Berkeley Symp. on Math. Stat. and Prob}
\bvolume{1}
\bpages{389--400}.
\end{binproceedings}
\endbibitem

\bibitem[\protect\citeauthoryear{{Wikipedia}}{2020}]{wiki:tf}
\begin{bmisc}[author]
\bauthor{\bsnm{{Wikipedia}}}
(\byear{2020}).
\btitle{Transcendental function}.
\bhowpublished{\url{https://en.wikipedia.org/wiki/Transcendental_function}}.
\bnote{Accessed 2020/02/24}.
\end{bmisc}
\endbibitem

\bibitem[\protect\citeauthoryear{Wong and Zidek}{2018}]{wong2018dimensional}
\begin{barticle}[author]
\bauthor{\bsnm{Wong},~\bfnm{Samuel~WK}\binits{S.~W.}} \AND
  \bauthor{\bsnm{Zidek},~\bfnm{James~V}\binits{J.~V.}}
(\byear{2018}).
\btitle{Dimensional and statistical foundations for accumulated damage models}.
\bjournal{Wood science and technology}
\bvolume{52}
\bpages{45--65}.
\end{barticle}
\endbibitem

\bibitem[\protect\citeauthoryear{Yang and Lin}{2021}]{yang2021note}
\begin{barticle}[author]
\bauthor{\bsnm{Yang},~\bfnm{Ching-Chi}\binits{C.-C.}} \AND
  \bauthor{\bsnm{Lin},~\bfnm{Dennis~KJ}\binits{D.~K.}}
(\byear{2021}).
\btitle{A note on selection of basis quantities for dimensional analysis}.
\bjournal{Quality Engineering}
\bvolume{33}
\bpages{240--251}.
\end{barticle}
\endbibitem

\bibitem[\protect\citeauthoryear{Yang, Zidek and Wong}{2018}]{yang2018bayesian}
\begin{barticle}[author]
\bauthor{\bsnm{Yang},~\bfnm{Chun-Hao}\binits{C.-H.}},
  \bauthor{\bsnm{Zidek},~\bfnm{James~V}\binits{J.~V.}} \AND
  \bauthor{\bsnm{Wong},~\bfnm{Samuel~WK}\binits{S.~W.}}
(\byear{2018}).
\btitle{Bayesian analysis of accumulated damage models in lumber reliability}.
\bjournal{Technometrics}.
\end{barticle}
\endbibitem

\bibitem[\protect\citeauthoryear{Zhai et~al.}{2012a}]{zhai2012a}
\begin{btechreport}[author]
\bauthor{\bsnm{Zhai},~\bfnm{Yongliang}\binits{Y.}},
  \bauthor{\bsnm{Pirvu},~\bfnm{Ciprian}\binits{C.}},
  \bauthor{\bsnm{Heckman},~\bfnm{Nancy}\binits{N.}},
  \bauthor{\bsnm{Lum},~\bfnm{Conroy}\binits{C.}},
  \bauthor{\bsnm{Wu},~\bfnm{Lang}\binits{L.}} \AND
  \bauthor{\bsnm{Zidek},~\bfnm{James~V}\binits{J.~V.}}
(\byear{2012}a).
\btitle{A review of dynamic duration of load models for lumber strength}
\btype{Technical Report} No. \bnumber{270},
\bpublisher{Department of Statistics, University of British Columbia}.
\end{btechreport}
\endbibitem

\bibitem[\protect\citeauthoryear{Zhai et~al.}{2012b}]{zhai2012b}
\begin{btechreport}[author]
\bauthor{\bsnm{Zhai},~\bfnm{Yongliang}\binits{Y.}},
  \bauthor{\bsnm{Heckman},~\bfnm{Nancy}\binits{N.}},
  \bauthor{\bsnm{Lum},~\bfnm{Conroy}\binits{C.}},
  \bauthor{\bsnm{Pirvu},~\bfnm{Ciprian}\binits{C.}},
  \bauthor{\bsnm{Wu},~\bfnm{lang}\binits{l.}} \AND
  \bauthor{\bsnm{Zidek},~\bfnm{James~V}\binits{J.~V.}}
(\byear{2012}b).
\btitle{Stochastic models for the effects of duration of load on lumber
  properties}
\btype{Technical Report} No. \bnumber{271},
\bpublisher{Department of Statistics, University of British Columbia}.
\end{btechreport}
\endbibitem

\bibitem[\protect\citeauthoryear{Zidek}{1969}]{zidek1969representation}
\begin{barticle}[author]
\bauthor{\bsnm{Zidek},~\bfnm{James~V}\binits{J.~V.}}
(\byear{1969}).
\btitle{A representation of Bayes invariant procedures in terms of Haar
  measure}.
\bjournal{Annals of the Institute of Statistical Mathematics}
\bvolume{21}
\bpages{291--308}.
\end{barticle}
\endbibitem

\end{thebibliography}


\end{document}